\documentclass{article}
\usepackage{graphicx} 
\usepackage{tikz}
\usepackage{tikz-cd}
\usepackage{quiver}
\usepackage{amssymb}
\usepackage{amsmath}
\usepackage{amsthm}
\usepackage{mathtools}
\usepackage{bm}
\usepackage{framed,enumitem} 
\usepackage{hyperref}
\usepackage{subfiles}
\usepackage{bbm}
\usepackage[noadjust]{cite} 
\usepackage{varwidth}
\usepackage{geometry}
\hypersetup{colorlinks=true,citecolor=blue,
linkcolor=blue,urlcolor=black,filecolor=black}

\newtheorem{theorem}{Theorem}
\numberwithin{theorem}{section}

\theoremstyle{definition}
\newtheorem{definition}[theorem]{Definition}

\newtheorem{proposition}[theorem]{Proposition}

\newtheorem{lemma}[theorem]{Lemma}

\newtheorem{remark}[theorem]{Remark}

\newtheorem{example}[theorem]{Example}

\newtheorem{corollary}[theorem]{Corollary}

\newtheorem*{notation}{Notation}

\newcommand{\N}{\mathbb{N}}
\newcommand{\Z}{\mathbb{Z}}

\newcommand{\C}{\mathbb{C}}
\newcommand{\K}{\mathbb{K}}
\newcommand{\1}{\mathbf{1}}
\newcommand{\Id}{\mathrm{Id}}
\newcommand{\rad}{\mathrm{rad}\;}

\newcommand{\matfont}[1]{\pmb{\mathsf{#1}}}
\newcommand{\mb}[1]{\mathbf{#1}}
\newcommand{\Eta}{\mathrm{H}}

\numberwithin{theorem}{section}
\newtheorem*{theoremA*}{\textbf{Main Theorem A}}

\newtheorem*{theoremC*}{\textbf{Main Theorem C}}

\newtheorem*{conjectureB*}{\textbf{Conjecture B}}

\newcommand{\Addresses}{{% additional braces for segregating \footnotesize
  \bigskip
  \footnotesize

  \textsc{Institute of Mathematics, University of Zurich, Winterthurerstrasse 190, CH-8057 Zurich, Switzerland}\par\nopagebreak
  \textit{E-mail address}: \texttt{maksymilian.manko@math.uzh.ch}

}}

\title{\normalsize \textbf{REPRESENTATION THEORY OF NON-FACTORIZABLE RIBBON HOPF ALGEBRAS}}

\author{Maksymilian Manko}
\date{\normalsize{March 2026}}

\begin{document}

\maketitle

\begin{abstract}
    In \cite{faes2025non} new examples of ribbon Hopf algebras based on the construction due to Nenciu \cite{nenciu_2004} were presented. This piece serves as a sequel where we study the representation theory of these new examples of ribbon Hopf algebras. We classify indecomposable projective and simple modules, find the Krull-Schmidt decomposition of the adjoint representation of Nenciu algebras, and prove fusion rules between its components. We also comment on the properties of M\"uger centres of their representation categories, in particular that they can be non-semisimple. Finally, we consider a new family of ribbon Hopf algebras over fields of prime characteristic $p>2$ in the context of 4-dimensional TQFTs presented in \cite{costantino2023skein} that constitute an improvement over examples given therein, although still seemingly falling short of producing powerful invariants of 4-manifolds.
\end{abstract}
% \keywords{Representation theory of Hopf algebras, Ribbon categories, TQFTs,  4-manifolds\\
% \textbf{Mathematics Subject Classification (2020)}: 16T05, 16G99, 18M15, 57K41, 57K16}

\tableofcontents
\section{Introduction}
Unimodular ribbon categories introduced by Turaev in \cite{turaev_94} constitute an important topic in quantum algebra, chiefly for their many applications in quantum topology. Since the 1990s, multiple constructions of 3-manifold and 4-manifold invariants have been studied using these as essential ingredients \cite{turaev_94}, \cite{kerler2001non}, \cite{crane1993categorical}. Ribbon Hopf algebras are a natural source of such categories through their representation theory. A ribbon Hopf algebra $H$ is a quasitriangular Hopf algebra (i.e. it comes with the extra data of an invertible element $R$ of $H\otimes H$ called the \textit{R-matrix}) together with a central element of $H$ called the ribbon element. Both $R$ and $v$ are required to satisfy certain axioms, see Section \ref{secprelims}.

In particular, one may define a $3$-manifold invariant from any so-called \emph{factorizable} unimodular ribbon Hopf algebra. Factorizability is a non-degeneracy condition on the \emph{monodromy matrix} $M:= R_{21}R$ of the algebra, where $R_{21}$ is the R-matrix with permuted entries (see Proposition \ref{factorizability}). Recall also that a Hopf algebra $H$ is unimodular if it admits an integral $\lambda \in H^*$ and a cointegral $\Lambda \in H$ satisfying certain properties, as detailed in Definition \ref{defunimod}.

More recently there has been a considerable interest in the unimodular, \textit{non-factorizable} ribbon Hopf algebras, as ingredients of non-semisimple 4-manifold invariants. In \cite{beliakova_derenzi_2023} Beliakova and De Renzi introduced Kerler-Lyubashenko functors from a certain category of 4-dimensional 2-handlebodies up to 2-deformations \cite{beliakova_derenzi_2023} based on the work of Bobtcheva and Piergallini \cite{bobtcheva_piergallini_2006}, and presented in a more streamlined way in a pre-print of Beliakova, Bobtcheva, De Renzi and Piergallini \cite{beliakova_bobtcheva_derenzi_piergallini_2023}. 

Several examples of non-factorizable unimodular ribbon Hopf algebras are known classically, including the small quantum group $u_q \mathfrak{sl}_2$ at an even root of unity $q$. It was shown in \cite{beliakova2022refined} that the scalar Kerler-Lyubashenko type invariant of \cite{beliakova_derenzi_2023} applied to a connected $4$-dimensional $2$-handlebody only sees its $3$-dimensional boundary and homological information. Another example is the Hopf algebra of \textit{symplectic fermions} $\operatorname{SF}_{2n}$, which admits a continuous family of non-factorizable ribbon structures. The generic case, for $n=1$, was shown by Kerler in \cite{kerler2003homology} to give essentially invariants of the $3$-dimensional boundary with this approach. In both cases, the Hopf algebras are, in some sense, ``close to factorizable". Indeed, in the unimodular case, factorizability is equivalent to requiring that applying the integral to one side of the monodromy matrix gives the cointegral. In these non-factorizable cases we mentioned, this equation is still verified up to an invertible element of the algebra.

Hence, in an attempt to avoid such ``degeneration to the boundary" in our recent paper with Q. Faes \cite{faes2025non}, henceforth referred to as part 1, we introduced the stricter notion of \textit{strong non-factorizability}. We demand that the result of applying the integral to either side of the monodromy matrix, with the antipode previously applied to its first tensor component, is zero, and hence cannot be related to the cointegral by any invertible element of the algebra:
$$\big((\lambda \circ S)\otimes \Id\big) M = \big( S \otimes \lambda \big)M= 0.$$ 
We introduced two extensive families of such Hopf algebras, we will refer to as \textit{Nenciu type} after \cite{nenciu_2004} and \textit{Nenciu biproducts}, respectively, as the second family is produced by taking semidirect products and coproducts between the Nenciu type Hopf algebras and the small quantum group $u_q \mathfrak{sl}_2$ at the unit of unity $q$ of order divisible by 4. 

In his paper we focus on the related representation-theoretic questions. From this perspective, module categories of such Hopf algebras are analogously unimodular, non-factorizable ribbon, and by the result of Shimizu \cite{shimizu_2016}, the factorizability condition can be equivalently expressed in terms of the \textit{M\"uger centre} - the full subcategory of modules which braid symmetrically with any other module in the category, see Definition \ref{Muger}. Precious few examples of such categories are known, even fewer if we additionally impose the condition that the M\"uger centre be non-semisimple (a notable example are \textit{higher Verlinde categories in the mixed case} defined in \cite{sutton2023sl} and studied in \cite{decoppet2024higher}). We show that the examples of Hopf algebras discussed in part 1 posses this property.
Non-factorizable braided tensor categories have been also recently studied in the context of Morita 4-category of braided tensor categories in \cite{Brochier_Jordan_Snyder_2021} and 4-category of braided, $\mathcal{E}$-enriched tensor categories in \cite{decoppet2025relative}, both works referencing the algebraic side of the topic. As also mentioned in part 1, such categories have applications in 4-dimensional topology, as the underlying algebraic data of the 4-dimensional Topological Quantum Field Theories (TQFTs) constructed by Constantino, Geer, Ha\"ioun and Patureau-Mirand in \cite{costantino2023skein}, which we will discuss towards the end of the paper. 
%\textcolor{red}{The need for such examples is further underlined by the results of Reutter \cite{Reutter_23} that semisimple 4d TQFTs cannot distinguish exotic structures, and thus non-semisimple counterparts are necessary for his task.} 
Additionally, in this case, it is expected that for the TQFT to detect 4-dimensional information the category has to be \textit{chromatic-non degenerate} but not \textit{chromatic compact} (see Definition \ref{ChromNonDegenCompact}). For representation categories, this is a property linking the action of the two-sided cointegral $\Lambda$ and of the element $\lambda(S(M'))M''$, where $\lambda$ is the left integral, on the projective cover of the monoidal unit (see Definition \ref{defprojcover}). Topologically this carries the information about handle attachments of index greater than 2 in the numbering convention of \cite{beliakova_derenzi_2023} we used in part 1 (but dual to the one in \cite{costantino2023skein}). We show that while this is not possible to achieve with the constructions of part 1, passing to a field of finite characteristic allows to produce such a category with our methods. This is an improvement on the strategy of \cite[Definition 6.6 and Proposition 6.7]{costantino2023skein}, producing, to the author's knowledge the only example of braided, but not symmetric unimodular ribbon category with non-semisimple M\"uger centre that is chromatic non-degenerate but not chromatic compact. We comment also on its potential power to detect exotic pairs of 4-manifolds.

\subsection{Main results}
The first result is the decompostion of the regular representation of any Nenciu-type Hopf algebra into indecomposable projective modules (Theorem \ref{regDecomp}). This is achieved by explicitly finding the full system of primitive idempotents. The idempotents are in one-to-one correspondence with the grouplike elements and we collect the eigenvalues of the action by grouplikes into tuples that identify the idempotents uniquely, in line with the conventions of part 1. We then do the same with the adjoint representation (Theorem \ref{adjointDecomp}), which beyond the indecomposable projective and corresponding simple modules contains also direct summands we identify with submodules of the indecomposable projectives. We proceed to study their fusion rules (Theorem \ref{fusionRulesNenciu}) and brading properties, in particular identifying when such modules are transparent (Proposition \ref{symBraiding}). We emphasize that most of the properties of this family of Hopf algebra is controlled by the data of the grouplike generators and their actions on the nilpotent generators. Finally, we consider further indecomposable modules deduced from the results of \cite{zaitsev_nikolenko_1971} (Example \ref{XiandEtamodules}). \\

Moving on to the case of Nenciu biproducts, we again find the full system of idempotents (Proposition \ref{productFullSystem}) and find an explicit basis for the corresponding projective modules (Definition \ref{productProjective}). We find that many properties of the modules over the biproducts are composite with regard to those of the constituents, as emphasized in Examples \ref{prodExample1} and $\ref{prodExample2}$. This includes the properties of the M\"uger centre, chiefly its non-semisimplicity (Propositions \ref{NenicuNonSSMuger} and \ref{extInfiniteMuger}). The results for both families are supported with plethora of examples.\\

In the final section we show that the representation categories of strongly non-factorizable Hopf algebras are neither chromatic non-degenerate nor chromatic compact (Proposition \ref{snfChromNonCompact}), so do not produce 4-dimensional TQFTs defined beyond 0-, 1- and 2-handles from the approach of \cite{costantino2023skein}. This holds, in particular, for strongly non-factorizable examples discussed in part 1. However, it is possible to modify our construction to the case of fields of finite characteristic to achieve chromatic non-degeneracy, but not chromatic compactness with our means (Proposition \ref{chromNon-Degen}), while retaining the non-semisimple M\"uger centre (Proposition \ref{UHGmuger}). Namely, we fix a field $\K$ of finite characteristic $\operatorname{char} \K = p>2$, such that it contains a root of unity $q\in \K$ and the order of $q$ is divisible by 4 and coprime to $p$. Let $S_p$ be the symmetric group on $p$ letters and $\K[S_p]$ its group algebra. Let also $\operatorname{SF_{2p}}$ be the Hopf algebra of symplectic fermions over $\K$ and $u_q \mathfrak{sl_2}$ the small quantum group at the root of unity $q$ over $\K$ (we explain why all necessary results regarding these hold in finite characteristic in Propositions \ref{nenciu_finitefield} and \ref{biproduct_finitefield}). Then we construct the following composite Hopf algebra.
\begin{theoremA*}
    Under an appropriate choice of parameters a semidirect biproduct $u_q \mathfrak{sl}_2 \ltimes \operatorname{SF}_{2p} \rtimes \K[S_p]$ can be simultaneously
    \begin{enumerate}[leftmargin=1cm]
        \item unimodular (Theorem \ref{UHGProp}, 1.),
        \item quasitriangular (Theorem \ref{UHGProp} 2.)
        \item ribbon (Theorem \ref{UHGProp} 3.),
        \item twist degenerate (Theorem \ref{UHGProp} 4.),
        \item chromatic non-degenerate but not chromatic compact (Proposition \ref{chromNon-Degen}).
    \end{enumerate}
\end{theoremA*}

However, we also show this is not enough to detect exotic pairs of 4-manifolds as the invariant will be susceptible to stabilisation by $S^2 \times S^2 $ (Propositon \ref{StabProblem}), according to the known result of Gompf (see Remark \ref{StabRem}).

  \subsection{Structure of the paper}
    In Section 2 we introduce algebraic and representation-theoretic preliminaries, while in Section 3 and 4 we recall some of the needed material from \cite{faes2025non}, regarding Nenicu Hopf algebras and biproducts with $u_q \mathfrak{sl}_2$, respectively. In Section 5 we study representation theory of Nencu algebras, including the classification of simple and projective modules, decomposition of the adjoint representation into indecomposable modules, relevant fusion rules and braiding properties. In Section 6 we do the same with the Nenciu biproducts with $u_q \mathfrak{sl}_2$, bar the decomposition of the adjoint representation. Finally, in Section 7 we consider an example in finite characteristic and put it in the context of the 4-dimensional TQFTs of \cite{costantino2023skein}.
\\\\
\noindent\textbf{Acknowledgements:} 
The author would like to thank Ivelina Bobtcheva, Marco de Renzi,  Azat Gainutdinov, Cris Negron, and Benjamin Ha\"ioun for fruitful discussions, and Quentin Faes, Matthieu Faitg and Thibault D\'ecoppet for helpful comments on the draft, as well as Riccardo Piergallini for making his Mathematica computer program available, and Anna Beliakova for her help and introducing him to the topic. The author was supported by Simons Collaboration on New Structures in Low-Dimensional Topology
and Grant 200020\_207374 of the Swiss National Science Foundation.

\section{Preliminaries}
\label{secprelims}
Let us recall some definitions and results regarding Hopf algebras, and fix the notation used in the rest of the paper.
\subsection{Notations for Hopf algebras}
In this paper, we consider Hopf algebras over the field of complex numbers $\C$, or a field of finite characteristic $p$, $\K$. We usually assume finite dimension, which implies, in particular, that the antipode is bijective \cite[Theorem 7.1.14 (b)]{radford_2012}.
\begin{definition}
\label{Hopf_alg_def}
    A \textit{Hopf algebra} $(H, \mu, \eta, \Delta,\epsilon, S)$ is the data of a vector space (over $\C$ or $\K$) $H$ together with a \textit{unit} $\1: \C \rightarrow H$ (resp. $\K$) and a \textit{product} $\mu: H\otimes H \rightarrow H,$  a \textit{counit} $\epsilon: H \rightarrow \C$ (resp. $\K$) and a \textit{coproduct} $\Delta: H \rightarrow H\otimes H$ as well as an invertible \textit{antipode} $S : H\rightarrow H$,
    that are $\C$-linear (resp. $K$-linear) maps satisfying certain axioms, see for instance \cite{radford_2012}.\\
% $$
% \begin{gathered}
% \mu \circ(\mu \otimes \mathrm{id})=\mu \circ(\mathrm{id} \otimes \mu), \\
% \mu \circ(\eta \otimes \mathrm{id})=\mathrm{id}=\mu \circ(\mathrm{id} \otimes \eta), \\
% (\Delta \otimes \mathrm{id}) \circ \Delta=(\mathrm{id} \otimes \Delta) \circ \Delta, \\
% (\varepsilon \otimes \mathrm{id}) \circ \Delta=\mathrm{id}=(\mathrm{id} \otimes \varepsilon) \circ \Delta, \\
% (\mu \otimes \mu) \circ(\mathrm{id} \otimes \tau \otimes \mathrm{id}) \circ(\Delta \otimes \Delta)=\Delta \circ \mu, \\
% \varepsilon \circ \mu=\varepsilon \otimes \varepsilon, \\
% \Delta \circ \eta=\eta \otimes \eta, \\
% \varepsilon \circ \eta=1, \\
% \mu \circ(S \otimes \mathrm{id}) \circ \Delta=\eta \circ \varepsilon=\mu \circ(\mathrm{id} \otimes S) \circ \Delta, 
% %\\S \circ S^{-1}=\mathrm{id}=S^{-1} \circ S.
% \end{gathered}$$
We denote $\tau:H\otimes H \rightarrow H \otimes H$ the \textit{flip map} exchanging the factors in the tensor product.
We will often abusively denote the Hopf algebra by $H$, and write $gh$ instead of $\mu(g,h)$ for $g,h \in H$. We also note that $H\otimes H$ can be endowed with an algebra structure, using the product
$ (a\otimes b)(c\otimes d):= ac \otimes bd$, where $a,b,c,d \in H$.
\end{definition}
\begin{remark}
    Hopf algebras may be defined in any braided monoidal category $\mathcal{C}$. The algebra structure of $H\otimes H$ is then defined using the braiding.
\end{remark}

\begin{notation}
     We use the Sweedler notation for the coproduct: if $\Delta(h) = \sum_{i =1}^n h_{(1)i} \otimes h_{(2)i} $ then we write
    $        \Delta(h) = h_{(1)} \otimes h_{(2)},$
    for some $n\in \Z_{>0}$.
By $\Delta^{cop}:=\tau \circ \Delta$ we denote the \textit{coopposite coproduct}, with $\tau:H\otimes H \rightarrow H \otimes H$ the flip map. Recall that if $H$ is a Hopf algebra, then $H^{cop}:= (H,\mu,\eta,\Delta^{cop}, \epsilon, S^{-1})$ is also a Hopf algebra, provided $S$ is bijective.
\end{notation}
Recall that if $H$ is a Hopf algebra, the dual vector space
        $H^*:= \operatorname{Hom}_\C (H, \C)$
    can also be provided with a Hopf algebra structure.

\begin{definition}
     An element $g\in H$ is called \textit{grouplike} if $\Delta(g)= g\otimes g$. The subset of all such elements is denoted $G(H)$. An element $X \in H$ is called \textit{skew-primitive with respect to} $L\in G(H)$ if $\Delta(X) = \1 \otimes X + X\otimes L$, and \textit{primitive} if $L=\1$. 
\end{definition}
\noindent The set $G(H)$ forms a group where the inverse of $g\in G(H)$ is given by $g^{-1} := S(g)$. It also generates a Hopf subalgebra of $H$, which is isomorphic to the group algebra of $G(H)$.
\begin{definition}
\label{adjointaction}
    The \emph{left adjoint action} of $H$ is the left $H$-module structure on $H$ given by
    $\rhd: H\otimes H \rightarrow H, \;\; a\rhd b :=  a_{(1)}\,b \,S(a_{(2)})$
    in the Sweedler notation.
\end{definition}

\subsection{Integrals and unimodularity}

Let $H$ be a finite dimensional Hopf algebra over an algebraically closed field. Then it admits a left (resp. right) integral $\lambda_L \in H^*$ (resp. $\lambda_R$) and cointegral $\Lambda_L\in H$ (resp. $\Lambda_R$), see \cite[Chapter 10.1]{radford_2012}. 
\begin{definition}
\label{defunimod}
    An element $\lambda_L \in H^*$ is a \textit{left} (resp. $\lambda_R \in H^*$ a \textit{right}) \textit{integral} of $H$ if for all $h\in H$,
\begin{align*}
    %\label{integralDef}
    \lambda_L(h_{(2)}) h_{(1)} = \lambda_L(h) \1 && \big(\text{resp.} && \lambda_R(h_{(1)}) h_{(2)}  = \lambda_R(h)\1 \big).
\end{align*} 
\noindent An element $\Lambda_L \in H$ is a left (resp. $\Lambda_R \in H$ a right) cointegral if for all $h\in H$, 
\begin{align*}
    %\label{cointegralDef}
     h\Lambda_L = \epsilon(h)\Lambda_L && \big(\text{resp.} && \Lambda_R h = \epsilon(h) \Lambda_R\big).
\end{align*}
An integral (resp. cointegral) is called two-sided if its both a left and right integral (resp. cointegral), or equivalently, if
\begin{align*}
    \lambda \circ S = \lambda && \big(\text{resp.} && S(\Lambda) = \Lambda\big).
\end{align*}
\end{definition}
\noindent It is also known that $\Lambda_L, \Lambda_R$ generate $1$-dimensional subalgebras of $H$ and $\lambda_L, \lambda_R$ generate $1$-dimensional subalgebras of $H^*$, thus different choices differ up to a constant in the ground field. We also note that in case $\lambda_L$ and $\lambda_R$ do not coincide, it still holds that $\lambda_R = \lambda_L \circ S$.
\begin{definition}[\cite{radford_2012} Definition 10.2.3]
    We will call $H$ \textit{unimodular} if $\Lambda$ is two-sided.
\end{definition}

\subsection{Quasitriangular and ribbon structures}
Let us recall the definition of the quasitriangular Hopf algebra. We need the following notation.

\begin{notation}
    For an element $R\in H\otimes H$ such as the \textit{R-matrix} of a quasitriangular Hopf algebra of Definition \ref{QTdef} below, we write $R = R' \otimes R'',$
where the summation is again implicit. Sometimes, we will need many copies of $R$, in which case they will be labeled with numbers. We also denote $R_{21}:= \tau(R)$.
\end{notation}

\begin{definition}
\label{QTdef}
    Let $H$ be a Hopf algebra and $R$ be an element of $H\otimes H$. The pair $(H,R)$ is quasitriangular if the following axioms are satisfied.
    \begin{center}
    \begin{varwidth}{\textwidth}
        \begin{itemize}[leftmargin=1cm]
            \item[(QT1)] $\Delta(R') \otimes R'' = R'_1 \otimes R'_2 \otimes R''_1 R''_2$
            \item[(QT2)] $\epsilon(R')R'' = \1$
            \item[(QT3)] $R' \otimes \Delta^{cop} (R'') = R'_1 R'_2 \otimes R''_1 \otimes R''_2$
            \item[(QT4)] $\epsilon(R'')R' = \1$
            \item[(QT5)] $\Delta^{cop} (h) R = R\Delta(h), \forall h \in H$.
        \end{itemize}
    \end{varwidth}
    \end{center}
    \noindent Then $R$ is called a  \textit{universal R-matrix} (where we usually drop the ``universal"), and it follows from (QT2) and (QT4) that it is invertible in the algebra $H\otimes H$, with the inverse denoted $R^{-1}$.
\end{definition} 
\noindent We also note that for any quasitriangular Hopf algebra $(H, R)$, the antipode $S$ is bijective and the equality $(S\otimes S)(R) = R$ holds.
\begin{definition}
    The \textit{monodromy matrix} is defined as
        $M:= R_{21} R = R''_1 R'_2 \otimes R'_1 R''_2.$ If $M =  \1\otimes \1$, the quasitriangular structure is called \textit{triangular}.
\end{definition}
\noindent 
In order to exhibit quasitriangular structures on Hopf algebras, we will strongly rely on the following result.
\begin{theorem}[\cite{radford_1994}, Section 2.1] \label{QTthm} 
\label{classicalTheorem}
    Let $H$ be a Hopf algebra. Any element $R \in H \otimes H$ satisfying (QT1)-(QT4) induces a Hopf algebra map 
         $f_R: H^* \rightarrow H^{cop}, h^* \mapsto h^*(R')R'',$
    where $H^{cop}$ is the \textit{coopposite Hopf algebra} of $H$, defined above.
\end{theorem}
\noindent By the contrapositive of Theorem \ref{classicalTheorem}, if there are no bialgebra maps between $H^*$ and $H^{cop}$, then there can be no quasitriangular structures on $H$.

Let $\omega: H \otimes H \mapsto \C$, be a billinear pairing. We call $\omega$ \textit{non-degenerate} if any element $h \in H$ such that $\omega(h, g) = 0$ for all $g\in H$, satisfies $h=0$.
\begin{proposition}
[\cite{shimizu_2016} Theorem 1.1, \cite{etingof_gelaki_nikshych_ostrik_2016} Exercise 8.6.4] \label{factorizability}
    For a finite-dimensional quasitriangular Hopf algebra $(H, R)$ the following two conditions are equivalent
    \begin{enumerate}[leftmargin=1cm]
        \item the monodromy matrix $M$ induces a non-degenerate pairing on $H^*\otimes H^* \rightarrow \C$ given by $h^*_1 \otimes h^*_2 \mapsto h^*_1 (M') h^*_2(M'')$, or
        \item the Drinfeld-Reshetikhin map 
             $H^* \rightarrow H, h^* \mapsto h^*(M')M''$
        is an isomorphism of vector spaces.
        \item The M\"uger centre $\mathcal{Z}_{(2)}(H\text{-mod})$ (see Definition \ref{Muger}) of the category of modules $H\text{-mod}$ is trivial (generated by the unit of $H\text{-mod}$). 
    \end{enumerate}
    If either item is fulfilled $(H, R)$ is called \textit{factorizable}, and otherwise \textit{non-factorizable}. If $\mathcal{Z}_{(2)}(H\text{-mod})$ is isomorphic to the category of finite dimensional super-vector spaces  $\operatorname{SVect}$, we will refer to $H$ as \textit{slightly factorizable}.
\end{proposition}
\noindent For the purpose the discussion, in \cite{faes2025non} we introduced a stronger notion
\begin{definition}
\label{SNFdef}
    A unimodular, quasitriangular Hopf algebra $(H, R)$ with a left integral $\lambda_L \in H^*$ is called \textit{strongly non-factorizable} if any of the two equivalent conditions is fulfilled
        \begin{align*}
         \lambda_L(S(M')) M'' = 0 && \text{or} && \lambda_L(M'') S(M')  = 0.
    \end{align*}
\end{definition}
\noindent It is easy to see from Definition \ref{factorizability} that any strongly non-factorizable Hopf algebra is non-factorizable.
\begin{remark}
    Let $(H, R)$ be unimodular triangular, so $M = \1 \otimes \1$, and let it have a left integral $\lambda_L$. Then, by \cite[Theorem 10.3.2]{radford_2012} it is cosemisimple if and only if $\lambda_L(\1) = 1$. Thus, if $H$ is non-cosemisimple, then $(H, R)$ is strongly non-factorizable.
\end{remark}

\noindent We now turn to the ribbon structure.

\begin{definition}
\label{ribbonAxioms} Let $(H, R)$ be a finite dimensional, quasitriangular Hopf algebra with monodromy matrix $M$. A \textit{ribbon element} $v\in Z(H)$ is a central element of $H$ such that
\begin{center}
\begin{varwidth}{\textwidth}
\begin{itemize}[leftmargin=1cm]
    \item[(R1)] $S(v) = v$
    \item[(R2)] $\epsilon(v)=1$
    \item[(R3)] $M\Delta(v) = v\otimes v$.
\end{itemize}
\end{varwidth}
\end{center}
The triple $(H, R, v)$ is called a \textit{ribbon Hopf algebra}.
\end{definition}

\begin{definition}
\label{drinfeldelem}
    The \textit{Drinfeld element} of a quasitriangular Hopf algebra $(H, R)$ is $u:= S(R'') R'.$
\end{definition} 

\noindent Moreover, $u$ is invertible with the inverse given by $u^{-1} = R'' S^2(R').$
\begin{proposition}[\cite{kassel_2012}, VIII.4.1]\label{DrinfeldProp1} Let $(H, R)$ be a quasitriangular Hopf algebra and $u$ be the corresponding Drinfeld element. For any $h \in H$ we have $ S^2(h) = u h u^{-1}.$
\end{proposition}
\noindent We call a Hopf algebra automorphism $f \in \operatorname{Aut}(A)$ \textit{inner} if there exists a fixed invertible element $u_f \in A$ such that for any $a \in A$, $f(a) = u_f a u^{-1}_f$ holds. It is easy to see that $f^{-1}$ is also inner with the fixed invertible element $u^{-1}_f$. In this sense $S^2$ is inner, which implies, in particular, that $S$ is bijective.
\begin{definition}[\cite{radford_1994}]
    Let $H$ be a Hopf algebra. A \textit{pivotal element} is a grouplike element $g\in \mathcal{G}(H)$ such that for any $h\in H$
    \[
        S^2(h) = g h g^{-1}.
    \]
\end{definition}
\noindent As explained in \cite[Prop. 2]{radford_1994} in a ribbon Hopf algebra, a ribbon element $v=g^{-1}u$ is the product of the Drinfeld element $u$ and the inverse of a pivotal element. However, not every pivotal element can be used to define a ribbon element. 
%\begin{definition}
 %   Let $H$ be  Hopf algebra admitting a left integral $\lambda_L \in H^*$. The the \textit{distinguished grouplike element} $a\in \mathcal{G}(H)$ is the unique element such that $\lambda_L \circ S(h) = \lambda_L(ah)$. 
%\end{definition}
%\begin{lemma}[\cite{radford_1994}, Proposition 2.(b)] \label{2SidedCoIntLemm} Let $(H, R)$ be a finite-dimensional quasitriangular Hopf algebra with a two-sided cointegral and denote the Drinfeld element by $u$. Let $g \in \mathcal{G}(H)$ be a pivotal element and $a$ be the distinguished grouplike element. Then $v:=g^{-1}u$ is a ribbon element if and only if 
%\[
 %   g^2 = a^{-1}.
%\]
%\end{lemma}

\begin{lemma}
\label{2SidedIntLemm}
    Let $(H, R)$ be a finite-dimensional quasitriangular Hopf algebra with a two-sided integral $\lambda$. If $g\in \mathcal{G}(H)$ is a pivotal element, then $g^2=\1$.
\end{lemma}
% \begin{proof}
%     If $\lambda$ is two-sided, then the distinguished grouplike element $a$ of Lemma \ref{2SidedCoIntLemm} satisfies $a=\1$. Thus, $g^2 = \1$.
% \end{proof}
\noindent We introduce one more property of the ribbon structure that bears significance in many topological constructions.
\begin{definition}
\label{DefAF}
    Let $H$ be a Hopf algebra carrying a (left) integral $\lambda$ and a ribbon element $v$. We call $H$ \textit{anomaly-free} (or \textit{twist non-degenerate} in \cite{costantino2023skein}) if  $\lambda(v) \neq 0.$
    Otherwise we will refer to it as \textit{anomalous} (or \textit{twist degenerate}).
\end{definition}
\subsection{Algebra generators and tuple notation}
All the constructions considered in this paper produce Hopf algebras generated by grouplike and skew-primitive generators. 
\begin{notation}
    We denote by $K_a$, $a=1, \dots, s$ the grouplike generators that commute with one another and are of respective orders $m_1, \dots, m_s$. By $X_k, X^\pm_k, Z_l^\pm$ for $k=1, \dots, t_1$ and $l=1, \dots, t_2$ we denote skew-primitive generators such that $(X_k)^2 = (X^\pm_k)^2 = (Z^\pm_l)^2=0$.
\end{notation}
\noindent We will repeatedly appeal to the following notion.
\begin{definition}
    Let $g, h\in H$. We say $g, h$ have a \textit{diagonal relation} if there exists $\gamma \in \C$ such that $gh = \gamma hg.$
\end{definition}
\begin{notation}
     We will use boldface for tuples of numbers and algebra generators, for instance, let $\mb{0}=(0, \dots, 0)$, and let
     \[
        \pmb{\xi} = (\xi_1, \dots, \xi_s):= (\exp^{2\pi/m_1}, \dots, \exp^{2\pi/m_s}) \in \C^s
     \]
     be the tuple of primitive roots of unity of orders $\mb{m} = (m_1, \dots, m_s) \in \Z^s_{>0}$. We similarly express the lists of algebra generators as $\mb{K}=(K_1, \dots, K_s)$, or $\mb{X} = (X_1, \dots, X_t)$.
\end{notation}
\noindent We will often use the following algebraic structures.
\begin{notation}
    Let $\Z_\mb{m} := \Z_{m_1}\times\dots\times \Z_{m_s}$, where $\Z_{m_a}:=\Z/m_a\Z$, and let $\C[\mb{K}] := \C[K_1, \dots, K_{m_s}] \subset H$ be the group subalgebra generated by $K_1, \dots, K_s$. It is isomorphic to $\C[\Z_\mb{m}]$ as a Hopf algebra.
\end{notation}
\noindent We also introduce some operations on tuples.
\begin{notation}
    For two tuples of the same shape, for instance $\mb{w}, \mb{v} \in \Z_\mb{m}$ we define the element-wise sum and product 
    \begin{align*}
         \mb{w} + \mb{v}:= (w_1+v_1, \dots, w_s+v_s) &&  \mb{w} \cdot \mb{v}:=(w_1 v_1, \dots, w_s v_s).
    \end{align*}
    By a power of a tuple by a tuple, for instance for $\mb{v} \in \Z_{\mb{m}}$, $\pmb{\xi}\in \C^s$ and $\mb{K}$, we mean
    \begin{align*}
         \pmb{\xi}^\mb{v} = (\xi_1)^{v_1}(\xi_1)^{v_2}\dots(\xi_s)^{v_s} && \text{and} && \mb{K}^\mb{v} = (K_1)^{v_1}(K_2)^{v_2}\dots(K_s)^{v_s}.
    \end{align*}
\end{notation}
\begin{notation}
    We will also encounter \textit{matrices}, denoted with boldface sans-serif font, for instance 
    \begin{align*}
        \pmb{\alpha} = (\alpha_{kl})_{1 \leq k, l \leq t} \in \operatorname{Mat}_{t\times t}(\C), && \matfont{d}  = (d_{ka})_{\substack{1 \leq a \leq s \\ 1 \leq k \leq t}}, && \matfont{u} = (u_{ka})_{\substack{1 \leq a \leq s \\ 1 \leq k \leq t}} \in \operatorname{Mat}_{t\times s}(\Z_\mb{m}).
    \end{align*}
     The rows $\matfont{d}_k = (d_{ka})_{0\leq a \leq s}$, $\matfont{u}_k=(u_{ka})_{0\leq a \leq s}$ are considered as tuples (and conform to the boldface notation). Products with no operation symbols mean the usual matrix multiplication, provided the participating matrices have appropriate dimensions and entries in the same ring. For $\mb{w}, \mb{v} \in \Z_\mb{m}$,  
\[
  \mb{w} \mb{v}^T = w_1 v_1+\dots+w_s v_s,
\]
is the inner product of row tuples, valued in $\Z_{\mathrm{lcm}(m_1, \dots, m_s)}$, where $\mathrm{lcm}(m_1, \dots, m_s)$ is the least common multiple of $m_1, \dots, m_s$. For a matrix $\matfont{z}\in \operatorname{Mat}_{s\times s}(\Z_\mb{m})$, the expression $\mb{v} \matfont{z}$ is a row tuple again, as is $(\mb{w}\mb{z})\cdot \mb{v}$, but $\mb{w}\matfont{z}\mb{v}^T$ is a scalar in $\Z_{\mathrm{lcm}(m_1, \dots, m_s)}$. 
\end{notation}

\subsection{The semi-direct biproduct structure}
Towards the goal of introducing Nenciu biproducts later on, we introduce versions of \textit{semi-direct} (or \textit{smash}) product and coproduct. While these constructions are well-known in the study of Hopf algebras \cite{molnar_1977}, we use slightly non-standard conventions (for instance in \cite{andruskiewitsch2006} the definitions have the sense of the action and coaction, as well as the order of tensor coefficients, exchanged) to achieve the accurate description of examples in the forthcoming sections.
\begin{definition}
    \label{smashBiproduct}
    Let $U, H$ be usual Hopf algebras. Let
    \[
        \lhd: H \otimes U \rightarrow H
    \]
    be a right algebra action of $U$ on $H$ (see for instance \cite[Chapter III.5]{kassel_2012}). We define their \textit{right semi-direct product} as the algebra over the vector space $U \otimes H$ with the unit $\1_U \otimes\1_H$ and the product
    \[
     (u\otimes h)(u' \otimes h') := u u'_{(1)} \otimes (h \lhd u'_{(2)}) h'.
    \]
   Let also
    \[
        \delta: H \rightarrow H \otimes U
    \]
    be a right coaction of $U$ on $H$ (see for instance \cite[Chapter III.6]{kassel_2012}), and we adapt the notation
    \[
        h\mapsto h^{[0]} \otimes h^{[1]}.
    \]
    We define their \textit{right semi-direct coproduct}, as the coalgebra over the vector space $U \otimes H$ with the counit $\epsilon_U \otimes \epsilon_H$, and the coproduct
    \[
        \Delta(u\otimes h) := (u_{(1)}\otimes h^{[0]}_{(1)})\otimes(h^{[1]}_{(1)} g_{(2)} \otimes h_{(2)}).
    \]
     Now if the semi-direct product and coproduct give rise to a Hopf algebra structure with the antipode $S:U \otimes H \rightarrow U \otimes H$
     \[
        S(u \otimes h) \mapsto S_U(u h^{[1]}) \otimes S_H(h^{[0]})
     \]
     and its inverse
     $S^{-1}:U \otimes H \rightarrow U \otimes H$
     \[
        S^{-1}(u \otimes h) \mapsto S^{-1}_U(u h^{[1]}) \otimes S^{-1}_H(h^{[0]}),
     \]
     then we will call it a \textit{semi-direct biproduct} of $U$ and $H$ and denote $U \ltimes H$.
\end{definition}%very ugly
\noindent Note that $\rhd$ always refers to the left \textit{adjoint} action, while $\lhd$ is always the right action of Definition \ref{smashBiproduct}.
\begin{remark}
    We emphasize that all tensor products between $U$ and $H$ in Definition \ref{smashBiproduct} are taken in the category $\text{\textbf{Vect}}_\C$ (resp. $\text{\textbf{Vect}}_\K$), which carries the trivial braiding. In particular, this is a different setting than the \textit{Radford biproduct} which is taken in a braided  \textit{category of Yetter-Drinfeld modules} over the acting algebra. 
\end{remark}

\subsection{Basic representation theory}
We first recall some basic notions in representation theory, see for instance \cite[Sections IV, VIII]{curtis_reiner_1966}, or \cite[Chapter 3]{drozd_kirichenko_2012}.
Let $A$ be a finite-dimensional associative algebra over an algebraically closed field, usually $\C$ for characteristic 0 below or $\K$ for prime characteristic, and $A\operatorname{-mod}$ be its category of modules. Recall $A\operatorname{-mod}$ contains all $A$-modules as objects and module maps as morphisms, but we will be concerned exclusively with finite-dimensional ones. 
\begin{notation}
    Let the left \textit{regular} representation of an associative algebra be denoted $(A, \cdot)$, or sometimes just $A$ if it is clear from the context. 
\end{notation}
\begin{definition}
    Let $H$ be a Hopf algebra and $W, V \in H\operatorname{-mod}$ be $H$-modules under some generic actions $\cdot_W$ and $\cdot_V$. Then $W\otimes V$ is an $H$-module with the action of $h \in H$
    \[
        h\cdot(W \otimes V):= h_{(1)} \cdot_{W} W\otimes h_{(2)} \cdot_{V} V.  
    \]
    The decomposition of such a module into a direct sum of indecomposable modules is called the \textit{fusion rule} of $W$ and $V$.
\end{definition}
 \begin{notation}
    Let $h\in H$ be an element of a Hopf algebra, and let $V\in H\operatorname{-mod}$. We say $h$ acts \textit{trivially} on $V$ if $h V = \epsilon(h) V$  
\end{notation}
\begin{definition}
    Define the \textit{indecomposable} module to be $V \in A\operatorname{-mod}$ that is one not isomorphic to a direct sum $V = V_1 \oplus V_2$, for some $V_1, V_2 \in A\operatorname{-mod}$. Indecomposable modules comprising the regular module are exactly the \textit{indecomposable projective} modules.
\end{definition}
\begin{definition}
    Let $e_a \in A$ be \textit{idempotent}, that is $e^2_i = e_i$. We call an idempotent $e_i$ \textit{primitive}, if it cannot be written as $e_i = e_j + e_k$, for some idempotents $e_j, e_k$.
    A \textit{full system of primitive orthogonal idempotents}. is a finite set of such $e_i$ for $i, j=1,\dots,s$ that $e_i e_j = \delta_{ij} e_i$, $e_i$ are all primitive and $\sum^s_{i=1} e_i = \1$.  
\end{definition}
\begin{proposition}
\label{regDecomposition}
    Let $(A, \cdot)$ be the regular representation, and $\{e_i\}_{1\leq i \leq n} \subset A$ for $n \in \Z_{>0}$ be a system full system of primitive orthogonal idempotents. Then we have the decomposition
        $(A, \cdot) = \bigoplus^n_i P_i,$
    where $P_i:=A e_i$ are the indecomposable projectives.
\end{proposition}
\noindent In our setting, any other projective module is always isomorphic to a direct sum of $P_i$.
\begin{definition}[\cite{etingof_gelaki_nikshych_ostrik_2016}, Definition 1.6.6]
\label{defprojcover}
    Let $V \in A\operatorname{-mod}$. The \textit{projective cover} of $V$, denoted $P_V$ is the projective module equipped with a surjection of modules $q_V: P_V \twoheadrightarrow V$ such that for any other projective $P$ and a surjection of modules $q: P \twoheadrightarrow V$, there is a surjection of modules $p: P \twoheadrightarrow P_V $, such that $q = q_V \circ p$. It is unique up to a unique isomorphism.
\end{definition}

\begin{notation}
    We will also encounter modules over a Hopf algebra $H$ under the \textit{left adjoint action} of Definition \ref{adjointaction}, which we will also denote with the symbol $\rhd$. In particular, by $(H, \rhd)$ we mean the representation of the Hopf algebra $H$ on itself with the left adjoint action.
\end{notation}

\noindent We will also use the following objects.
\begin{definition}[\cite{curtis_reiner_1966},  Definition (58.1)]
    The \textit{socle} of the module $V$ is the union of the semisimple submodules of $V$
    \[
        \mathrm{Soc}(V)\subset V := \left\{\bigcup M | M\subset V \;\mathrm{semisimple}\right\} 
    \]
\end{definition}
 
 \begin{definition}
     An ideal $I\subset A$ is called \textit{nilpotent} if $I^n = 0$ for some $n\in \N$
 \end{definition}
 
\begin{definition}[\cite{drozd_kirichenko_2012}, Corollary 3.1.9]
    The \textit{radical} of the algebra $A$, $\text{Rad}(A)$ is the nilpotent ideal containing all nilpotent left and right ideals. 
\end{definition}
\noindent In particular for $A$ finite dimensional $\text{Rad}(A)$ comprises all nilpotent elements of $A$.
\begin{definition}[\cite{drozd_kirichenko_2012}, Section 3.1]
    The \textit{radical} of the module $V\in A\operatorname{-mod}$, is $\text{Rad}(V):= \text{Rad}(A) \cdot V$, where $\cdot$ is the action of $A$ on the module $V$. 
\end{definition}
\begin{definition}[\cite{benson_2008}, Section 1.1]
     The \textit{head} of $V$ is the semisimple quotient module 
        $\Bar{V}:=V/\mathrm{Rad}(V).$
      It is simple if $V$ is indecomposable.
\end{definition}
\begin{notation}
     \noindent We will often consider modules generated by a single element, for instance $v\in V$, that we will call the \textit{generator} of the module. To indicate this we will write $V = A\cdot v = \langle v \rangle$.
     % In those cases, the head will be the simple module constructed as $\C\text{-}$vector space $\C\langle \omega_\mb{p}\rangle$ or $\C\langle g \rangle$ built upon the generator. In particular, we will sometimes speak of the socle and head of the algebra $A$ as the regular module.
\end{notation}
\noindent Let $\mathbbm{1}$ be the monoidal unit of the monoidal category $H\operatorname{-mod}$.
\begin{definition}
    \label{invModule}
    We call a module $V \in H\operatorname{-mod}$ \textit{invertible} if there exists a module $W \in  H\operatorname{-mod}$ such that
       $ V \otimes W = W \otimes V = \mathbbm{1},$
    where $\mathbbm{1} \in H\operatorname{-mod}$ is the trivial representation.
\end{definition}

Let $(H, R, v)$ be a ribbon Hopf algebra. Then $H\operatorname{-mod}$ is in particular rigid, so every object $V \in H\operatorname{-mod}$ has a two-sided dual denoted by $V^* \in H\operatorname{-mod}$ (see \cite[Chapter 2.10]{etingof_gelaki_nikshych_ostrik_2016}). We can use this to give a categorical equivalent of Definition \ref{defunimod}.
\begin{definition}
    The ribbon Hopf algebra $(H, R, v)$ is called unimodular if $P_{\mathbbm{1}} \cong P^*_{\mathbbm{1}}$ as $H$-modules.
\end{definition}
\noindent Equivalently, in such case $H\operatorname{-mod}$ is a unimodular category.

  % \begin{definition}
  %     Let $H$ be a Hopf algebra and $\Bar{H}$ its coradical, that is the head of the regular representation of $H$ regarded as a Hopf algebra. Then $H$ is called \textit{pointed} if $\Bar{H}\subset H$ as a Hopf subalgebra. 
  % \end{definition}
  % \noindent Examples of pointed Hopf algebras include the Lusztig $u_q \mathfrak{sl}_2$ at a root of unity and any instance of the Nenciu construction.\\\\
\begin{proposition}
 Let $(H, R)$ be a quasitriangular Hopf algebra. Then the category $H\operatorname{-mod}$ is braided monoidal. Let $V_1, V_2 \in H\operatorname{-mod}$. Then the braiding is given by the isomorphism
    $$c(V_1 \otimes V_2):=R'' \cdot V_2 \otimes R' \cdot V_1 \cong V_2 \otimes V_1,$$
where $\cdot$ is the generic module operation relevant to each module in question. 
\end{proposition}
\begin{definition}%[\cite{shimizu_2016}, Theorem 1.1]
\label{Muger}
    We call an object $V_1 \in H\operatorname{-mod}$ \textit{transparent}, if for any $V_2 \in H\operatorname{-mod}$, we have \\ $c^2(V_1 \otimes V_2) = V_1 \otimes V_2.$
The transparent objects of $H\operatorname{-mod}$ form the full subcategory 
\[
    \mathcal{Z}_{(2)}(H\operatorname{-mod}):=\{V_1 \in  H\operatorname{-mod} \,|\, c^2(V_1 \otimes V_2) = V_1 \otimes V_2,\;\; V_2 \in H\operatorname{-mod}\} \subset H\operatorname{-mod}
\]
called the \textit{M\"uger center}.
\end{definition}
\noindent Recall from Definition \ref{factorizability} that a quasitriangular Hopf algebra $(H, R)$ is factorizable if $\mathcal{Z}_{(2)}(H\operatorname{-mod})$ is trivial. Otherwise it is non-factorizable.

\section{The Nenciu construction}
In this section we recall some of the results from \cite{faes2025non}, and refer there for the corresponding proofs. The Nenciu construction of \cite{nenciu_2004} produces finite-dimensional ribbon Hopf algebras generated by grouplike generators collected in a tuple $\mb{K} = (K_1, \dots, K_s)$ and nilpotent generators collected in a tuple $\mb{X}=(X_1, \dots, X_t)$, such that the group subalgebra $\C[\mb{K}]$ is abelian, the nilpotent generators are skew-primitive with respect to $\C[\mb{K}]$ and all commutation relations between generators are diagonal. The choice of the generators and relations is strongly constrained, and the sufficient conditions for the result to be a ribbon Hopf algebra were determined in \cite{nenciu_2004} and collected in Theorem 3.19 \cite[Theorem 3.19]{faes2025non}.
\subsection{Hopf algebra structure}
We start by recalling Nenciu's construction: the fact that Definition \ref{genDef} indeed defines a Hopf algebra is a Theorem of \cite{nenciu_2004}.
\begin{definition}
\label{genDef}
    Let $\mb{m} \in \Z_{>0}^s$ be a tuple of positive integers of length $s$, and $t\in \Z_{>0}$ and $\matfont{d}$, $\matfont{u}$ be dimension $t\times s$, $\Z$-valued matrices, such that
    \begin{equation}
\label{invConstr}
    \pmb{\xi}^{\matfont{d}_k \cdot \matfont{u}_l}\pmb{\xi}^{\matfont{d}_l \cdot \matfont{u}_k}=1 \;\; \text{ and } \;\; \pmb{\xi}^{\matfont{d}_k \cdot \matfont{u}_k}=-1,
\end{equation}
    where $\pmb{\xi}:=(\exp^{2\pi/m_1}, \dots, \exp^{2\pi/m_s})$ is the tuple of primitive roots of unity of orders given by the tuple $\mb{m}\in \Z^s$.\\
    Define $H(\mb{m}, t, \matfont{d}, \matfont{u})$ to be the Hopf algebra generated by grouplike generators $\mb{K}=(K_1, \dots, K_s)$, and skew-primitive generators $\mb{X}=(X_1, \dots, X_t)$, and with the relations 
    \begin{align*}
        K^{m_a}_a=\1, && K_a K_b = K_b K_a, && K_a X_k := \xi_a^{d_{ka}} X_k K_a,
        && X^2_k = 0,  && X_l X_k := \pmb{\xi}^{\matfont{d}_k \cdot \matfont{u}_l} X_k X_l,
    \end{align*}
    coalgebra structure determined by
    \begin{align*}
        \epsilon(K_a) &:= 1, &\Delta(K_a) &:= K_a \otimes K_a, \\
        \epsilon(X_k) &:= 0, &\Delta(X_k) &:= \1 \otimes X_k + X_k, \otimes \mb{K}^{\matfont{u}_k},
    \end{align*}
    and antipode determined by
    \begin{align*}
        S(K_a) = K^{-1}_a := K^{m_a-1}_a, && S(X_k) := -X_k \mb{K}^{-\matfont{u}_k},
    \end{align*}
    for $a, b = 1, \dots, s$, $s:=|\mb{m}|$ and $k, l = 1, \dots, t$. Here $\mb{m}$ is the tuple of orders of the grouplike elements, $\matfont{d}$ encodes the diagonal relations for the grouplike elements, and $\matfont{u}$ the relations for the skew-primitive elements.
\end{definition}
\begin{remark}
    Note that in Definition \ref{genDef} the relations $X^2_k = 0$, for $k=1, ..., t$, result from the requirement that
        $\pmb{\xi}^{\matfont{d}_k \cdot \matfont{u}_k}=-1,$
    holds, rather than being imposed independently. We include them in the definition for clarity.
\end{remark}
% \begin{remark}
%     It is easy to see that 
%     \begin{align*}
%         S^2(K_a) = K_a, && S^2(X_k) = -X_k,
%     \end{align*}
%     that is the antipode restricted to any $K_a$ is of order 2, and to any $X_k$ of order 4.
% \end{remark}
\noindent Sometimes, we will choose several variables $X_k$ to have the same commutation properties. We refer to this as \textit{type}, and in examples we will indicate each type with a different letter $X, Y, Z$ etc. if necessary. 
\begin{definition}
    By the \textit{type} of $X_k$ we mean the prescription of $\matfont{d}_k$ and $\matfont{u}_k$, that is $X_l$ is of the same type as $X_k$ if $\matfont{d}_k=\matfont{d}_l$ and $\matfont{u}_k=\matfont{u}_l$. Since all relations involving $X_k$ and $X_l$ are diagonal, we will call the two generators of \textit{opposite type} if for $a=1, \dots, s$
    \begin{align*}
         d_{ka} \equiv -d_{la} \;\; \mod \; m_a && \text{ and } &&  u_{ka} \equiv -u_{la} \;\; \mod \; m_a.
    \end{align*}
This essentially means that generators of opposite type commute with all other generators over reciprocal constants. We will sometimes indicate variables of opposite type using superscripts, for instance $X^+_k$, $X^-_k$.
\end{definition}
\begin{proposition}
     Let $X^+_k, X^-_k$ be nilpotent generators of opposite type, then $\{X^+_k, X^-_k\}=0$.
\end{proposition}
\begin{proof}
We have directly
    $
    X^+_k X^-_k= \pmb{\xi}^{\matfont{d}_k \cdot \matfont{u}_k}X^-_k X^+_k = (-1)^{-1} X^-_k X^+_k = -X^-_k X^+_k.
    $
\end{proof}
\begin{remark}
    When making general statements about $H(\mb{m}, t, \matfont{d}, \matfont{u})$ that do not invoke the types we will use the notation $X_k$ for any nilpotent generator, with no superscripts. In such sense, two nilpotent generators with distinct indices, say $X_k$, $X_l$, can be of different, opposite or same type. 

% \begin{enumerate}[leftmargin=1cm]
%     \item The key point of the construction is the constraint \eqref{invConstr}, expressing the compatibility between the order of grouplike generators and the nilpotent generators. Once this is satisfied, the Hopf algebra axioms follow readily.
%     \item It is easy to see that $S^2(X_k) = -X_k$ by the prescribed commutation relations. So the antipode restricted to skew-primitive generators is of order 4.
% \end{enumerate}
\end{remark}
\noindent We can define a \textit{monomial basis} for $H(\mb{m}, t, \matfont{d}, \matfont{u})$.
\begin{proposition}
    \label{monomialBasis} The Hopf algebra $H(\mb{m}, t, \matfont{d}, \matfont{u})$  has a monomial basis
    $
        \{\mb{K}^\mb{v} \mb{X}^\mb{r} |\mb{v}\in \Z_\mb{m},  \mb{r} \in \Z^t_2\}.
    $
    % That is, any element $h \in H(\mb{m}, t, \matfont{d}, \matfont{u})$ can be uniquely written as 
    % \[
    %     h = \sum_{\mb{v}\in \Z_\mb{m}, \mb{r} \in \Z^t_2} \gamma_{\mb{v}, \mb{r}} \mb{K}^\mb{v} \mb{X}^\mb{r} 
    % \]
    % for some $\gamma_{\mb{v}, \mb{r}} \in \C$. 
\end{proposition}

\begin{definition}
    We call $T := \mb{X}^\mb{r}$ for $\mb{r}=(1, \dots, 1)$ the \textit{top element}. 
\end{definition}
\noindent It will play an important role in the unimodular structure we will sometimes define on $H(\mb{m}, t, \matfont{d}, \matfont{u})$.
\noindent We will also need to use the adjoint action of generators on one another (see Definition \ref{adjointaction}). By direct computation, we get:
\begin{proposition}
\label{adjointAct}
     The adjoint action of the generators of $H$ is given by the following formulas.
     %consider removing this part
     A grouplike $\mb{K}^\mb{w}$, $\mb{w} \in \Z_{\mb{m}}$ acts on the skew-primitive generator $X_k$ as
    \[
        \mb{K}^\mb{w} \rhd X_k = \mb{K}^\mb{w} X_k (\mb{K}^\mb{w})^{-1} = \pmb{\xi}^{\mb{w} \cdot \matfont{d}_k} X_k,
    \]
    a skew-primitive $X_k$ acts on a grouplike $\mb{K}^\mb{w}$ as
    \[
        X_k \rhd \mb{K}^\mb{w} = \mb{K}^\mb{w} S(X_k) + X_k \mb{K}^\mb{w} \mb{K}^{-\matfont{u}_k} = (1-\pmb{\xi}^{\mb{w} \cdot \matfont{d}_k})X_k \mb{K}^{-\matfont{u}_k}\mb{K}^{\mb{w}},
    \]
    and a skew-primitive $X_k$ acts on a skew-primitive $X_l$ as 
    \[
        X_k \rhd X_l = X_l S(X_k) + X_k X_l \mb{K}^{-\matfont{u}_k} = (\pmb{\xi}^{\matfont{u}_k \cdot \matfont{d}_l}-1)X_k \mb{K}^{-\matfont{u}_k}X_l.
    \]
\end{proposition}
\subsection{Quasitriangular structure}
As we will consider the bariding structuire of module catagories later on, for the reader's convenience, we recall the form of the R-matrices that arise from the Nenciu construction.
\begin{notation}
    Let $h, g\in H$ be elements of a Hopf algebra. We the define the \textit{exponential} in $H\otimes H$ by
    \[
        \exp (h \otimes g) := \sum^\infty_{i=0} \frac{1}{i!} (h\otimes g)^i, 
    \]
    where the products are taken in $H\otimes H$, as remarked in Definition \ref{Hopf_alg_def}. We will often encounter exponentials of sums, and then we require the summands to commute in $H\otimes H$.
\end{notation}
\begin{theorem}[\cite{nenciu_2004}, Theorem 4.1]
    \label{R-matrixConstr} 
     Let $R \in H(\mb{m}, t, \matfont{d}, \matfont{u}) \otimes H(\mb{m}, t, \matfont{d}, \matfont{u})$ be  the element defined by $R:= R_{\matfont{z}} R_{\pmb{\alpha}},$ for
    \begin{align*}
         R_{\matfont{z}} = \frac{1}{\prod^s_{a=1} m_a}\sum_{\mb{v}, \mb{w} \in \Z_\mb{m}} \pmb{\xi}^{-\mb{v}\cdot \mb{w}} \mb{K}^{\mb{w}} \otimes \mb{K}^{\mb{v}\matfont{z}}, && \text{ and } && R_{\pmb{\alpha}} = \exp \left( \sum^{t}_{k, l=1} \alpha_{kl} X_k \otimes \mb{K}^{-\matfont{u}_l}X_l \right), 
    \end{align*}
    where $\matfont{z}\in \operatorname{Mat}_{s\times t}(\Z_\mb{m})$ and $\pmb{\alpha} \in \operatorname{Mat}_{t\times t}(\C)$, are matrices. Then $R$ is an R-matrix, if and only if the constraints of \cite[Theorem 3.19]{faes2025non} are satisfied.
\end{theorem}
\noindent For unimodular and ribbon properties see \cite[Section 3.2]{faes2025non}.

\subsection{Examples}
     In this subsection we recall from \cite{faes2025non} some explicit examples of the Nenciu construction, and their unimodular, strongly non-factorizable quasitriangular and ribbon structures. 
 \begin{notation}
     In this section the nilpotent component of the R-matrix of Theorem \ref{R-matrixConstr},  $R_{\pmb{\alpha}}$ will recur in the same form given by
     \[
        R_{\pmb{\alpha}} := \exp\left(\sum^{t_2}_{l = 1} \alpha_l( Z^+_l \otimes L Z^-_l - Z^-_l \otimes L Z^+_l)\right).
     \]
    Here we stick to the notation for the nilpotent geenrators emphasizing the types. In particular, let $Z^+_l$ and $Z^-_l$ be nilpotent generators in $H(\mb{m}, t, \matfont{d}, \matfont{u})$ of opposite type, with neighbouring indices, that is $Z^+_l = X_k$ and $Z^-_l = X_{k+1}$ for some $k$. By a slight abuse of notation, in the expression above we denoted $\alpha_{k, k+1} = -\alpha_{k+1, k}=:\alpha_l$. We collect $\alpha_l$ in the tuple $\pmb{\alpha}$. 
 \end{notation} 
We start with a (not generically strongly) non-factorizable Hopf algebra, which can be retrieved from Nenciu's construction, and appears in \cite{nenciu_2004} as $E(n)$, first studied in \cite{panaite_oystaeyen_1999}. Here we will refer to it as \textit{symplectic fermions}, after \cite{gainutdinov_runkel_2017}, also appearing in \cite{faitg_gainutdinov_schweigert_2024}. The notation using the variables $L, Z^\pm_l$ is introduced to be used in further examples.
 
 \begin{definition}
 \label{SFdef}
     By the \textit{Hopf algebra of symplectic fermions}, $\operatorname{SF}_{2n}$ we mean the Hopf algebra generated by $K_1 = L$ and, $X_{2l-1} = Z^+_l$, $X_{2l} = Z^-_l$, for $k, l=1, ..., n$, subject to the following relations
    \[
         L^2 = \1, \;\; L Z^\pm_l = - Z^\pm_l L, \;\; \{Z^\pm_l, Z^\pm_k\}=\{Z^\pm_l, Z^\mp_k\}=0.
    \]
    The Hopf structure is defined by
    \[
        \epsilon(L) = 1, \;\; \epsilon(Z^{\pm}_l)= 0,
    \]
    \[
        \Delta(L) = L \otimes L, \;\; \Delta(Z^{\pm}_l) = \1 \otimes Z^{\pm}_l + Z^{\pm}_l \otimes L,
    \]
    \[
        S(L)=L^{-1}=L, \;\; S(Z^{\pm}_l) = -Z^{\pm}_l L.
    \]
    % The dimension of the algebra is $2^{2n+1}$. The corresponding Nenciu data is 
    % \begin{itemize}
    %     \item $\mb{m}=2$ is just a number, so $s=1$
    %     \item $t=2n$
    %     \item $\matfont{d} = \matfont{u} = (1, ..., 1)^T$ are column vectors whose each row $d_k = 1$ and $u_k = 1$, $k=1, ..., t$
    %      \begin{align*}
    %         \matfont{d}=\begin{pmatrix}
    %             1\\
    %             \multicolumn{1}{c}{$\vdots$} \\
    %             1\\
    %         \end{pmatrix},&&
    %         \matfont{u}=\begin{pmatrix}
    %             1\\
    %             \multicolumn{1}{c}{$\vdots$} \\
    %             1\\
    %         \end{pmatrix}
    %     \end{align*}
    % \end{itemize}
 \end{definition}
 \begin{notation}
     In line with the usual terminology, we will call monomials in $\operatorname{SF}_{2n}$ commuting with $L$ \textit{even} and anti-commuting \textit{odd}.
 \end{notation}
\begin{proposition}
\label{SFprop}
    The algebra $\operatorname{SF}_{2n}$ of Definition \ref{SFdef} is
    \begin{enumerate}
        \item unimodular, with a two-sided cointegral
        \[
            \Lambda := (\1+L)\prod^{n}_{l=1}Z^+_l
            \prod^{n}_{l=1}Z^-_l.
        \]
        and a two-sided integral defined on the monomial basis by
        \[
            \lambda(L^v \mb{X}^\mb{r}):=
            \begin{cases}
                1 \text{\, if \,} v = 0, \mb{r} = (1, 1,...,1) \\
                0 \text{\, otherwise}.
            \end{cases}
        \]
 \item quasitriangular, with the R-matrix defined for the grouplike generators as 
        \[
            R_{\matfont{z}} := \frac{1}{2}\sum_{v, w \in \Z_2} (-1)^{-vw} L^{w} \otimes L^{v z} \exp\left(\sum^{n}_{l = 1} \alpha_l( Z^+_l \otimes L Z^-_l - Z^-_l \otimes L Z^+_l)\right),
        \]
        where $z=1$ and $\pmb{\alpha}=(\alpha_1, ..., \alpha_{t_2}) \in \C^{2t_2}$,
        \item ribbon, with the ribbon element
        \[
            v:= \exp \left( - 2\sum^{n}_{k = l} \alpha_l Z^+_l Z^-_l  \right),
        \]
        corresponding to the pivotal element $g = L$, which in this example is the only such element available. 
    \end{enumerate}
\end{proposition}

The following example corresponds to $\operatorname{N}_4$  \cite[Example 4.17]{faes2025non}, but for present purposes we relabel it as $\operatorname{N}_1$. This is due to the fact that from representation-theoretic perspective $\operatorname{N}_1$ of \cite[Example 3.22]{faes2025non} and $\operatorname{N}_4$ of \cite[Example 4.17]{faes2025non} share many characteristics.
\begin{example}
\label{auxExample}
    Let $\operatorname{N}_1:=H(\mb{m}, t, \matfont{d}, \matfont{u})$ be a Nenciu type algebra generated by grouplike $K_1, K_2$ and nilpotent $X^\pm$ generators, and the following relations
    \[
        K^4_1 = K^4_2 = \1, \;\;(X^\pm)^2 = 0, \;\; 
    \]
    \[
        K X^\pm K^{-1} =q^{r'}X^\pm = -X^\pm, \;\; K_1 X^{\pm} K^{-1}_2 = K_2 X^{\pm} K^{-1}_2 = \pm i X^\pm. 
    \]
    Let $L:=K^{r'/2} K_1 K_2$ for convenience. The Hopf structure is 
    \[
        \epsilon(K_1) = \epsilon(K_2)=1, \;\; \epsilon(X^\pm) = 0, 
    \]
    \[
        \Delta(K_1) = K_1 \otimes K_1, \;\; \Delta(K_2) = K_2 \otimes K_2, \;\; \Delta(X^\pm) = \1 \otimes X^\pm + X^\pm \otimes L^{\pm 1},
    \]
    and
    \[
        S(K_1) = K_1^{-1}, \;\; S(K_2) = K_2^{-1}, \;\; S(X^\pm) = - X^\pm L^{\mp 1}.
    \]
\end{example}
\begin{proposition}
     The Hopf algebra $\operatorname{N}_1$ of Example \ref{auxExample} is
    \begin{enumerate}
        \item unimodular, with a two-sided cointegral
        \[
            \Lambda := \left(\sum^3_{a,b=0} K^a_1 K^b_2 \right) X^+X^- 
        \]
        and a two-sided integral defined on the monomial basis by
        \[
            \lambda(\mb{K}^\mb{v}\mb{X}^\mb{r} ):=
            \begin{cases}
                1 \text{\, if \,} \mb{v} = (0, 0), \mb{r} = (1, 1) \\
                0 \text{\, otherwise}.
            \end{cases}
        \]
    \item triangular, with the R-matrix defined for the grouplike elements as 
         \[
             R_{\matfont{z}} := \frac{1}{16}\sum_{\mb{v}, \mb{w} \in \Z^2_4} i^{-\mb{v}\mb{w}^T} \mb{K}^{\mb{w}} \otimes \mb{K}^{\mb{v} \matfont{z}}, 
         \]
         where $\mb{K} = (K_1, K_2)$, $\pmb{\xi} = (i, i)$, and $\matfont{z} =\begin{pmatrix} 2 && 3\\ 1 &&0 \end{pmatrix}$,
        \item trivially ribbon, with the ribbon element $\1$, corresponding to the unique pivotal element $g = K^2_1$. 
    \end{enumerate}
\end{proposition}
\begin{example}
    \label{ex2}
    Let $\operatorname{N}_2$ be the Hopf algebra generated by $K_a$, $X_k^{\pm}$ and $Z^{\pm}_l$, for $a=1, 2, 3$, $j, k=1, ..., t_1$; $ t_1\in \N$ and $l, m =1, ..., t_2$; $t_2\in \N$, subject to the following relations
    \[
        K_a^4 = \1, \;\;  \;\; K_a X^{\pm}_k = \pm i X^{\pm}_k K_a
    \]
    \[
       K_1 Z^{\pm}_l = Z^{\pm}_l K_1, \;\; K_2 Z^{\pm}_l = - Z^{\pm}_l K_2, \;\;  K_3 Z^{\pm}_l = \pm i Z^{\pm}_l K_3
    \]
    \[
        \{X^{\pm}_j, X^{\pm}_k\} = \{X^{\pm}_j, X^{\mp}_k\} = \{Z^{\pm}_l, X^{\pm}_k\} = \{Z^{\pm}_l, X^{\mp}_k\} = \{Z^{\pm}_l, Z^{\pm}_m\}=\{Z^{\pm}_l, Z^{\mp}_m\} =  0.
    \]
    Let also $L:=K_3^2$ as a shorthand, note that $L^2 = \1$. The Hopf structure is defined by
    \[
        \epsilon(K_a) = 1, \;\; \epsilon(X^{\pm}_k)=\epsilon(Z^{\pm}_l) = 0,
    \]
    \[
        \Delta(K_a) = K_a \otimes K_a, \;\; \Delta(X^{\pm}_k) = \1 \otimes X^{\pm}_k + X^{\pm}_k \otimes (K_1 K_2)^{\pm 1}, \;\; \Delta(Z^{\pm}_l) = \1 \otimes Z^{\pm}_l + Z^{\pm}_l \otimes L
    \]
    \[
       S(K_a)=K^{-1}_a \;\; S(X^{\pm}_k) = -X^{\pm}_k (K_1 K_2)^{\mp1} \;\;  S(Z^{\pm}_l) = -Z^{\pm}_l L.
    \]
\end{example}

\begin{proposition}
\label{ex2prop}
     The Hopf algebra $\operatorname{N}_2$ of Example \ref{ex2} is
    \begin{enumerate}
        \item unimodular, with a two-sided cointegral
        \[
            \Lambda := \left(\sum^3_{a,b, c=0} K^a_1 K^b_2 K^c_3\right) \prod^{t_1}_{k=1}X^+_k \prod^{t_1}_{k=1}X^-_k \prod^{t_2}_{l=1}Z^+_l
            \prod^{t_2}_{l=1}Z^-_l.
        \]
        and a two-sided integral defined on the monomial basis as
        \[
            \lambda(\mb{K}^\mb{v}\mb{X}^\mb{r} ):=
            \begin{cases}
                1 \text{\, if \,} \mb{v} = (0, 0, 0), \mb{r} = (1, 1,...,1) \\
                0 \text{\, otherwise}.
            \end{cases}
        \]
        % \[
        %     \lambda\left(\prod^{t_1}_{k=1}(X^+_k)^{a_k}\prod^{t_1}_{k=1}(X^-_k)^{b_k} \prod^{t_2}_{l=1}(Z^+_l)^{c_l}
        %     \prod^{t_2}_{l=1}(Z^-_l)^{c_l}\right):= \prod^{t_1}_{k=1}\delta_{1, a_k} \prod^{t_1}_{k=1}\delta_{1, b_k}
        %     \prod^{t_2}_{l=1}\delta_{1, c_l}\prod^{t_2}_{l=1}\delta_{1, d_l}.
        % \]
        \item quasitriangular, with the R-matrix 
        \[
            R_{\matfont{z}} R_{\pmb{\alpha}} := \frac{1}{64}\sum_{\mb{v}, \mb{w} \in \Z^3_4} i^{-\mb{v}\mb{w}^T} \mb{K}^{\mb{w}} \otimes \mb{K}^{\mb{v} \matfont{z}} \exp\left(\sum^{t_2}_{l = 1} \alpha_l( Z^+_l \otimes L Z^-_l - Z^-_l \otimes L Z^+_l)\right),
        \]
         where $\mb{K} = (K_1, K_2, K_3)$, $\matfont{z} =\begin{pmatrix} 0 && 3 && 2 \\ 1 && 0 && 0 \\ 2 && 0 && 2  \end{pmatrix}$ and $\pmb{\alpha}=(\alpha_1, ..., \alpha_{t_2}) \in \C^{2t_2}$,
        \item ribbon, with the ribbon element
        \[
            v:= \exp \left( - 2\sum^{t_2}_{l = 0} \alpha_l Z^+_l Z^-_l  \right),
        \]
        corresponding to the unique pivotal element $g = L =K^2_3$.
    \end{enumerate}
\end{proposition}
\begin{remark}
    Note that Theorem \ref{QTthm} allows multiple quasitriangular structures for $H(\mb{m}, t, \matfont{d}, \matfont{u})$. Indeed, it is known from \cite{gainutdinov_runkel_2017} that $\operatorname{SF}_2$ admits a quasitriangular structure where $R:=R_z R_{\pmb{\alpha}}$, with $R_z$ as in Proposition \ref{SFprop}, but  
    \[
        R_{\pmb{\alpha}} = \exp \left( \alpha^+ Z^+ \otimes L Z^+ + \alpha^- Z^- \otimes L Z^- \right)
    \]
    where $\pmb{\alpha} = (\alpha^+, \alpha^-) \in \C^2$. This choice leads to a symmetric braiding in $\operatorname{SF}_2\operatorname{-mod}$ and admits a trivial ribbon structure $v=\1$.
\end{remark}

\begin{proposition}
\label{ExAreSNF}
    The Hopf algebras of Examples \ref{auxExample} and \ref{ex2} admit only strongly non-factorizable quasitriangular structures.
\end{proposition}

We finish off this section by recalling from \cite{faes2025non} the monodromy matrix corresponding to the class of R-matrices $R_{\matfont{z}}R_{\pmb{\alpha}}$ that appeared extensively in the examples.
\begin{proposition}
    \label{monodromyProp}
    Let $H(\mb{m}, t, \matfont{d}, \matfont{u})$ be a Hopf algebra containing nilpotent generators $Z_l^\pm$, for $l=1, ..., t_2$. Moreover let $Z^+_l$ and $Z^-_l$ be of opposite type and $\{Z_{l_1}^\pm, Z_{l_2}^\pm\}=\{Z_{l_1}^\pm, Z_{l_2}^\mp\}=0$. Denote $L_l:=\mb{K}^{\matfont{u}_l}$  Let $R:=R_{\matfont{z}}R_{\pmb{\alpha}}$ be an R-matrix of the form as in Theorem \ref{R-matrixConstr}, where 
    \[
        R_{\matfont{z}}=\frac{1}{\prod^s_{a=1} m_a}\sum_{\mb{v}, \mb{w} \in \Z_\mb{m}} \pmb{\xi}^{-\mb{v}\cdot \mb{w}} \mb{K}^{\mb{w}} \otimes \mb{K}^{\mb{v}\matfont{z}}
    \]
    and let
    \[
        R_{\pmb{\alpha}} = \exp\left(\sum^{t_2}_{l=1} \alpha_l (Z_l^+ \otimes L_l Z_l^{-} - Z_l^- \otimes L_l Z_l^{+})\right),
    \]
    where $\alpha_l \in \C$.
    Then the monodromy matrix is given by
    \[
        M:=R_{21}R =\exp\left(2\sum^{t_2}_{l=1} \alpha_l (Z_l^+ \otimes L_l Z_l^{-} - Z_l^- \otimes L_l Z_l^{+})\right).
    \] 
\end{proposition}

\section{Non-factorizable biproducts of $u_q \mathfrak{sl}_2$}

In this section we construct a family of algebras, where the $u_q \mathfrak{sl}_2$, for $q$ a primitive root of unity and $q^r = 1$, $r \equiv 0 \mod 8$ will be augmented with a Hopf algebra of Nenciu type. The idea of extending $u_q \mathfrak{sl}_2$ using nilpotent generators is inspired by the work of Majid, for instance \cite{majid_2000}.

\subsection{The small quantum group $u_q \mathfrak{sl}_2$}
We use this section to set the conventions for the small quantum group. We follow mostly \cite{beliakova_derenzi_2023}, and according to their conventions we let $r'=r/2$ and $r'' = r/4$.

\begin{definition}
\label{smallQG}
    The Hopf algebra $u_q \mathfrak{sl}_2$ for $q^r = 1$, $r \equiv 0 \mod 4$ is generated by the grouplike $K$ and nilpotent $E, F$, satisfying the relations 
    \[
        K^{r'} = \1, \;\; E^{r'} = F^{r'} = 0
    \]
    \[
        KEK^{-1} = q^2 E, \;\; KFK^{-1} = q^{-2} F, \;\; [E, F] = \frac{K-K^{-1}}{q-q^{-1}}. 
    \]
    The Hopf structure is given by
    \[
        \epsilon(K) = 1, \;\; \epsilon(E) = \epsilon(F) = 0,
    \]
    \[
        \Delta(K) = K \otimes K, \;\; \Delta(E) = \1 \otimes E + E \otimes K, \;\; \Delta(F) = F \otimes \1 + K^{-1} \otimes F
    \]
    and 
    \[
        S(K) = K^{-1} = K^{r'-1}, \;\; S(E) = -E K^{-1}, \;\; S(F) = -K F.
    \]
\end{definition}
\noindent Recall also the quantum integers be for $k \in \Z$
\[
    \{k\}:= q^k - q^{-k},\;\; [k]:=\frac{\{k\}}{\{1\}},\;\; [k]!:= [k][k-1]...[1].
\]
\begin{proposition}
\label{uqsl2prop}
    The Hopf algebra  $u_q \mathfrak{sl}_2$ is
    \begin{enumerate}
        \item unimodular with a two-sided cointegral
        \[
            \Lambda:=\frac{\{1\}^{r'-1}}{\sqrt{r''}[r'-1]!}\sum^{r'-1}_{a=0}  E^{r'-1} F^{r'-1} K^a,
        \]
        and a left integral
        \[
            \lambda(E^b F^c K^a) := \frac{\sqrt{r''}[r'-1]!}{\{1\}^{r'-1}} \delta_{a, r'-1}\delta_{b, r'-1}\delta_{c, r'-1},
        \]
        \item quasitriangular with $R:= D \Theta$, where
        \[
            D:= \frac{1}{r'}\sum^{r'-1}_{b, c=0} q^{-2bc} K^b \otimes K^c,
        \]
        and 
        \[
            \Theta:= \sum^{r'-1}_{a=0} \frac{\{1\}^a}{[a]!}q^{\frac{a(a-1)}{2}} E^a \otimes F^a,
        \]
        \item ribbon, with the ribbon element 
        \[
            v := \frac{1-i}{\sqrt{r'}} \sum^{r'-1}_{a=0}\sum^{r''-1}_{b=0} \frac{\{-1\}^a}{[a]!} q^{-\frac{(a+3)a}{2} + 2b^2} E^aF^a  K^{-a-2b-1}.
        \]
        corresponding to the pivotal element $g=K$.
    \end{enumerate}
\end{proposition}
\begin{remark}
\label{DasRz}
    To make a connection with the previous notation for $D$, note that $q^{2r'}=1$, but $K^{r'}=\1$. Thus, we can write
    \[
        D =  \frac{1}{r'}\sum_{v, w \in \Z_{r'}} (q^2)^{-wv} K^w \otimes K^{zv}
    \]
    for $z=1$.
\end{remark}

% \begin{remark}
%     We do not attempt to determine the conditions under which the semi-direct product and semi-direct coproduct in Definition \ref{smashBiproduct} give rise to the semi-direct biproduct. We merely present in the sequel a construction where the conditions are indeed satisfied.
% \end{remark}
% \begin{remark}
%     This is not a typical bosonisation - the compatibility is 
%     \[
%         g_{(1)} h^{[1]}\otimes h^{[0]}\lhd g_{(2)} = g_{(1)}(h\lhd g_{(2)})^{[1]}\otimes (h \lhd g_{(2)})^{[0]}.
%     \]
% \end{remark}

\subsection{Non-factorizable biproducts with $u_q \mathfrak{sl}_2$}
In this section we build biproducts of $u_q \mathfrak{sl}_2$ and $H(\mb{m}, t, \matfont{d}, \matfont{u})$, and we show they can have the properties of unimodularity, quasitriangularity and ribbonness, and can be (strongly) non-factorizable. For clarity of discussion and the applications we keep in mind, we restrict the order of $q$, $r$ to be a multiple of 8.
\begin{definition}
\label{extHopfDef}
    Let $U= u_q \mathfrak{sl}_2$ and $H = H(\mb{m}, t, \matfont{d}, \matfont{u})$, where $q$, is a primitive $r$-th root of unity, $r \equiv 0 \mod 8$ and $r'=r/2$, $r''=r/4$. We define $U \ltimes H$, where
    \begin{itemize}[leftmargin=1cm]
        \item the right action $\lhd: H(\mb{m}, t, \matfont{d}, \matfont{u}) \otimes  u_q \mathfrak{sl}_2 \rightarrow H(\mb{m}, t, \matfont{d}, \matfont{u})$ for $K_a \in H$ and  $X_k \in H$ is defined by
         \begin{align*}
              K_a \lhd K = K_a, \;\; K_a \lhd E = K_a \lhd F = 0, &&  X_k \lhd K = -X_k, \;\; X_k \lhd E = X_k \lhd F = 0, 
         \end{align*}
        \item the right coaction $\delta: H(\mb{m}, t, \matfont{d}, \matfont{u}) \rightarrow H(\mb{m}, t, \matfont{d}, \matfont{u}) \otimes  u_q \mathfrak{sl}_2$  for $K_a \in H$ and $X_k \in H$ is defined by
        \begin{align*}
            K_a \mapsto K_a \otimes \1_U, && X_k \mapsto X_k \otimes K^{r''},
        \end{align*}
        and on an arbitrary monomial $\mb{K}^\mb{v} \mb{X}^\mb{r} \in H$ for $\mb{v}\in \Z_\mb{m}$ and $\mb{r}\in \Z^t_2$ by
        \[
            \mb{K}^\mb{v} \mb{X}^\mb{r} \mapsto \mb{K}^\mb{v} \mb{X}^\mb{r} \otimes K^{|\mb{r}|r''}.
        \]
    \end{itemize}
    
\end{definition}
\begin{theorem}
\label{extHopfThm}
    Let $u_q \mathfrak{sl}_2 \ltimes H(\mb{m}, t, \matfont{d}, \matfont{u})$ be as in Definition \ref{extHopfDef}. Then it is a Hopf algebra with the following semi-direct biproduct structure:
    \begin{itemize}[leftmargin=1cm]
        \item the new unit $\1 = \1_U \otimes \1_H$ and counit $\epsilon := \epsilon_U \otimes \epsilon_H$,
        \item the algebra structure, where we suppress the tensor product structure 
    \begin{align*}
       K:= K \otimes \1_H, && E:= E \otimes \1_H, && F:= F \otimes \1_H, && K_a := \1_U \otimes K_a,  && X_k := \1_U \otimes X_k,
    \end{align*}
    so that for any $a, b, c, \in \Z_{r'}$, $\mb{v}\in \Z_\mb{m}$ and $\mb{r}\in \Z^t_2$
    \[
        E^a F^b K^c \mb{X}^\mb{r}\mb{K}^\mb{v}:= E^a F^b K^c \otimes \mb{X}^\mb{r}\mb{K}^\mb{v}
    \]
    and find the complete set of relations between the the elements of $U$ and $H$ determining  the structure of the new algebra to be 
    \begin{align*}
        [K, K_a] = [E, K_a] = [F, K_a] = 0, && \{K, X_k\} = \{E, X_k\} = [F, X_k]=0. 
    \end{align*} 
    \item new Hopf structure for $X_k$ given by
    \begin{align*}
        \Delta(X_k) = \1 \otimes X_k + X_k \otimes K^{r''}\mb{K}^{\matfont{u}_k} && S(X_k)  =  -X_k K^{r''}\mb{K}^{-\matfont{u}_k},
    \end{align*}
    while for the remaining generators $E, F, K$ and $K_a$ in this notation the coproduct $\Delta$ and antipode $S$ are unchanged.
    \end{itemize}
\end{theorem}
\begin{proof}
    See \cite[Theorem 4.8]{faes2025non}.
\end{proof}

\begin{proposition}
\label{extMonomialBasis}
The Hopf algebra $u_q \mathfrak{sl}_2 \ltimes H(\mb{m}, t, \matfont{d}, \matfont{u})$ has a monomial basis
    \[
       \{E^e F^f K^k \mb{K}^\mb{w}\mb{X}^\mb{r}|e, f, k =1, \dots, r', \mb{w} \in \Z_\mb{m}, \mb{r} \in \Z^t_2\}.
    \]
\end{proposition}
\begin{proof}
    See \cite[Proposition 4.9]{faes2025non}.
\end{proof}
\begin{remark}
\label{SESproposition}
    Let $U,H$ and $U\ltimes H$ be as in Definition \ref{extHopfDef}, in particular let the pivotal element of $H$ be $g_H$. Then we claim the following maps are Hopf algebra maps, which do \emph{not} give a short exact sequence:
    \[
        U\hookrightarrow U\ltimes H \twoheadrightarrow H.
    \]
    The map $\iota : U\hookrightarrow U\ltimes H$ is given by $u \mapsto u \otimes \1_H$ for $u \in U$, and has a left inverse $j: U \ltimes H\twoheadrightarrow U$ given by
    \begin{align*}
        j(K) = K, && j(E) = E, && j(F) = F, && j(K_a) = \1_U, && j(X_k) = 0. 
    \end{align*}
    The map $p : U\ltimes H \twoheadrightarrow H$, is given by 
    \begin{align*}
        p(K) = g_H, && p(E) = p(F) = 0, && p(K_a) = K_a, && p(X_k) = X_k.  
    \end{align*}
\end{remark}
\begin{remark}
\label{NonSplit}
    The map $p$ in Remark \ref{SESproposition} has a right inverse as an algebra map, given by the map $$q :H \hookrightarrow H\ltimes U, \;\; q(h)=\1_U\otimes h,$$ but it is \textit{not} a Hopf algebra map. For instance $S(X_k) = K^{r''} \mb{K}^{\matfont{u}_k}X_k$, so the image of $H$ under $q$ is not closed under taking antipodes (as it is not under taking coproducts).
    Nevertheless, it makes sense to speak, for instance, of elements $u, h \in U\ltimes H$ as belonging to $u \in U$ and $h\in H$, and we will often slightly abuse the notation in this way.
\end{remark}
% \begin{remark}
%     The map $p$ of Proposition \ref{SESproposition} is well defined only if $H$ carries a pivotal element . 
% \end{remark}
\begin{remark}
    Let $U,H$ and $U\ltimes H$ be as in Theorem \ref{extHopfThm}. Then in light of Remark \ref{SESproposition} $U\ltimes H$ is not isomorphic to $U\otimes H$ as Hopf algebra in the usual sense. In particular, $U\otimes H$ has no non-trivial commutation relations between the elements belonging to $U$ and $H$. Similarly, $U\ltimes H$ is not isomorphic to $U\oplus H$, that is the short exact sequence is not fully split, as per Remark \ref{NonSplit}.
\end{remark}

Now we state the main theorem of the section, for the proof we again refer to \cite[Theorem 4.13]{faes2025non}.
\begin{theorem}
\label{extMainThm}
    Let $U = u_q \mathfrak{sl}_2$ with 
    the two-sided cointegral $\Lambda_U$, the left integral $\lambda_U$, 
    the R-matrix $R_U=D\Theta$  as in Proposition \ref{smallQG}.\\
    Let $H = H(\mb{m}, t, \matfont{d}, \matfont{u})$ with 
%a two-sided cointegral $\Lambda_H$, a two-side integral $\lambda_U$
    an R-matrix $R_H := R_{\matfont{z}}R_{\pmb{\alpha}}$ where
        \[
            R_{\matfont{z}}=\frac{1}{\prod^s_{a=1} m_a}\sum_{\mb{v}, \mb{w} \in \Z_\mb{m}} \pmb{\xi}^{-\mb{v}\cdot \mb{w}} \mb{K}^{\mb{w}} \otimes \mb{K}^{\mb{v}\matfont{z}}
        \]
        \[
            R_{\pmb{\alpha}} = \exp\left( \sum^{t_2}_{l=0} \alpha_l(Z^+ \otimes LZ^- - Z^- \otimes LZ^+) \right),
        \]
     and a ribbon element
         \[
             v_H:= \exp \left( -2 \sum^{t_2}_{k = 1} \alpha_l Z^+_l Z^-_l  \right).
         \]
    Then, $U\ltimes H$ carries the 
    \begin{enumerate}
         \item a two-sided cointegral
         \[
             \Lambda:=\Lambda_U \Lambda_H = \frac{\{1\}^{r'-1}}{\sqrt{r''}[r'-1]!}\sum^{r'-1}_{a=0}  E^{r'-1} F^{r'-1} K^a 
              \prod^s_{a=1} \left(\sum^{m_a-1}_{b=0} K^b_a \right) \prod^t_{k=1} X_k,
         \]
         \item a left integral $\lambda\in (U\ltimes H)^*$ defined on the monomial basis by
         \begin{multline*}
              \lambda(E^b F^c K^a\mb{K}^\mb{v}\mb{X}^\mb{r}) := \\
             \begin{cases}
                 \frac{\sqrt{r''}[r'-1]!}{\{1\}^{r'-1}} \ \text{\, if \,} \mb{v} = (0, ..., 0), \mb{r} = (1,...,1), a=b=c=r'-1 \\
                 0 \text{\, otherwise},
             \end{cases}
         \end{multline*}
        \item an R-matrix $R:=DR_{\matfont{z}}\Theta\Bar{R}_{\pmb{\alpha}}$, where 
        \[
            \Bar{R}_{\pmb{\alpha}}:= \exp\left( \sum^{t_2}_{l=0} \alpha_l(Z^+ \otimes \Bar{L}Z^- - Z^- \otimes \Bar{L}Z^+) \right)
        \]
        for $\Bar{L} = K^{r'/2} L$,
         \item a ribbon element 
         \begin{multline*}
             v:=v_U u_H = \\
             \frac{1-i}{\sqrt{r'}} \sum^{r'-1}_{a=0}\sum^{r''-1}_{b=0} \frac{\{-1\}^a}{[a]!} q^{-\frac{(a+3)a}{2} + 2b^2} E^aF^a  K^{-a-2b-1} L\exp \left( -2 \sum^{t_2}_{k = 1} \alpha_l Z^+_l Z^-_l  \right)
         \end{multline*}
         with the pivotal element $g= K^{-1}$.
     \end{enumerate}
\end{theorem}

\begin{corollary}
    Let $u_q \mathfrak{sl}_2 \ltimes H(\mb{m}, t, \matfont{d}, \matfont{u})$ be as above. If $H(\mb{m}, t, \matfont{d}, \matfont{u})$ carries a strongly non-factorizable ribbon structure, then it induces a non-factorizable ribbon structure on  $u_q \mathfrak{sl}_2 \ltimes H(\mb{m}, t, \matfont{d}, \matfont{u})$.
\end{corollary}
% \begin{proof}
%     By the reasoning of Proposition \ref{ExAreSNF}, if $H = H(\mb{m}, t, \matfont{d}, \matfont{u})$ is strongly non-factorizable, the top element $\mb{X}^{(1, ..., 1)}$ is missing from the monodromy matrix $M_H \in H\otimes H$. Since all new relations are also diagonal, it will be missing from the new monodromy matrix $M \in (U \ltimes H)\otimes(U \ltimes H)$, so by the same argument the result follows.
% \end{proof}

\subsection{Examples}
As in the preceding section, we restrict to $r\equiv 0 \mod 8$, $r\geq 8$ so we retain a non-trivial relation of $[E, F]$ in the $u_q \mathfrak{sl}_2$. We start with an example related to the Hopf algebra $\operatorname{SF}_2$ (the version involving $\operatorname{SF}_{2n}$ is analogous), but one which generically admits non-strongly non-factorizable quaistriangular structures.

\begin{example}
\label{SFextEx}
    Let $u_q \mathfrak{sl}_2\ltimes H$, where $H = \operatorname{SF}_2$ be a Hopf algebra generated by $K, E, F$ with relations and morphisms of Definition \ref{smallQG}, as well as $L$, $Z^\pm$ with the following relations
    \[
        L^2= \1, \;\;(Z^\pm)^2 = 0, \;\; 
    \]
    \[
        K Z^\pm K^{-1} =q^{r'}Z^\pm = -Z^\pm, \;\; L Z^{\pm} K^{-1}_2 =- Z^\pm,
    \]
    \[
        [L, K] = [L, E] = [L, F]=0,
    \]
    \[
        EZ^\pm = q^{r'} Z^\pm E= - Z^\pm E, \;\; [Z^\pm, F]=0. 
    \]
    Let $\Bar{L}:=K^{r''} L$ for convenience. The Hopf structure is 
    \[
        \epsilon(L) = 1, \;\; \epsilon(Z^\pm) = 0, 
    \]
    \[
        \Delta(L) = L \otimes L, \;\;  \Delta(Z^\pm) = \1 \otimes Z^\pm + Z^\pm \otimes \Bar{L}^{\pm 1},
    \]
    and
    \[
        S(L) = L^{-1}, \;\; S(Z^\pm) = - Z^\pm \Bar{L}^{\mp 1}.
    \]
\end{example}
\noindent The smallest instance with a non-trivial commutation relation $[E, F]$ occurs at $r=8$. The dimension is then $4^3\times 2^3  = 2^{8}=512$.
\begin{proposition}
\label{SFextProp}
    The algebra $u_q \mathfrak{sl}_2 \ltimes \operatorname{SF}_2$ is 
    \begin{enumerate}
        \item unimodular with a two-sided cointegral
        \[
            \Lambda:=\frac{\{1\}^{r'-1}}{8\sqrt{r''}[r'-1]!}\sum^{r'-1}_{a=0}\sum^2_{b}   E^{r'-1} F^{r'-1}K^a L^b Z^+ Z^-,
        \]
        and a left integral expressed on the monomial basis by 
        \[
            \lambda(  E^e F^f K^a L^b (Z^+)^g (Z^-)^h) :=
            \\ \frac{\sqrt{r''}[r'-1]!}{\{1\}^{r'-1}} \delta_{a, r'-1}\delta_{b, 0}\delta_{e, r'-1}\delta_{f, r'-1}\delta_{g, 1}\delta_{h, 1},
        \]
        \item quasitriangular, with the R-matrix $R:=R_{\matfont{z}}D\Theta \Bar{R}_{\pmb{\alpha}}$, where $D, \Theta$ were defined in Proposition \ref{uqsl2prop},  $R_{\matfont{z}}$ in Proposition \ref{SFprop}, and
        \[
            \Bar{R}_{\pmb{\alpha}} := \exp\left( \alpha( Z^+ \otimes \Bar{L} Z^- - Z^- \otimes \Bar{L} Z^+)\right)
        \]
        for $\Bar{L} = K^{r''}L$,
        \item ribbon, with the ribbon element
        \[
            v := \frac{1-i}{\sqrt{r'}} \sum^{r'-1}_{a=0}\sum^{r''-1}_{b=0} \frac{\{-1\}^a}{[a]!} q^{-\frac{(a+3)a}{2} + 2b^2}  E^aF^a K^{-a-2b-1} L\exp \left( -2 \alpha Z^+ Z^-  \right)
        \]
        corresponding to the pivotal element $g=K$.
    \end{enumerate}
\end{proposition}
\begin{remark}
\label{SF-notSNF}
    The presence of $R_{\pmb{\alpha}}$ can lead to $\lambda(M')\otimes M'' \neq 0$ being non-zero. In particular, the resulting expression is dependent on $\alpha$ and turns out to be non-zero as long as $\alpha\neq 0$. Thus, $u_q \mathfrak{sl}_2\ltimes \operatorname{SF}_2$ is not generically strongly non-factorizable.
\end{remark}
 \noindent We now construct an example admitting only strongly non-factorizable quasitriangular structures.

\begin{example}
\label{extexample1}
    Let $u_q \mathfrak{sl}_2 \ltimes \operatorname{N}_1$, be the Hopf algebra generated by $K, E, F$ with relations and morphisms of Definition \ref{smallQG}, as well as $K_1, K_2$, $X^\pm$ with the following relations
    \[
        K^4_1 = K^4_2 = \1, \;\;(X^\pm)^2 = 0, \;\; 
    \]
    \[
        K X^\pm K^{-1} =q^{r'}X^\pm = -X^\pm, \;\; K_1 X^{\pm} K^{-1}_2 = K_2 X^{\pm} K^{-1}_3 = \pm i X^\pm, 
    \]
    \[
        [K_1, K_2] = [K_1, K] = [K_1, E] = [K_1, F] = [K_2, K] = [K_2, E]=[K_2, F]=0,
    \]
    \[
        EX^\pm = q^{r'} X^\pm E= - X^\pm E, \;\; [X^\pm, F]=0. 
    \]
    Let $\Bar{L}:=K^{r''} K_1 K_2$ for convenience. The Hopf structure is 
    \[
        \epsilon(K_1) = \epsilon(K_2)=1, \;\; \epsilon(X^\pm) = 0, 
    \]
    \[
        \Delta(K_1) = K_1 \otimes K_1, \;\; \Delta(K_2) = K_2 \otimes K_2, \;\; \Delta(X^\pm) = \1 \otimes X^\pm + X^\pm \otimes \Bar{L}^{\pm 1},
    \]
    and
    \[
        S(K_1) = K_1^{-1}, \;\; S(K_2) = K_2^{-1}, \;\; S(X^\pm) = - X^\pm \Bar{L}^{\mp 1}.
    \]
\end{example}

\noindent The smallest instance with a non-trivial commutation relation $[E, F]$ occurs at $r=8$. The dimension is then $4^3\times 2^4\times2^2 = 2^{12}=4096$. We have to modify the unimodular and ribbon structure accordingly. %make precise
\begin{proposition}
\label{extProp1}
    The algebra $u_q \mathfrak{sl}_2 \ltimes \operatorname{N}_1$ is quasitriangular, with the R-matrix $R:=R_{\matfont{z}}D\Theta$, where $D, \Theta$ were defined in Proposition \ref{uqsl2prop}, and 
        \[
            R_{\matfont{z}} := \frac{1}{16}\sum_{\mb{v}, \mb{w} \in \Z^2_4} i^{-\mb{v}\mb{w}^T} (K_1, K_2)^{\mb{w}} \otimes (K_1, K_2)^{\mb{v} \matfont{z}}, 
        \] where $\matfont{z} =\begin{pmatrix} 2 && 3\\ 1 && 0 \end{pmatrix}$.
\end{proposition}
\begin{proof}
    See \cite[Proposition 4.21]{faes2025non}.
\end{proof}

It is also possible to include the exponential piece of the R-matrix as in Example \ref{ex2}, while retaining strong non-factorizability.

\begin{example}
\label{extexample2}
    Let $u_q \mathfrak{sl}_2\ltimes \operatorname{N}_2$, be the Hopf algebra generated by $K, E, F$ with relations and morphisms of Definition \ref{smallQG}, as well as $K_1, K_2, K_3$, $X^\pm$, $Z^\pm$ with the following relations, for $a=1, 2, 3$
    \[
      [K_a, K] = [K_a, E] = [K_a, F]=0,\;\;  K X^\pm K^{-1} = -X^\pm, \;\;K_a^4 = \1, \;\;  \;\; K_a X^{\pm} = \pm i X^{\pm} K_a,
    \]
    \[
         K Z^\pm K^{-1} = -Z^\pm\;\;K_1 Z^{\pm} K^{-1}_2= Z^{\pm} , \;\; K_2 Z^{\pm} K^{-1}_3= - Z^{\pm}, \;\;  K_3 Z^{\pm} K^{-1}_4= \pm i Z^{\pm},
    \]
    \[
        \{X^{\pm}, X^{\pm}\} = \{X^{\pm}, X^{\mp}\} = \{Z^{\pm}, X^{\pm}\} = \{Z^{\pm}, X^{\mp}\} = \{Z^{\pm}, Z^{\pm}\}=\{Z^{\pm}, Z^{\mp}\} =  0,
    \]
    \[
        EX^\pm = - X^\pm E,\;\; EZ^\pm = - Z^\pm E, \;\; [X^\pm, F]=[Z^\pm, F]=0. 
    \]
    Let also $\Bar{L}:=K^{r''} K_3^2$ as a shorthand. The Hopf structure is defined by
    \[
        \epsilon(K_a) = 1, \;\; \epsilon(X^{\pm})=\epsilon(Z^{\pm}) = 0,
    \]
    \[
        \Delta(K_a) = K_a \otimes K_a, \;\; \Delta(X^{\pm}) = \1 \otimes X^{\pm} + X^{\pm} \otimes (K^{r''}K_1 K_2)^{\pm 1}, \;\; \Delta(Z^{\pm}) = \1 \otimes Z^{\pm} + Z^{\pm} \otimes \Bar{L},
    \]
    \[
       S(K_a)=K^{-1}_a \;\; S(X^{\pm}) = -X^{\pm} (K^{r''} K_1 K_2)^{\mp1} \;\;  S(Z^{\pm}) = -Z^{\pm} \Bar{L}.
    \]
\end{example}
\noindent At $r=8$ the dimension is $4^4\times 2^4\times2^2 \times 2^2= 2^{16}=65536$. 
\begin{proposition}
\label{extprop2}
   The algebra $u_q \mathfrak{sl}_2\ltimes \operatorname{N}_2$ is quasitriangular, with the R-matrix $R:=R_{\matfont{z}}D\Theta R_{\pmb{\alpha}}$, where $D, \Theta$ were defined in Proposition \ref{uqsl2prop}, and 
   \begin{align*}
        R_{\matfont{z}} := \frac{1}{64}\sum_{\mb{v}, \mb{w} \in \Z^3_4} i^{-\mb{v}\mb{w}^T} (K_1, K_2, K_3)^{\mb{w}} \otimes (K_1, K_2, K_3)^{\mb{v} \matfont{z}},  &&  R_{\alpha} := \exp\left(\alpha( Z^+ \otimes \Bar{L} Z^- - Z^- \otimes \Bar{L} Z^+)\right),
   \end{align*}
 for $\matfont{z} =\begin{pmatrix} 0 && 3 && 2 \\ 1 && 0 && 0 \\ 2 && 0 && 2  \end{pmatrix}$.
        
\end{proposition}
% \begin{proof}
%     See \cite[Proposition 4.25]{faes2025non}.
%     \begin{proposition}
%     The algebra $u_q \mathfrak{sl}_2 \ltimes \operatorname{SF}_2$ is 
%     \begin{enumerate}
%         \item unimodular with a two-sided cointegral
%         \[
%             \Lambda:=\frac{\{1\}^{r'-1}}{8\sqrt{r''}[r'-1]!}\sum^{r'-1}_{a=0}\sum^2_{b}   E^{r'-1} F^{r'-1}K^a L^b Z^+ Z^-,
%         \]
%         and a left integral expressed on the monomial basis by 
%         \[
%             \lambda(  E^e F^f K^a L^b (Z^+)^g (Z^-)^h) :=
%             \\ \frac{\sqrt{r''}[r'-1]!}{\{1\}^{r'-1}} \delta_{a, r'-1}\delta_{b, 0}\delta_{e, r'-1}\delta_{f, r'-1}\delta_{g, 1}\delta_{h, 1},
%         \]
%         \item quasitriangular, with the R-matrix $R:=R_{\matfont{z}}D\Theta \Bar{R}_{\pmb{\alpha}}$, where $D, \Theta$ were defined in Proposition \ref{uqsl2prop},  $R_{\matfont{z}}$ in Proposition \ref{SFprop}, and
%         \[
%             \Bar{R}_{\pmb{\alpha}} := \exp\left( \alpha( Z^+ \otimes \Bar{L} Z^- - Z^- \otimes \Bar{L} Z^+)\right)
%         \]
%         for $\Bar{L} = K^{r''}L$,
%         \item ribbon, with the ribbon element
%         \[
%             v = \frac{1-i}{\sqrt{r'}} \sum^{r'-1}_{a=0}\sum^{r''-1}_{b=0} \frac{\{-1\}^a}{[a]!} q^{-\frac{(a+3)a}{2} + 2b^2}  E^aF^a K^{-a-2b-1} L\exp \left( -2 \alpha Z^+ Z^-  \right)
%         \]
%         corresponding to the pivotal element $g=K$.
%     \end{enumerate}
% \end{proposition}
%\end{proof}

\section{Representation theory of the Nenciu family}
In this section we consider some aspects of the representation theory of the algebras constructed with Nenciu's approach. We find the decompositions of the regular and adjoint representations into indecomposable modules, their fusion rules and braiding in the category $H\operatorname{-mod}$.

\subsection{Regular representation of $H(\mb{m}, t, \matfont{d}, \matfont{u})$}
First of all, we describe a choice of a full system of idempotents for $H(\mb{m}, t, \matfont{d}, \matfont{u})$.
\begin{definition}
\label{subIdems}
    Let for a grouplike generator $K_a\in  H(\mb{m}, t, \matfont{d}, \matfont{u})$, and the corresponding root of unity $\xi_a$, 
    $
        \varphi^a_p := \frac{1}{m_a}\sum^{m_a}_{j=0} \xi_a^{-jp}K^j_a.
    $
\end{definition}
\begin{proposition}
\label{propIdems}
    For all $p, q \in \Z_{m_a}$ and fixed $a=1, \dots, s$;  $\varphi^a_p$ are orthogonal idempotents, that is $\varphi^a_p \varphi^a_q = \delta_{pq} \varphi^a_p$.
\end{proposition}
\begin{proof}
    This is well-known, see for instance \cite{cibils_1993}, \cite{yang_2004}, or \cite{curtis_reiner_1966} (Equation 33.8).
\end{proof}

 \noindent We will now use this fact to prove the same for the full system of idempotents we are about to construct.
\begin{definition}
    Let for any $\mb{p} \in \Z_{\mb{m}}$,
    $
        \omega_{\mb{p}} := \prod^s_{a=1}  \varphi^a_{p_a},
    $
\end{definition}
\begin{proposition}
\label{systemIndems}
    For all $\mb{p}\in \Z_\mb{m}$, $\omega_{\mb{p}}$ form a full system of idempotents in $H(\mb{m}, t, \matfont{d}, \matfont{u})$. 
\end{proposition}
\begin{proof}
     As we mentioned, it is a classical fact, that $\varphi^a_{p_a}$ form full systems of primitve idempotents for all $a=1, \dots, s$ and $p_a=1, \dots, m_a$ for the group subalgebras $\C[\Z_{m_a}] \subset \C[\Z_\mb{m}]$. Now, it follows from \cite[Theorem 1]{Janusz66}, that any product $\omega_\mb{p} \in \C[\Z_\mb{m}]$ is primitive as well. It can be checked directly that $\omega_\mb{p} \omega_\mb{q}=\prod^s_{a=1} \delta_{p_a, q_a}$, and $\sum_{\mb{p}\in \Z_\mb{m}}\omega_\mb{p} = \1$.
\end{proof}

\noindent We can easily find the regular action of the generators $K_a$ on the idempotents. 
\begin{proposition}
\label{eigenprop}
    Let $\omega_{\mb{p}}$ be an idempotent as above. Then for each $a=1, \dots, s$, we have
    $
        K_a \omega_{\mb{p}} = \xi_a^{-p_a} \omega_{\mb{p}},
    $
    so that all $\omega_\mb{p}$ are eigenvectors of all $K_a$.
\end{proposition}
\begin{proof}
Recall that $\omega_{\mb{p}} = \prod^s_{a=1}  \varphi^a_{p_a}.$
For each $a, b=1, ..., s$ the multiplication by $K_a$ affects the $\varphi^b_{p_b}$ as 
    $
        K_a \varphi^b_{p_b} = \xi^{-p_a}_a \varphi^a_{p_a}
    $
    if $a=b$, 
    $
        K_a \varphi^b_{p_b} =  \varphi^b_{p_b} K_a
    $
    otherwise. Since the product contains exactly one $\varphi^a_{p_a}$, the result follows.
    % \[
    %     \omega_{\mb{p}} = \frac{1}{\prod^s_{a=1}  m_a}\sum_{\mb{j} \in \Z_{\mb{m}}} \pmb{\xi}^{-\mb{j} \cdot \mb{p}} \mb{K}^{\mb{j}},
    % \]
    % this is equivalent to the shift $\mb{j}\mapsto (j_1, \dots, j_{a-1}, j_{a}+1, \dots, j_s)$ for each $\mb{j} \in \Z_{\mb{m}}$ that yields the required constant. 
\end{proof}
\begin{notation}
Because of the above, we will somewhat abusively call $\mb{p}$ the \textit{eigenvalue} of the idempotent, as it determines the eigenvalues of the actions of all $K_a$. If we wish to distinguish it from the classical notion, we will call it the \textit{tuple eigenvalue} as apposed to a \textit{scalar eigenvalue}.
\end{notation}

\begin{definition}
\label{NenciuProj}
    Let $H=H(\mb{m}, t, \matfont{d}, \matfont{u})$ be a Nenciu type Hopf algebra, and $\{\omega_\mb{p}\}_{\mb{p}\in \Z_\mb{m}}$, its system of idempotents. Define $P_\mb{p} := H\omega_\mb{p}$ to be the module of the shape

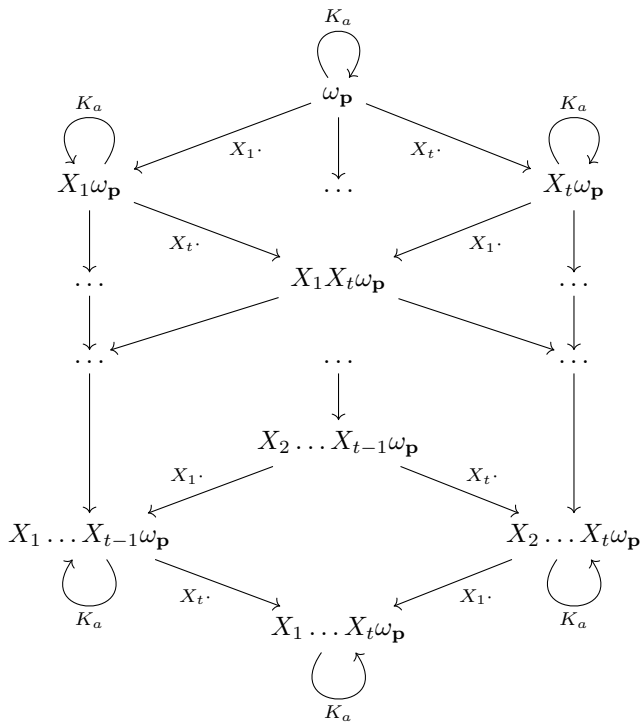
\begin{center}
% https://q.uiver.app/#q=WzAsMTYsWzEsMCwiXFxvbWVnYV9cXG1hdGhiZntwfSJdLFsxLDEsIi4uLiJdLFswLDEsIlhfMVxcb21lZ2FfXFxtYXRoYmZ7cH0iXSxbMiwxLCJYX3RcXG9tZWdhX1xcbWF0aGJme3B9Il0sWzEsMiwiWF8xIFhfdCBcXG9tZWdhX1xcbWF0aGJme3B9Il0sWzIsMiwiLi4uIl0sWzEsNCwiWF8yLi4uWF97dC0xfVxcb21lZ2FfXFxtYXRoYmZ7cH0iXSxbMSw2LCJYXzEuLi5YX3QgXFxvbWVnYV9cXG1hdGhiZntwfSJdLFsxLDUsIi4uLiJdLFsyLDUsIlhfMi4uLlhfe3R9XFxvbWVnYV9cXG1hdGhiZntwfSJdLFswLDUsIlhfMS4uLlhfe3QtMX1cXG9tZWdhX1xcbWF0aGJme3B9Il0sWzAsMiwiLi4uIl0sWzAsMywiLi4uIl0sWzIsMywiLi4uIl0sWzEsMywiLi4uIl0sWzIsNl0sWzAsMiwiWF8xXFxjZG90Il0sWzAsMywiWF90XFxjZG90IiwyXSxbMiw0LCJYX3RcXGNkb3QiLDJdLFszLDQsIlhfMSBcXGNkb3QiXSxbOCw3XSxbOSw3LCJYXzFcXGNkb3QiXSxbMTAsNywiWF90XFxjZG90IiwyXSxbMCwxXSxbMiwxMV0sWzMsNV0sWzEyLDEwXSxbNSwxM10sWzExLDEyXSxbMTMsOV0sWzYsMTAsIlhfMVxcY2RvdCIsMl0sWzYsOSwiWF90XFxjZG90Il0sWzYsOF0sWzQsMTJdLFs0LDEzXSxbMTQsNl0sWzAsMCwiS19hIl0sWzIsMiwiS19hIiwyLHsicmFkaXVzIjotMywiYW5nbGUiOi0xODB9XSxbMywzLCJLX2EiXSxbNyw3LCJLX2EiLDIseyJyYWRpdXMiOi0zfV0sWzksOSwiS19hIiwyLHsicmFkaXVzIjotM31dLFsxMCwxMCwiS19hIiwwLHsiYW5nbGUiOi0xODB9XV0=
\begin{tikzcd}
	& {\omega_\mathbf{p}} \\
	{X_1\omega_\mathbf{p}} & {\dots} & {X_t\omega_\mathbf{p}} \\
	{\dots} & {X_1 X_t \omega_\mathbf{p}} & {\dots} \\
	{\dots} & {\dots} & {\dots} \\
	& {X_2\dots X_{t-1}\omega_\mathbf{p}} \\
	{X_1\dots X_{t-1}\omega_\mathbf{p}} & {} & {X_2\dots X_{t}\omega_\mathbf{p}} \\
	& {X_1\dots X_t \omega_\mathbf{p}} & {}
	\arrow["{K_a}", from=1-2, to=1-2, loop, in=55, out=125, distance=10mm]
	\arrow["{X_1\cdot}", from=1-2, to=2-1]
	\arrow[from=1-2, to=2-2]
	\arrow["{X_t\cdot}"', from=1-2, to=2-3]
	\arrow["{K_a}"', from=2-1, to=2-1, loop, in=125, out=55, distance=10mm]
	\arrow[from=2-1, to=3-1]
	\arrow["{X_t\cdot}"', from=2-1, to=3-2]
	\arrow["{K_a}", from=2-3, to=2-3, loop, in=55, out=125, distance=10mm]
	\arrow["{X_1 \cdot}", from=2-3, to=3-2]
	\arrow[from=2-3, to=3-3]
	\arrow[from=3-1, to=4-1]
	\arrow[from=3-2, to=4-1]
	\arrow[from=3-2, to=4-3]
	\arrow[from=3-3, to=4-3]
	\arrow[from=4-1, to=6-1]
	\arrow[from=4-2, to=5-2]
	\arrow[from=4-3, to=6-3]
	\arrow["{X_1\cdot}"', from=5-2, to=6-1]
	
	\arrow["{X_t\cdot}", from=5-2, to=6-3]
	\arrow["{K_a}", from=6-1, to=6-1, loop, in=235, out=305, distance=10mm]
	\arrow["{X_t\cdot}"', from=6-1, to=7-2]
	
	\arrow["{K_a}"', from=6-3, to=6-3, loop, in=305, out=235, distance=10mm]
	\arrow["{X_1\cdot}", from=6-3, to=7-2]
	\arrow["{K_a}"', from=7-2, to=7-2, loop, in=305, out=235, distance=10mm]
\end{tikzcd}
\end{center}
generated by a single element $\omega_\mb{p}$. The loops corresponding to the actions of the $K_a$ are present for each element and for all $a=1, ..., s$, but are not represented for readability. The quotient module generated by the actions of grouplike generators on the single element $\omega_{\mb{p}}$, is the head. Also, $\mathrm{Soc}(P_\mb{p}) = \langle X_1 \dots X_t\omega_\mb{p} \rangle$, is the 1-dimensional submodule comprising the bottom node of the diagram.
\end{definition}
\noindent Now, we can state the main theorem of this section.
\begin{theorem}
\label{regDecomp}
    Let $H =  H(\mb{m}, t, \matfont{d}, \matfont{u})$, be a Nenciu type Hopf algebra with the usual data. Then the regular representation $(H, \cdot)$ decomposes as
    \[
        (H, \cdot) \cong \bigoplus_{\mb{p} \in \Z_{\mb{m}}} P_{\mb{p}}
    \]
    where each $P_{\mb{p}} = H \omega_{\mb{p}}$ is a $2^t$-dimensional projective module, generated by acting on the fixed idempotent with the nilpotent generators. 
\end{theorem}
\begin{proof}
    We established in the the Proposition \ref{regDecomposition} that we can decompose $H$ into the direct sum of indecomposable protective modules using the full system of idempotents $\{\omega_\mb{p}\}_{\mb{p}\in \Z_\mb{m}} \subset H$.  According to Definition \ref{NenciuProj} each $P_{\mb{p}}$ is a module generated by a primitive idempotent $\omega_\mb{p}$. Indeed, there are $\prod^s_{a=1}m_a$ idempotents so the total dimension is $2^{t}\prod^s_{a=1}m_a $ as required. 
\end{proof}

\begin{remark}
\label{reorderingRem}
    It should be remarked that the diagram in Definition \ref{NenciuProj} is schematic and does reflect the basis of a $P_\mb{p}$, but not the exact action on each element. Once an ordering of the generators $X_k$ in the monomial basis of $H(\mb{m}, t, \matfont{d}, \matfont{u})$ is fixed, the subsequent actions of various $X_k$ on $\omega_\mb{p}$ will incur commutations constants accordingly. These are the following
    $
        X_k \mb{X}^\mb{r}\omega_\mb{p} = \pmb{\xi}^{\matfont{u}_k\cdot \sum^{k}_{j=1}r_j \matfont{d}_j}\mb{X}^{\Bar{\mb{r}}}\omega_\mb{p}
    $
    where $\Bar{\mb{r}}$ is such that $\Bar{r}_k = r_k + 1$, and $\Bar{r}_j = r_j$ for $j\neq k$ if $r_k=0$, otherwise the expression vanishes.
\end{remark}
We also describe the corresponding simple modules in $H\operatorname{-mod}$.
\begin{definition}
    Let $H:=H(\mb{m}, t, \matfont{d}, \matfont{u})$ and $P_\mb{p} = H\omega_\mb{p}$ as above. Then let 
    $
        S_\mb{p} = P_\mb{p}/\text{Rad }P_\mb{p},
    $
    be the corresponding simple submodule.
\end{definition}
\begin{proposition}
    All simple modules $S_{\mb{p}}$ are 1-dimensional and, in particular, invertible.
\end{proposition}
\begin{proof}
    The fact that all $S_\mb{p}$ are 1-dimensional follows from the diagonality of relations and because all $X_k$, $k=1, \dots, t$ are nilpotent of order 2. This implies they are all invertible, and to see this explicitly, fix $S_{\mb{p}}$ for some $\mb{p}\in \Z_{\mb{m}}$. Let also $q = (-p_1, ..., -p_a)$. Then it is easy to check that $S_\mb{p}\otimes S_\mb{q} = S_\mb{q}\otimes S_\mb{p}=S_{\mb{0}}$ as required.
\end{proof}
\begin{example}
\label{SFProjectives}
    Consider the symplectic fermions $\operatorname{SF}_{2n}$ of Definition \ref{SFdef} with $2n$ nilpotent generators. The full system of idempotents is just $\omega_0 = \frac{\1+L}{2}$, $\omega_1 = \frac{\1-L}{2}$. The regular representation decomposes as
    $
        (\operatorname{SF}_{2n}, \cdot) = P_0 \oplus P_1
    $
    into two $2^{2n}$-dimensional projectives $P_0 = \operatorname{SF}_{2n} \omega_0$ and $P_1 = \operatorname{SF}_{2n} \omega_1$ with the scalar eigenvalues $1$ and $-1$, respectively. 
\end{example}

\begin{example}
    Consider $\operatorname{N}_2$ from Example \ref{ex2} at $t_1 = t_2  = 1$. Then $\xi_1 = \xi_2 = \xi_3 = i$, so the full system of idempotents is $\{\omega_{\mb{p}}\}_{\mb{p} \in \Z^3_4}$. There are four nilpotent generators $X^\pm$ and $Z^\pm$. Thus we have $4^3 = 64$ projective modules $P_{\mb{p}}$ of dimension $2^4=16$, decomposing the 1024-dimensional Hopf algebra, each with a different eigenvalue tuple
    $
        (\operatorname{N}_2, \cdot) = \oplus_{\mb{p} \in \Z^3_4} P_{\mb{p}}.
    $
\end{example}

\noindent We finish off this section by introducing certain class of indecomposable modules over $ H(\mb{m}, t, \matfont{d}, \matfont{u})$ that are submodules of $P_\mb{p}$ under the regular action. These will appear when we consider $(H, \rhd)$, the decomposition of $H$ under the adjoint action. We can use a tuple $\mb{r} \in \Z^t_2$ to express a subset of the generators $X_1, \dots, X_t$. Thus we adapt the following notation.

\begin{definition}
\label{indecomps}
Let $\mb{p} \in \Z_{\mb{m}}$ and $\mb{r} = (r_1, \dots, r_t)\in \Z^t_2$. Define $V_{\mb{p}}^{\mb{r}}$ to be a submodule of $P_\mb{p}$ generated by the element $\mb{X}^\mb{r}\omega_\mb{p}$. Then if $r_k \neq 0$ for $k=1, \dots, t$, then $X_k V_\mb{p}^{\mb{r}}=0$ under the regular action. In particular $P_\mb{p} = V^\mb{0}_\mb{p}$ and $S_\mb{p} = V^{(1, \dots, 1)}_\mb{p}$.\\
For instance, for $\mb{r}=(0, 1, 1, \dots, 1, 0)$ the shape of $V^\mb{r}_\mb{p}$ (red) inside a $P_\mb{p}$ is
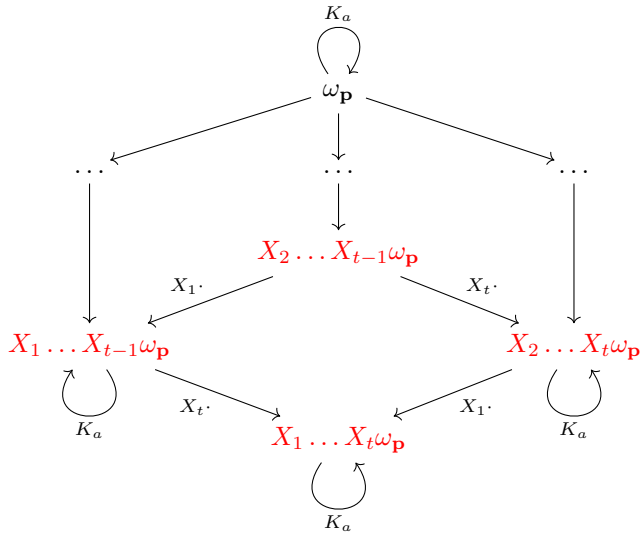
\begin{center}
   % https://q.uiver.app/#q=WzAsOSxbMSwwLCJcXG9tZWdhX1xcbWF0aGJme3B9Il0sWzEsMiwiXFx0ZXh0Y29sb3J7cmVkfXtYXzIuLi5YX3t0LTF9XFxvbWVnYV9cXG1hdGhiZntwfX0iXSxbMSw0LCJcXHRleHRjb2xvcntyZWR9e1hfMS4uLlhfdCBcXG9tZWdhX1xcbWF0aGJme3B9fSJdLFsyLDMsIlxcdGV4dGNvbG9ye3JlZH17WF8yLi4uWF97dH1cXG9tZWdhX1xcbWF0aGJme3B9fSJdLFswLDMsIlxcdGV4dGNvbG9ye3JlZH17WF8xLi4uWF97dC0xfVxcb21lZ2FfXFxtYXRoYmZ7cH19Il0sWzAsMSwiLi4uIl0sWzIsMSwiLi4uIl0sWzEsMSwiLi4uIl0sWzIsNF0sWzMsMiwiWF8xXFxjZG90Il0sWzQsMiwiWF90XFxjZG90IiwyXSxbNSw0XSxbNiwzXSxbMSw0LCJYXzFcXGNkb3QiLDJdLFsxLDMsIlhfdFxcY2RvdCJdLFs3LDFdLFswLDAsIktfYSJdLFsyLDIsIktfYSIsMix7InJhZGl1cyI6LTN9XSxbMywzLCJLX2EiLDIseyJyYWRpdXMiOi0zfV0sWzQsNCwiS19hIiwwLHsiYW5nbGUiOi0xODB9XSxbMCw1XSxbMCw2XSxbMCw3XV0=
\begin{tikzcd}
	& {\omega_\mathbf{p}} \\
	{\dots} & {\dots} & {\dots} \\
	& {\textcolor{red}{X_2\dots X_{t-1}\omega_\mathbf{p}}} \\
	{\textcolor{red}{X_1 X_2\dots X_{t-1}\omega_\mathbf{p}}} && {\textcolor{red}{X_2\dots X_{t-1} X_t\omega_\mathbf{p}}} \\
	& {\textcolor{red}{X_1\dots X_{t}  \omega_\mathbf{p}}} & {}
	\arrow["{K_a}", from=1-2, to=1-2, loop, in=55, out=125, distance=10mm]
	\arrow[from=1-2, to=2-1]
	\arrow[from=1-2, to=2-2]
	\arrow[from=1-2, to=2-3]
	\arrow[from=2-1, to=4-1]
	\arrow[from=2-2, to=3-2]
	\arrow[from=2-3, to=4-3]
	\arrow["{X_1\cdot}"', from=3-2, to=4-1]
	\arrow["{X_t\cdot}", from=3-2, to=4-3]
	\arrow["{K_a}", from=4-1, to=4-1, loop, in=235, out=305, distance=10mm]
	\arrow["{X_t\cdot}"', from=4-1, to=5-2]
	\arrow["{K_a}"', from=4-3, to=4-3, loop, in=305, out=235, distance=10mm]
	\arrow["{X_1\cdot}", from=4-3, to=5-2]
	\arrow["{K_a}"', from=5-2, to=5-2, loop, in=305, out=235, distance=10mm]
\end{tikzcd}
\end{center}
The loops from $K_a$ actions are present for each element, not represented for readability. The quotient module generated by the single element $g:=\mb{X}^\mb{r}\omega_{\mb{p}} \in V_{\mb{p}}^{\mb{r}}\subset P_\mb{p}$ under the action of grouplike generators is the head.  Also, note that $\mathrm{Soc}(V^\mb{r}_\mb{p}) =  \langle X_1 \dots X_t\omega_\mb{p}\rangle$ is the same as $\mathrm{Soc}(P_\mb{p})$.
\end{definition}
\noindent Intuitively, we have already acted with these $X_{r_1},\dots,X_{r_s}$ inside $P_\mb{p}$, and what is left is a submodule, where the remaining generators act. It is worth remarking, that these form just a small subset of indecomposable modules in $H\operatorname{-mod}$, if $t>1$. More questions regarding their classification will be addressed in later sections. 
\begin{proposition}
    For any $\mb{p}\in \Z_\mb{m}$ and $\mb{r}\in \Z^t_2$, $V_{\mb{p}}^{\mb{r}}$ is an indecomposable module.
\end{proposition}
\noindent The following proof was suggested to the author by M. Faitg.
\begin{proof}
We will show that any submodule of $V_{\mb{p}}^{\mb{r}}$ intersects the $\text{Soc}(V_\mb{p}^\mb{r})$, and thus $V_{\mb{p}}^{\mb{r}}$ cannot be written as a direct sum.\\
Consider any element $v\in V_{\mb{p}}^{\mb{r}}$, and action of all $X_k$ such that $r_k=0$ on $v$ to generate a submodule $M\subset V_{\mb{p}}^{\mb{r}}$. We can write $v$ in the monomial basis of Proposition \ref{monomialBasis} as
    \[
        v = \sum_{\mb{s} \in \Z^t_2|\mb{X}^\mb{s}\mb{X}^\mb{r} \neq 0} \gamma_{\mb{s}} \mb{X}^\mb{s} \mb{X}^\mb{r} \omega_\mb{p},
    \]
    for $\gamma_\mb{s} \in \C$. By diagonality of relations, and linear independence of the monomial basis, we can act by all $X_k$ such that $r_k\neq0$, collected in $\mb{X}^\mb{o}$, for $\mb{o} = (1, ..., 1)-\mb{r}$. Then $\mb{X}^\mb{o}v$ is such that $X_k \mb{X}^\mb{o}v=0$ for all $k=1, \dots, t$. \\
    On the other hand, if we write
    \[
         X_k \mb{X}^\mb{o}v = \sum_{\mb{s} \in \Z^t_2|\mb{X}^\mb{s}\mb{X}^\mb{r} \neq 0} \gamma_{\mb{s}} X_k \mb{X}^\mb{o}\mb{X}^\mb{s} \mb{X}^\mb{r} \omega_\mb{p} = 0,
    \]
    then the sum vanishes term-wise by diagonality of relations and linear independence of the elements of the monomial basis. Thus, $\mb{X}^\mb{o}v = \gamma \mb{X}^{(1, \dots, 1)}\omega_\mb{p}$, for some $\gamma\in \C$, so $\mb{X}^\mb{o}v \in \text{Soc}(V_\mb{p}^\mb{r})$. Thus, for each $v$ generating $M$, we found a $\mb{X}^\mb{o}v
    \in M\cap \text{Soc}(V_\mb{p}^\mb{r})$, and so $M\cap \text{Soc}(V_\mb{p}^\mb{r})\neq 0$, as we wished.
\end{proof}

%\textcolor{red}{Complete using Zaitsev-Nikolenko for SF. We always have a copy of SF action for even-order groups, by restricting the action of Ks to $\Z_2$ subgroup (check if we can have a canonical K easily).}
%\begin{proposition}
% \label{submodProp}
%     Any submodule of $V^\mb{r}_\mb{p}$ has the form $V^\mb{r'}_\mb{p}$, where $\mb{r}\cdot \mb{r}' = \mb{r}$. In particular 
%     $\mathrm{Soc}(V^\mb{r}_\mb{p}) = \mathrm{Soc}(V^\mb{r'}_\mb{p})$, as elements of $H$.
% \end{proposition}
% \begin{proof}
%     Firstly, the relation $\mb{r}\cdot \mb{r}' = \mb{r}$ means that every element that acts non-trivially on $V^\mb{r'}_\mb{p}$, acts non-trivially also on $V^\mb{r}_\mb{p}$. We prove the claim inductively, on $|r'|$, that is the number of generators $X_k$ that act trivially on $V^\mb{r'}_\mb{p}$. \\
%     \underline{Base}: the base case is the head $g\in V^\mb{r}_\mb{r}$, which of course generates the whole module.
%     \underline{Step}: Assume that any submodule where $l$ of $X_k$ act trivially has the shape $V^\mb{r'}_\mb{p}$ for some $|r'|=l\leq t$. Any module where $l+1$ of $X_k$ act trivially will have the head $\mb{X}^\mb{r'} g$, as it is necessarily an element of $\mb{P}_\mb{p}$. Then the action of the remaining $X_k$ produces clearly a submodule $V^\mb{r''}_\mp{p}$ with $|\mb{r''}|=l+1$. We continue the induction up to when $|\mb{r}^{(l)}|=t$, that is until we reach $\mathrm{Soc}(V^\mb{r}_\mb{p})$, which clearly all the submodules share, as claimed.
% \end{proof}
\noindent We will also use the following fact.
\begin{lemma}
    If $V^\mb{r}_\mb{p}\cap V^\mb{s}_\mb{q} \neq 0$, then $P_\mb{p} = P_\mb{q}$.
\end{lemma}
\begin{proof}%Matthieu's proof
    Consider $V^\mb{r}_\mb{p}, V^\mb{s}_\mb{q}\subset A$ as submodules of the regular representation. Then since $V^\mb{r}_\mb{p}\cap V^\mb{s}_\mb{q} \neq 0$, we find $V^\mb{r}_\mb{p} \subset P_\mb{q}$ and $V^\mb{s}_\mb{q} \in P_\mb{p}$. But recall that $P_\mb{p}, P_\mb{q}$ are direct summands of $A$, and we have a non-trivial module $V^\mb{r}_\mb{p}\cup V^\mb{s}_\mb{q} \subset P_\mb{p}$ and $V^\mb{r}_\mb{p}\cup V^\mb{s}_\mb{q} \in P_\mb{q}$. Hence, it must be $P_\mb{p} = P_\mb{q}$. 
\end{proof}
\noindent We also have an obvious corollary
\begin{corollary}
\label{submodCor}
    If $V^\mb{r}_\mb{p}\cap V^\mb{s}_\mb{q} \neq 0$, then $\mathrm{Soc}(V^\mb{r}_\mb{p}) = \mathrm{Soc}(V^\mb{s}_\mb{q})$. 
\end{corollary}
\subsection{Adjoint representation}
In this section we consider the decomposition of a Nenciu type algebra $H= H(\mb{m}, t, \matfont{d}, \matfont{u})$ under the adjoint action $(H, \rhd)$, as introduced in Proposition \ref{adjointAct}. 
\begin{notation}
    For any nilpotent generator $X_k \in H$, denote $L_k:= \mb{K}^{\matfont{u}_k}$. 
\end{notation}
\noindent We note the following useful fact.
\begin{proposition}
\label{XLcommutes}
    Let $X_k, X_l$ be nilpotent generators, then
    $
        [X_l, X_k L_k^{-1}] =0.
    $
\end{proposition}
\begin{proof}
    We have directly
    $
        [X_l, X_k L^{-1}] = X_l X_k L^{-1} - X_k L^{-1} X_l = X_l X_k L^{-1} - \xi^{(\matfont{u}_k - \matfont{u}_k)\cdot \matfont{d}_l}X_l X_k L^{-1} = 0.
    $
\end{proof}

\noindent As a result of the associativity of the adjoint action, for any $h\in H$, 
    $
         X_l \rhd X_k \rhd h = \pmb{\xi}^{\matfont{d}_k\cdot \matfont{u}_l}  X_k \rhd X_l \rhd h  
    $
 and we observe the following.
\begin{proposition}
\label{monomialForm}
Any monomial $\mb{X}^\mb{a} \mb{K}^\mb{v} \in H$ (note the change in ordering against Proposition \ref{monomialBasis}) can be written as 
\[
    \mb{X}^{\mb{a}} \mb{L}^{-\mb{a}}\mb{K}^\mb{w}=\gamma (X_{1}L^{-1}_{1})^{a_1}\dots (X_{t}L^{-1}_{t})^{a_t} \mb{K}^\mb{w}, 
\]
where $\mb{K}^\mb{w}=\mb{K}^{\mb{v}} L_{a_1}\dots L_{a_t}$ and $\gamma\in \C$ is a constant obtained by commuting all $L^{-1}_{a_k}$ out of the product of $X_{a_k}$. Moreover, such a decomposition is unique. Thus we have a re-expression of the monomial basis of $H$ as 
$
    \{\mb{X}^{\mb{a}} \mb{L}^{-\mb{a}}\mb{K}^\mb{w}| \mb{w} \in \Z_{\mb{m}} \text{,  } \mb{a} \in \Z^t_2\}.
$
\end{proposition}
\noindent Let us keep this notation in this section. Then, we have the following corollary.
\begin{corollary}
\label{corCommut}
    Let $x = \mb{X}^\mb{a} \mb{L}^{-\mb{a}} \mb{K}^\mb{w} \in H$ be any monomial in the Hopf algebra as in Proposition \ref{monomialForm} and $X_k \in H$ a nilpotent generator. Then $[X_k, x]=0$ if and only if:
    \begin{itemize}
        \item[1.] either $X_k x = 0$, i.e. $X_k$ is one of the $X_{a_l}$,
        \item[2.] or $[X_k, \mb{K}^\mb{w}]=0$. 
    \end{itemize}
\end{corollary}
\begin{proof}
    Clearly if $X_k x = 0$, then $[X_k, x]=0$. Similarly, if $x = \mb{X}^\mb{a} \mb{L}^{-\mb{a}} \mb{K}^\mb{w}$, then from Proposition \ref{XLcommutes} we see that any $X_k$ commutes with $\mb{X}^\mb{a}\mb{L}^{-\mb{a}}$, so we compute
    $
        [X_k, \mb{X}^\mb{a} \mb{L}^{-\mb{a}} \mb{K}^\mb{w}] = \mb{X}^\mb{a} \mb{L}^{-\mb{a}} [X_k, \mb{K}^\mb{w}]=0
    $
    as long as $[X_k, \mb{K}^\mb{w}]=0$.\\
    Conversely, assume $ [X_k, \mb{X}^\mb{a} \mb{L}^{-\mb{a}} \mb{K}^\mb{w}]=0$, then we again find 
    $
        0=[X_k, \mb{X}^\mb{a} \mb{L}^{-\mb{a}} \mb{K}^\mb{w}] = \mb{X}^\mb{a} \mb{L}^{-\mb{a}} [X_k, \mb{K}^\mb{w}]
    $
    so either $[X_k, \mb{K}^\mb{v}]=0$, or the expression contains a pair of zero divisors, that is $X_k \mb{X}^\mb
{a}= 0$.
\end{proof}
\noindent Furthermore, we have the following simple proposition with an immediate proof.
\begin{proposition}
     If $X^+_k$ and $X^-_k$ are nilpotent generators of opposite type, then if $w \in C(X^+_k)$, then $w \in C(X^-_k)$ Thus, if $X^+_k \rhd w = 0$, then $X^-_k\rhd w = 0$.
\end{proposition}
With these observations, let now $W\subset (H, \rhd)$ be a submodule under the adjoint action. Let $C(h)$ be the centraliser of the element $h\in H$. If we write $w=\mb{X}^{\mb{a}} \mb{L}^{-\mb{a}}\mb{K}^\mb{w}$, then $X_k \rhd w = 0$ if $\mb{K}^\mb{v}\in C(X_k)$ or $X_k$ appears in $\mb{X}^{\mb{a}}$. We can extend this idea in order to relate $(H, \rhd)$ to the previous discussion of $(H, \cdot)$.

\begin{definition}
    Denote by $W^\mb{a}_{\mb{K}^\mb{w}} \subset (H, \rhd)$ a submodule of $(H, \rhd)$ generated by a single element $g:=\mb{X}^{\mb{a}} \mb{L}^{-\mb{a}}\mb{K}^\mb{w}$ under adjoint action.
\end{definition}
\noindent Then we have the following correspondence.
\begin{proposition}
\label{submodCor}
    Let $W^\mb{a}_{\mb{K}^\mb{w}}$ be a $H$-module under the adjoint action, generated by the element $g:=\mb{X}^{\mb{a}} \mb{L}^{-\mb{a}} \mb{K}^\mb{v} \in H$. Then $W^\mb{a}_{\mb{K}^\mb{w}}\cong V^{\mb{r}}_{\mb{p}}$, where $\mb{r} \in \Z^t_2$ is such that whenever $r_k\neq0$ then $X_{k}\in C(g)$, and $\mb{p}, \mb{q}\in \Z_\mb{m}$ are such that $K_a \mb{X}^\mb{r} = \xi^{q_a}\mb{X}^\mb{r} K_a$, and $K_a \mb{X}^\mb{a} = \xi_a^{-p_a+q_a}\mb{X}^\mb{a} K_a$.
\end{proposition}
\begin{proof}
    We exhibit the isomorphism on the monomial bases of both modules. First of all, $\mb{q}$ and $\mb{r}$ are determined by $\mb{K}^\mb{w}$. We choose $\mb{r}$ so that the regular action of $X_k$ on $\mb{X}^\mb{r}$ be non-zero exactly when the adjoint action of $X_k$ on $g$ is non-zero. Then $\mb{q}$ is determined by this choice and $\mb{p}$ can be chosen accordingly to fix the eigenvalue. We consider the map defined on the monomial bases (given the modules are submodules of $(H, \lhd)$ and $(H, \cdot)$ respectively) by 
    \begin{align*}
    f: W^\mb{a}_{\mb{K}^\mb{w}} &\rightarrow V^{\mb{r}}_{\mb{p}}, \\
         \mb{X}^{\mb{a}} \mb{L}^{-\mb{a}}\mb{K}^\mb{w} &\mapsto \mb{X}^\mb{r}\omega_\mb{p}, \\
          [X_k, \mb{X}^{\mb{a}} \mb{L}^{-\mb{a}}\mb{K}^\mb{w}] L^{-1}_k &\mapsto X_k \mb{X}^\mb{r} \omega_\mb{p}
    \end{align*}
    with the inverse map obtained by simply inverting all the arrows. We assume $f$ and its inverse are $\C$-linear. By the fact that all the relations are diagonal and all $X_k^2=0$, and by the choice of $\mb{r}\in \Z^t_2$ we can easily see the assignment is both injective and surjective. In particular any element $w\in W^\mb{a}_{\mb{K}^\mb{w}}$ can be written as a linear combination of elements of the form $\mb{X}^\mb{b} \rhd g$, for some $\mb{b} \in \Z^t_2$. By linearity, it is enough to establish that the map sends $g \mapsto \mb{X}^\mb{r}\omega_\mb{p}$ and then proceed by induction on $|\mb{b}|$, the number of actions by nilpotent generators performed away from the generator for each monomial.\\
    For the induction base we have to verify that for any $h = K_a, X_k \in H$, 
    $
        f(h\rhd \mb{X}^{\mb{a}} \mb{L}^{-\mb{a}} \mb{K}^\mb{w}) = h \mb{X}^\mb{r} \omega_\mb{p}. 
    $
    \noindent By linearity for any $a=1, \dots, s$
    \[
        f(K_a \rhd \mb{X}^{\mb{a}} \mb{L}^{-\mb{a}}\mb{K}^\mb{w}) = f(K_a \mb{X}^{\mb{a}} \mb{L}^{-\mb{a}}\mb{K}^\mb{w} K^{-1}_a) = \xi^{-p_a + q_a}_a f(\mb{X}^{\mb{a}} \mb{L}^{-\mb{a}}\mb{K}^\mb{w})=\xi^{-p_a + q_a}_a \mb{X}^{\mb{r}}\omega_\mb{p}, 
    \]
    but we know that
    \[
         K_a \mb{X}^{\mb{r}}\omega_{\mb{p}} = \xi^{-p_a + q_a}_a \mb{X}^{\mb{r}}\omega_{\mb{p}}.
    \]
    Similarly if $r_k \neq0$
    \[
        X_k \rhd \mb{X}^{\mb{a}} \mb{L}^{-\mb{a}}\mb{K}^\mb{w}= X_k  \mb{X}^\mb{r} \omega_\mb{p} = 0.
    \]
    However, if $r_k = 0$, by construction
    \[
        f(X_k \rhd \mb{X}^{\mb{a}} \mb{L}^{-\mb{a}}\mb{K}^\mb{w}) = f([X_k, \mb{X}^{\mb{a}} \mb{L}^{-\mb{a}}\mb{K}^\mb{w}] L^{-1}_k) = X_k \mb{X}^\mb{r}\omega_\mb{p} = X_k f(\mb{X}^{\mb{a}} \mb{L}^{-\mb{a}}\mb{K}^\mb{w}).
    \]
    For the induction step, assume now $w \in W^\mb{a}_{\mb{K}^\mb{w}}$ is of the form $\mb{X}^\mb{b} \rhd g$ for some $|\mb{b}|=b<t$, so that $f(w) = \mb{X}^\mb{b} \mb{X}^\mb{r} \omega_\mb{p}$. By construction we have that
    \[
        f(K_a \rhd w) 
        % = \xi^{- l_a}_a f(\mb{X}^\mb{b} \rhd K_a \rhd w) = \xi^{- l_a}_a \mb{X}^b K_a \mb{X}^\mb{r}\omega_\mb{p}
         = K_a f(w) = \xi^{-p_a + q_a - l_a}_a  f(w),
    \]
    for $\mb{l} \in \Z_\mb{m}$ such that $K_a \mb{X}^\mb{b} = \xi^{- l_a}_a\mb{X}^\mb{b}K_a$, so this action is compatible with $f$. 
    Then for the action of $X_k$ we assume that
    \[
        f(X_k \rhd w) = f([X_k, w] L^{-1}_k) = X_k f(w).
    \]
    Now, we consider the element $X_k \rhd w$. 
    If $K_b$ is another grouplike generator we have 
    \[
        f(K_b\rhd X_k \rhd w) = \xi^{d_{kb}} X_k f(K_b\rhd w) = \xi^{d_{kb}} X_k K_b f(w) = K_b X_k f(w), 
    \]
    by associativity of both actions.\\
    Also, if $X_l$ is another nilpotent generator, by associativity of the adjoint action in Hopf algebras, we have 
    \[
        f(X_lX_k  \rhd w) = f(X_l \rhd X_k \rhd w) = X_l X_k f(w) = \pmb{\xi}^{\matfont{u}_l \cdot \matfont{d}_k} X_k X_l f(w) = \pmb{\xi}^{\matfont{u}_l \cdot \matfont{d}_k}f(X_k X_l \rhd w).
    \]

    This means that after we constructed $w$ with $|\mb{b}|$ steps away from the generator, the steps $|\mb{b}|+1$ and $|\mb{b}|+2$ can be applied and their order of multiplication is respected by the map.
    The argument for the inverse statement that is analogous.
\end{proof}
\begin{corollary}
     For the setup of Proposition \ref{submodCor} if $C(g) \cap \rad H = 0$, then $W^\mb{\mb{0}}_{\mb{K}^\mb{w}} \cong P_{\mb{0}}$.
\end{corollary}
\begin{proof}
    If $C(g) \cap \rad H = 0$, and $\mb{a}=\mb{0}$, then $W^\mb{\mb{0}}_{\mb{K}^\mb{w}} \cong P_{\mb{0}}$, as every action $X_k\rhd$ is non-vanishing, and the eigenvalue measured at the generator is $\mb{0}$ because $K_a$ all commute.
\end{proof}
% \begin{remark}
%     For any $w \in W^\mb{a}_{\mb{K}^\mb{w}}$ we have in general 
%     \[
%         (X_l X_k)\rhd w \neq [X_l X_k, w]L^{-1}_k L^{-1}_l. 
%     \]
%     Nevertheless, $f$  respects the associativity of the action on both sides as 
%     \[
%         [X_l, [X_k, w]L^{-1}_k]L^{-1}_l = (X_l X_k)_{(1)} w S((X_l X_k)_{(2)}) = (X_l X_k) \rhd w
%     \]
%     for which we established
%     \[
%         f((X_l X_k) \rhd w) = X_l X_k f(w) = X_l(X_k f(w)).
%     \]
% \end{remark}
Hence, the notation $W^\mb{a}_{\mb{K}^\mb{w}}$ is actually redundant, and we will continue using $V^\mb{r}_\mb{p}$ instead, but still sometimes invoke the equivalence of the regular and adjoint picture. Before stating the decomposition of $(H, \rhd)$ into indecomposables, we need another piece of notation.
\begin{notation}
\begin{enumerate}
    \item Let $\mathcal{L}:=\{\mb{(\mb{r}}\matfont{d}) \cdot \mb{v}| \mb{v} \in \Z_{\mb{m}}, \mb{r} \in \Z^t_2\}$. That is, $\{\pmb{\xi}^{\mb{l}}, \mb{l} \in \mathcal{L}\}$ is the set of all constants obtainable by commuting any $\mb{K}^{\mb{v}}$ past any $\mb{X}^{\mb{r}}$, with multiplicities.
   
    \item Let $\mathcal{L}_{\mb{s}}:=\{((\mb{r}\cdot \mb{s}) \matfont{d})\cdot \mb{v}| \mb{v} \in \Z_{\mb{m}}, \mb{r} \in \Z^t_2\}$ for a fixed $\mb{s}\in \Z^t_2$. This means we choose only such elements of $\mathcal{L}$ that can be obtained given the generators participating in $\mb{X}^\mb{s}$.
    \item Let $\mb{r}(\mb{w})$ be the vector in $\Z^t_2$ selecting all of nilpotent generators $X_k\in C(\mb{K}^\mb{w})$. In particular, $\mathcal{L}_{\mb{r}(\mb{0})} = \mathcal{L}_{\mb{0}}= \mathcal{L}$.
    
\end{enumerate}
\end{notation}

\begin{theorem}
\label{adjointDecomp}
     Let $(H, \rhd)$ be the adjoint representation of a Nenciu type Hopf algebra $H =  H(\mb{m}, t, \matfont{d}, \matfont{u})$. Then we have the following decomposition into indecomposable modules
     \[
         (H, \rhd)   \cong \bigoplus_{\mb{v} \in \Z_\mb{m}} \bigoplus_{\mb{l}\in \mathcal{L}_{\mb{r}(\mb{v})}} V^{\mb{r}(\mb{w})}_{-\mb{l}}.
     \]
     In particular any $\mb{K}^\mb{w} \notin \bigcup^t_{k=1} C(X_k)$ contributes a single projective $P_{\mb{0}}$ and any $\mb{K}^\mb{w} \in Z(H)$ contributes $2^{t}$ simples $S_{-\mb{l}}$ for all $\mb{l} \in \mathcal{L}$.
\end{theorem}
\begin{proof}
    We use the adjoint picture for the proof, so that we consider any $V^\mb{r(v)}_\mb{-l}$ as the corresponding $W^\mb{a}_{\mb{K}^\mb{w}}$. Consider a generating element $g = \mb{K}^\mb{w}$. We treat three cases separately.
    \begin{enumerate}
        \item If $g \notin \bigcup^t_{k=1} C(X_k)$, consider a submodule $H \rhd g$. Then we can perform all successive actions $X_k \rhd g $ for $k-1, ..., t$ to produce the $2^t$-dimensional projective module $P_{\mb{p}}$ for some $\mb{p} \in \Z_\mb{m}$. Now, since all nilpotent generators are used to build the module, the eigenvalue from acting by any $K_a$ (measured at the generator) is the constant resulting from the commutation of $K_a$, for $a=1, ..., s$ with the grouplike $g$, so $\mb{p} = \mb{0}$, as all elements of $\C[\mb{K}]$ commute, and we retrieve $P_\mb{0}$.
        
        \item If $g \in \bigcap^t_{k=1} C(X_k) = Z(H)$, then $g$ generates a trivial module, as no non-zero action $X_k \rhd g$ can be performed without killing all elements, while the actions of all $K_a$ result in the eigenvalue $\mb{0}$. There are also monomials $\mb{X}^\mb{b} \mb{L}^{-\mb{b}}g \in H$ that could be written as $\mb{X}^\mb{b} \mb{L}^\mb{-b}g$, where $\mb{b}\in \Z^t_2$. By the choice of $g$ this action is always zero and none of the monomials can be constructed by acting on $g$, and moreover $X_k \rhd \mb{X}^\mb{b} \rhd g=0$ for any $X_k$. Thus, the entire vector subspace underlying the module $P_{\mb{0}}$ (but not with the action of this module) that would be build over $g$ in item 1. is split into a direct sum  $2^t$ 1-dimensional simple modules $S_{\mb{-l}}$, each with one of the such $\mb{X}^\mb{b} \mb{L}^\mb{-b}g$ as the generating element, such that $K_a \rhd \mb{X}^\mb{b} \mb{L}^\mb{-b}g = \xi^{-l_a}_a$.  The available eigenvalues are given by all $\mb{l} \in \mathcal{L}$, since any $\mb{X}^\mb{b} \mb{L}^\mb{-b}g$ is allowed as a generator. Note that this allows for repeated eigenvalues, and multiplicities of the modules in the decomposition are easy to determine for fixed Nenciu data, as demonstrated in examples below.
        
        \item We are left with the intermediate case, that is $g \in \bigcap_{k\in \kappa} C(X_k)$ for some indexing set $\kappa \subsetneq \{1, ..., t\}$. Moreover, let the nilpotent generators acting trivially $g$ be collected according to $\mb{r}(\mb{w})$, or equivalently $\mb{a}$ in the adjoint picture, and consider elements of the form $h=\mb{X}^\mb{r} \mb{L}^\mb{-r} \mb{K}^\mb{w}$, such that if $r_k\neq 0$, then $(\mb{r}(\mb{w}))_k\neq 0$. To each such $\mb{r}$, there corresponds an $\mb{l}\in \mathcal{L}_{\mb{r}(\mb{w})}$, so each $h$ generates a submodule $V^\mb{r(v)}_{-\mb{l}}$.  By the reasoning of part 2. restricted to this subset we take the corresponding space underlying a module $P_\mb{0}$ that could be built over $\mb{K}^\mb{v}$. By construction, none of $\mb{X}^\mb{r} \mb{L}^\mb{-r} \mb{K}^\mb{w}$ can be reached from $\mb{K^\mb{w}}$ by the action of any $X_k$, or from another such element. Hence, each produces a generator $h$ for a disjoint  $V^\mb{r(w)}_{-\mb{l}}$, exhausting all elements of this $P_\mb{0}$. Obviously, none of the elements can belong to a projective or simple of items 1. and 2., by the uniqueness of this form of the generator monomials.
        \end{enumerate}
    Since in each case, for each grouplike we have a $P_{\mb{0}}$-worth of elements, with each copy of the projective clearly disjoint, we have in total $2^{t+s}$ elements, and we decomposed the entire module $ (H, \rhd)$.
\end{proof}
\noindent With this proof in mind, we introduce a handy piece of notation.
\begin{notation}
    For $h\in H$, let $C_{\mathcal{G}}(h):= C(h) \cap \C[\mb{K}]$ indicate the centraliser of $h$ restricted to $\C[\mb{K}]$, so containing only elements of the group subalgebra $\C[\mb{K}]$ that commute with $h$.
\end{notation}

\begin{example}
    Consider again $\operatorname{SF}_{2n}$. There are only two grouplikes that could serve as generators of adjoint submodules, $\1$ and $K$. The former, $\1,$ produces $2^{2n-1}$ simple modules $S_{0}$ and $2^{2n-1}$ simple modules $S_{1}$, each built on an even or odd element, respectively. The latter, $K$, gives the projective $P_{0}$, whose basis elements are monomials containing an odd number of $Z_k$ generators so in total
    \[
        (\operatorname{SF}_{2n}, \rhd) \cong 2^{2n-1} S_{0} \oplus 2^{2n-1} S_{1} \oplus  P_{0}.
    \]
\end{example}
\begin{example}
    Consider $\operatorname{N}_2$ of Example \ref{ex2} at $t_1=t_2=1$. For any $h\in H$, let $C_{\mathcal{G}}(h):= C(h) \cap \C[K_1, K_2, K_3]$. Let also $\mb{r}$ be ordered as $(X^+, X^-, Z^+, Z^-)$. Then for $\mb{w} \in \Z^3_4$ we have the following centralisers.
    \begin{itemize}
        \item  For $X^\pm$
            \[
            C_{\mathcal{G}}(X^\pm) = \{\mb{K}^{\mb{w}}| w_1+w_2 +w_3 \equiv 0 \mod 4\}
            \]
    with $|C_{\mathcal{G}}(X^\pm)|=12$. The base elements for each $g \in C_{\mathcal{G}}(X^\pm)$ are $g, X^\pm \mb{K}^{\pm \matfont{u}_1} g, X^+X^-g$, so the corresponding eigenvalues are \\$\mb{l} \in\{(0, 0, 0), (\pm 1, \pm 1, \pm1)\} $. This gives $4 \cdot 12 = 48$ modules of dimension 4, so in total
    \[
        24 V^{(1, 1, 0, 0)}_{(0, 0, 0)} \oplus 12 V^{(1, 1, 0, 0)}_{(1, 1, 1)} \oplus 12 V^{(1, 1, 0, 0)}_{(-1, -1, -1)}.
    \]
    \item  For $Z^\pm$ 
    \[
        C_{\mathcal{G}}(Z^\pm) = \{\mb{K}^{\mb{w}}| 2w_2 +w_3 \equiv 0 \mod 4\}
    \]
    with $|C_{\mathcal{G}}(Z^\pm)|=12$. The available generators for each $g \in C_{\mathcal{G}}(Z^\pm)$ are $g, Z^\pm L g, Z^+Z^- g$, so the corresponding eigenvalues are $\mb{l} \in\{(0, 0, 0), (0, 2, \pm  1)\} $. This gives $4 \cdot 12 = 48$ modules of dimension 4, so in total
    \[
        24 V^{(0, 0, 1, 1)}_{(0, 0, 0)} \oplus 12 V^{(0, 0, 1, 1)}_{(0, 2, 1)} \oplus 12 V^{(0, 0, 1, 1)}_{(0, 2, -1)}.
    \]
    \item For both $X^\pm, Z^\pm$
    \[
         C_{\mathcal{G}}(X^\pm) \cap C_{\mathcal{G}}(Z^\pm) =\{\1, K_1 K_2 K^2_3, K^2_1 K_2^2, K^3_1 K^3_2 K^2_3\}
    \]
    which gives $16 \cdot 4 = 64$ trivial modules. The generating elements cover any combination of nilpotent generators for any $g$, so the eigenvalues are $\mb{l}\in\{(0, 0, 0), (\pm 1, \pm 1, \pm 1), (0, 2, \pm  1), (\pm1, \mp 1, 2), (\pm 1, \mp 1, 0)\}$. In total we have
    \[
        16 S_{(0, 0, 0)} \oplus 8 S_{(1, 1, 1)} \oplus 8 S_{(-1, -1, -1)} \oplus 8 S_{(1, -1, 2)} \oplus 8 S_{(-1, 1, 2)} \oplus 8 S_{(1, -1, 0)} \oplus 8 S_{(-1, 1, 0)}.
    \]
    \item The remaining $64-24 = 36$ grouplikes which belong to neither centraliser give $36$ dimension $16$ projective modules, all with the eigenvalue $\mb{0}$,
    $
        36 P_{\mb{0}}.
    $
    \end{itemize}
   
\end{example}
\subsection{Fusion rules for $V^\mb{r}_\mb{p}$}
Now that we reduced both decompositions of interest to enumeration of modules $V^{\mb{r}}_{\mb{p}}$, we can find their fusion rules under the regular action. For this, we keep in mind Theorem 3.1 of \cite{yang_2004}, which describes the rules for the case $s=1$. 
\begin{theorem}
\label{fusionRulesNenciu}
    Let $V^{\mb{r}}_{\mb{p}}$ and $V^{\mb{s}}_{\mb{q}}$ be indecomposable $H$-modules, where $H = H(\mb{m}, t, \matfont{d}, \matfont{u})$. Then 
    \[
        V^{\mb{r}}_{\mb{p}} \otimes V^{\mb{s}}_{\mb{q}} \cong \bigoplus_{\mb{l}\in \mathcal{L}_{(1, \dots, 1)-(\mb{r}+ \mb{s})-\mb{r}\cdot \mb{s}}} V^{\mb{r}\cdot\mb{s}}_{\mb{p}+\mb{q} -\mb{a}-\mb{b}-\mb{l}},
    \]
    where $g=\mb{X}^\mb{r}\omega_\mb{p}$, $h=\mb{X}^\mb{s}\omega_\mb{q}$ are the generators of $V^\mb{r}_\mb{p}$ and $V^\mb{s}_\mb{q}$ respectively, and $\mb{a}, \mb{b} \in \Z_\mb{m}$ are such that $K_a \mb{X}^\mb{r} = \xi^{a_a} \mb{X}^\mb{r} K_a$, and $K_a \mb{X}^\mb{s} = \xi^{b_a} \mb{X}^\mb{s} K_a$.
\end{theorem}
\begin{proof}
    We first explain how the eigenvalues are measured in the direct summands. Let $\mb{c}\in \Z_\mb{m}$ be such that $K_a \mb{X}^{\mb{r}\cdot \mb{s}} = \xi_a^{c_a} \mb{X}^{\mb{r}\cdot \mb{s}} K_a$. Note first that for any $K_a$
        \[
            \Delta(K_a) (g\otimes h) = K_a g \otimes K_a h = 
            \xi_a^{-p_a -q_a +a_a+b_a} g\otimes h,
        \]
    since both $\mb{X}^\mb{a}$ and $\mb{X}^\mb{b}$ contribute to the constant in front. Now, it is clear that if $r_k \neq 0$ and $s_k \neq 0$, then $\Delta(X_k)(g\otimes h)=0$. Otherwise, 
    \[
            \Delta(X_k) (g \otimes h) =g \otimes X_k h + X_k g \otimes L_k h \neq 0.
    \]
    Thus, we wish to relate this module generated by $(g\otimes h)$ to a module of the shape $V^{\mb{r}\cdot \mb{s}}_\mb{p'}$ for some $\mb{p}'$. Here we have 
    \[
        K_a \mb{X}^{\mb{r}\cdot \mb{s}} \omega_\mb{p'} = \xi_a^{c_a - p'_a} \mb{X}^{\mb{r}\cdot \mb{s}} \omega_\mb{p'}.
    \]
    Now, we need to have $\mb{p}' = \mb{p}+\mb{q}$, but this does not conform to the constant $\xi_a^{-p_a -q_a +a_a+b_a}$. Thus, we shift the eigenvalue of the underlying idempotent generator to $\mb{p} + \mb{q} - \mb{a} - \mb{b}$, resulting in the eigenvalues of the generators $g\otimes h$ and $\mb{X}^{\mb{r}\cdot \mb{s}}\omega_{\mb{p}+\mb{q}-\mb{a}-\mb{b}}$ being the same.\\
    
    Now, we show the direct sum decomposition of $V^{\mb{r}}_{\mb{p}} \otimes V^{\mb{s}}_{\mb{q}}$, by finding the corresponding generators and showing the submodules do not intersect. Let $\mb{o} = (1, \dots, 1) - (\mb{r}+\mb{s})-\mb{r}\cdot \mb{s}$ be the tuple selecting all such $X_k$ that $X_k g \neq 0$ and $X_k h \neq 0$. Let $\mb{l}\in \mathcal{L}_{\mb{o}}$, and let $Y = \mb{X}^\mb{o}$ be constructed using such $X_k$, so that we have $K_a \mb{X}^\mb{o} = \xi_a^{l_a} \mb{X}^\mb{o} K_a$ for $\mb{o}\in \Z^t_2$. Then a generator $(Y g) \otimes h$ gives rise to a submodule $V^\mb{\mb{r}\cdot\mb{s}}_{\mb{p}+\mb{q}-\mb{a}-\mb{b} - \mb{l}}$, since any of the $X_k$ such that $(\mb{r}\cdot \mb{s})_k = 0$ acts non-trivially, and the eigenvalue of the idempotent is obtained as discussed above. \\
    
    There are exactly $2^{t-|\mb{r}+\mb{s}| - |\mb{r}\cdot \mb{s}|}$ such $Y$ and since each submodule will have dimension $2^{t-|\mb{r}\cdot\mb{s}|}$, their total dimension is $2^{(t-|\mb{r}|)(t-|\mb{s}|)}$, the same as the dimension of $ V^{\mb{r}}_{\mb{p}} \otimes V^{\mb{s}}_{\mb{q}}$. This is because $|\mb{r}|+|\mb{s}| = |\mb{r}+\mb{s}|+2|\mb{r}\cdot \mb{s}|$, as the entries of the tuples are in the set $\{0, 1\}$. \\
    Therefore, we have left to show that for distinct $Y_1$, $Y_2$ the resulting modules $V^{\mb{r}\cdot \mb{s}}_{\mb{p}+\mb{q}-\mb{a}-\mb{b} - \mb{l}}$ and $V^{\mb{r}\cdot \mb{s}}_{\mb{p}+\mb{q}-\mb{a}-\mb{b} - \mb{l}'}$ are disjoint, which will then imply their sum is direct. Let $Y_\alpha'\in H$, for $\alpha=1, 2$ be a monomial composed of the remaining $X_k$ that act non-trivially on $g\otimes h$. These include such $X_k$ that $r_k =  0$ or $s_k= 0$. Then,
    \[
        \Delta(Y'_\alpha)((Y_\alpha g) \otimes h )= \mathrm{Soc}(V^\mb{\mb{r}\cdot\mb{s}}_{\mb{p}+\mb{q}-\mb{a}-\mb{b} - \mb{l}_\alpha})
    \]
    as modules. Since the modules $ \mathrm{Soc}(V^\mb{\mb{r}\cdot\mb{s}}_{\mb{p}+\mb{q}-\mb{a}-\mb{b}-\mb{l}_\alpha})$ are 1-dimensional, if $Y_1\neq Y_2$ then 
    \[
        \Delta(Y'_1)((Y_1 g) \otimes h)\neq \Delta(Y'_2)((Y_2 g) \otimes h).
    \] 
    This is because there exists at least one $X_l$ such that $r_l =  0$ and $s_l= 0$, so at least one action $\Delta(X_l) ((Y_\alpha g)\otimes h)$ results in both terms after the action of $\Delta(X_l)$ being non-zero, and thus providing linearly independent terms for $\alpha=1, 2$. Then\[
    \mathrm{Soc}(V^\mb{\mb{r}\cdot\mb{s}}_{\mb{p}+\mb{q}-\mb{a}-\mb{b} - \mb{l}_1}) \neq \mathrm{Soc}(V^\mb{\mb{r}\cdot\mb{s}}_{\mb{p}+\mb{q}-\mb{a}-\mb{b} - \mb{l}_2}).
    \] 
    Thus, by the contrapositive of Corollary \ref{submodCor}, we have that 
    \[
        V^\mb{\mb{r}\cdot\mb{s}}_{\mb{p}+\mb{q}-\mb{a}-\mb{b} - \mb{l}_1} \cap V^\mb{\mb{r}\cdot\mb{s}}_{\mb{p}+\mb{q}-\mb{a}-\mb{b} - \mb{l}_2}=0,
    \]
    and so the decomposition of the tensor product into submodules is direct.
\end{proof}

\begin{example}
    Consider $\operatorname{SF}_{2n}$, and its two indecomposable projectives $P_0, P_1$.  The fusion rules between them are
    \[
        P_0 \otimes P_0 = P_1 \otimes P_1 = 2^{2n-1} P_0 \oplus 2^{2n-1} P_1, 
    \]
    \[
        P_0 \otimes P_1 = 2^{2n} P_1, 
    \]
    where the difference in the final rule comes from the fact that while all nilpotent elements commute with $K^0 = K^{1+1} = \1$, only even elements commute with $K^1 = K$. So $\mathcal{L}_{\mb{r}(1)} = \{0, \dots, 0\}$ with multiplicity $2^{2n}$, as opposed to $\mathcal{L}_{\mb{r}(0)} = \{0, 1, \dots, 0, 1\}$ both with half the multiplicity.
\end{example}
\noindent We consider also another example.
\begin{example}
Consider $\operatorname{N}_2$ of Example \ref{ex2} at $t_1=t_2=1$, and two projectives $P_{(0, 0, 0)}, P_{(1, 0, 0)} \in \operatorname{N}_2\operatorname{-mod}$. Since the resulting module is projective any element of $\mathcal{L}$ is realised in the modules $V^\mb{0}_{-\mb{l}}$ of the decomposition.  Recall from \cite[Example 3.11]{faes2025non} that the commutation data for this Hopf algebra is $\matfont{d}_{1} = (1, 1, 1)$, $\matfont{d}_{2} = (-1, -1, -1)$ and $\matfont{d}_{3} = (0, 2, 1)$, $\matfont{d}_{4} = (0, -2, -1)$, so we can parametrise 
$$\mathcal{L} = \{(a-b,a-b+2c-2d, a-b+c-d),|a, b, c, d \in \Z_2\},$$
where $a, b, c, d$ count the presence of the generators $X^+, X^-, Z^+, Z^-$, respectively. This is in total $2^4=16$ (possibly pairwise isomorphic, but always disjoint) submodules of $\operatorname{N}_2$, with each resulting in a projective module of dimension 16. Thus, we retrieve the correct dimensions, and the decomposition stands 
\[
    P_{(0, 0, 0)} \otimes P_{(1, 0, 0)} = \bigoplus_{a, b, c, d\in\Z_2} P_{(1-a+b,-a+b-2c+2d, -a+b-c+d)}.
\]
Consider now a pair of submodules of the indecomposable projectives, say $V^\mb{r}_\mb{p} = V^{(1, 1, 1, 0)}_{(0, 0, 0)}$ and $V^\mb{s}_\mb{q} =V^{(1, 0, 1, 0)}_{(1, 0, 0)}$. We have $\mb{r}\cdot \mb{s} = (1, 0, 1, 0)$, so this will be the shape of the summand modules in the decomposition. Moreover, $\mb{a} = (0, 2, 1)$ and $\mb{b} = (1, 3, 2)$. We now take $\mathcal{L}_{(1, 1, 1, 1) - (\mb{r}+\mb{s})-\mb{r}\cdot \mb{s}} = \mathcal{L}_{(0, 0, 0, 1)}$, which is just $\mathcal{L}_{(0, 0, 0, 1)}=\{(0, 0, 0), (0, 2, 3)\}$. Thus, the decomposition is
\[
    V^{(1, 1, 1, 0)}_{(0, 0, 0)} \otimes V^{(1, 0, 1, 0)}_{(1, 0, 0)} = V^{(1, 0, 1, 0)}_{(0, 0, 0)} \oplus V^{(1, 0, 1, 0)}_{(0, 1, 2)}. 
\]
The dimensions on both sides can be checked to be the consistent $2\cdot 4 = 8$. 
\end{example}
% \begin{remark}
%     \textcolor{red}{The reason acting once never produces a zero terms is that the resulting terms are all linearly independent - different order of action gives different value - this is essentially the acting by the product and then indeed no cancellations can occur. \\
%     The reason acting twice always gives zero, is firstly by the fact the adjoint it a Hopf action so is consistent with the commutation in the acting term. Secondly, because terms of shape $(\1 \otimes X_k)$ and $(X_k \otimes L_k)$ (resp. $X_l, L_l$)all commute with one another the same way -  just as $X_k X_l$ would - no matter the ordering of action, the commutation is consistent. Then the two $(\1 \otimes X_k)(X_k \otimes L_k) = -(X_k \otimes L_k)(\1 \otimes X_k)$, to produce the vanishing cancellation.\\
%     The modules in the decomposition are the same, because for a chosen fusion rule, the same nilpotents act on everything, but the starting terms are different. Since they are all linearly independent at the start and always the entire module is produced, they will be linearly independent, so disjoint.}
% \end{remark}

\subsection{Some braiding properties of $V^{\mb{r}}_{\mb{p}}$.}
In this section, we investigate the braiding of the modules $V^{\mb{r}}_{\mb{p}}$ for a given Nenciu type Hopf algebra $ H(\mb{m}, t, \matfont{d}, \matfont{u})$, where we assume the following:
\begin{enumerate}
    \item The R-matrix has the shape $R_{\matfont{z}}R_{\pmb{\alpha}}$ where $R_{\matfont{z}}$ satisfies $R_{\matfont{z}} R_{\matfont{z}, 21} = \1 \otimes \1$ and is constructed with only grouplike generators, and $R_{\pmb{\alpha}}$ has the shape 
    \[
        R_{\pmb{\alpha}} = \exp\left(\sum^{t_2}_{l=1} \alpha_l (Z_l^+ \otimes L_l Z_l^{-} - Z_l^- \otimes L_l Z_l^{+})\right)
    \]
    where $\alpha_l \in \C$ and $Z_l^\pm$ are of opposite type for $l_1, l_2=1, \dots, t_2$. Moreover $\{Z_{l_1}^\pm, Z_{l_2}^\pm\}=\{Z_{l_1}^\pm, Z_{l_2}^\mp\}=0$.
    \item From Proposition \ref{monodromyProp} the monodromy matrix is
    \[
        M:=R_{21}R =\exp\left(2\sum^{t_2}_{l=1} \alpha_l (Z_l^+ \otimes L_l Z_l^{-} - Z_l^- \otimes L_l Z_l^{+})\right).
    \] 
    \item There is a (possibly empty) subset of nilpotent generators, $X_k$, $k=1, \dots, t_1$ that do not participate in $R_{\pmb{\alpha}}$ and such that $\{X_k, Z^\pm_{l}\}=0$. 
\end{enumerate}
With these assumptions, we have the following theorem.
\begin{theorem}
\label{symBraiding}
    Let $(H, R)$ be a Nenciu type, quasitriangular Hopf algebra and $V^{\mb{r}}_\mb{p}, V^{\mb{s}}_\mb{q} \in H\operatorname{-mod}$. Then
    \[
        c^2(V^{\mb{r}}_\mb{p} \otimes V^{\mb{s}}_\mb{q}) =  (V^{\mb{r}}_\mb{p} \otimes V^{\mb{s}}_\mb{q}), 
    \]
    if and only if for each entry $Z^\pm_l$, either $\alpha_l=0$ or the corresponding $r_l \neq 0$ or $s_l \neq 0$. Otherwise, the order of the braiding map $c$ is infinite.
\end{theorem}
\begin{proof}
    Assume first that for each entry $Z^\pm_l$, either $\alpha_l=0$ or the corresponding $r_l \neq 0$ or $s_l \neq 0$. This means that either a particular pair of generators $Z^\pm_l$ is excluded from the $R_{\pmb{\alpha}}$ by setting the parameter $\alpha_l=0$, or choosing both $Z^\pm_l$ to act trivially on either $V^{\mb{r}}_\mb{p}$ or $V^{\mb{s}}_\mb{q}$. It is enough that one of the modules has zero action for each term $Z_l^\pm \otimes L Z_l^\mp$. 
    % These two combined mean that the only piece of $R_{\pmb{\alpha}}$ that acts on the module with a non-zero value is $\1\otimes \1$, so this is equivalent to the action of $R_{\matfont{z}}$ in total. But this satisfies $R_{\matfont{z}} R_{\matfont{z}, 21} = \1\otimes \1$, so the result is $c^2 =\Id\otimes \Id$.  
    Note that $c^2$ is equivalent to the action of $M$ by associativity of the action, so then
    \[
        M' V^{\mb{r}}_\mb{p} \otimes  M''V^{\mb{s}}_\mb{q} =\1 V^{\mb{r}}_\mb{p} \otimes  \1 V^{\mb{s}}_\mb{q} = V^{\mb{r}}_\mb{p} \otimes  V^{\mb{s}}_\mb{q},
    \]
    as $R_\mb{z}R_{\mb{z}, 21} =\1 \otimes \1$.\\
    Coversely, assume now that $ c^2(V^{\mb{r}}_\mb{p} \otimes V^{\mb{s}}_\mb{q}) =  (V^{\mb{r}}_\mb{p}\otimes V^{\mb{s}}_\mb{q})$. By above, equivalently
    \[
        M' V^{\mb{r}}_\mb{p} \otimes  M''V^{\mb{s}}_\mb{q} = V^{\mb{r}}_\mb{p} \otimes  V^{\mb{s}}_\mb{q}.
    \]
    We wish to show this is not true unless the condition is fulfilled. It is enough to consider the product of the generators $g \otimes h$, then 
    \[
        M' g \otimes M'' h = g \otimes h + 2\sum^{t_2}_{l=1} \alpha_l (Z_l^+ \otimes L_l Z_l^{-} - Z_l^- \otimes L_l Z_l^{+}) (g \otimes h) +\mathcal{O}((Z_l^{\pm})^2),
    \]
    where the ''big O'' notation indicates terms at least quadratic in the generators $Z^\pm_l$. By extending the monomial basis of $H$ to $H\otimes H$ in the obvious way, we see that since all relations are diagonal, the action of any skew-primitive $X_k$ or $Z^\pm_l$ cannot result in distinct terms in the sum becoming linearly dependent. Thus, the sum does not vanish, unless for each $(Z_l^\pm \otimes L_l Z_l^\mp)$, either $\alpha_l=0$, or both $Z_l^\pm$ act trivially on the respective modules. In the latter case, it is necessary for both $Z^\pm_l$ so that $(Z_l^\pm \otimes L_l Z_l^\mp)$ and $(Z_l^\mp \otimes L_l Z_l^\pm)$ vanish. Hence, we retrieve the desired condition. \\\\
    For the final statement, note that by associativity for $2m=0, 2, \dots$
    \[
        c^{2m}(V^{\mb{r}}_\mb{p} \otimes V^{\mb{s}}_\mb{q}) = (M^{m})' V^{\mb{r}}_\mb{p} \otimes (M^{m})''V^{\mb{s}}_\mb{q} 
    \]
     with $M^0:= \1 \otimes \1$, and
    \[
        c^{2m+1}(V^{\mb{r}}_\mb{p} \otimes V^{\mb{s}}_\mb{q}) = R_{\mb{z}, 21} R_{\pmb{\alpha}, 21}\left((M^{m})'' V^{\mb{s}}_\mb{q} \otimes (M^{m})' V^{\mb{r}}_\mb{p}\right). 
    \]
    Since $m$ is even, clearly this is the action of 
    \begin{align*}
        \exp\left(2m\sum^{t_2}_{l=1} \alpha_l (Z_l^+ \otimes L_l Z_l^{-} - Z_l^- \otimes L_l Z_l^{+}) \right) =
         \1\otimes \1 + 2m\sum^{t_2}_{l=1} \alpha_l (Z_l^+ \otimes L_l Z_l^{-} - Z_l^- \otimes L_l Z_l^{+}) + 
        \mathcal{O}((Z_l^{\pm})^2),
    \end{align*}
    \noindent and as above the sum vanishes if and only if the condition of the theorem is fulfilled. In the odd case, it is enough to verify that for any $2m=0, 2, \dots $, the $2m+1$st iteration is non-zero. But we are considering the following action on the generator $h \otimes g$
    \begin{align*}
         R_{\mb{z}, 21} &\exp\left(\sum^{t_2}_{l=1} \alpha_l (L_lZ_l^- \otimes Z_l^{+} - L_lZ_l^+ \otimes Z_l^-)\right)
         \exp\left(2m\sum^{t_2}_{l=1} \alpha_l (L_lZ_l^- \otimes Z_l^{+} - L_lZ_l^+ \otimes Z_l^-)\right) (h \otimes g)=\\
     =R_{\mb{z}, 21}&\exp\left((2m+1)\sum^{t_2}_{l=1} \alpha_l (L_lZ_l^- \otimes Z_l^{+} - L_lZ_l^+ \otimes Z_l^-)\right)(h \otimes g)\\
     =R_{\mb{z}, 21}&\left( h \otimes g + (2m+1)\sum^{t_2}_{l=1} \alpha_l (L_lZ_l^- \otimes Z_l^{+} - L_lZ_l^+ \otimes Z_l^-) (h \otimes g)
     +\mathcal{O}((Z_l^{\pm})^2)\right),
    \end{align*}
    where the sum is again non-vanishing, unless the condition is fulfilled, by the same reasoning. This completes the proof.
\end{proof}
\noindent Using our running examples we see the following.
\begin{example}
    Consider $\operatorname{SF}_{4n}$. Fix also all entries of $\pmb{\alpha}$ to be non-zero. Then, considering the projectives $P_0$, $P_1$, and their submodules, the only transparent modules are the simples $S_0$, $S_1$.  \\
    However, a pair of modules can braid trivially with each other but not other modules, such as the $P_0$ and $P_1$. Consider $V^\mb{r}_0$, with $r_l = 1$ if $l\leq n$, that is one where $Z^+_1, Z^-_1,\dots,Z^+_n, Z^-_n$ act trivially, but $Z^+_{n+1}, Z^-_{n+1},\dots,Z^+_{2n}, Z^-_{2n}$ act non-trivially. Consider also  $V^\mb{s}_0$, with $s_l = 1$ if $n<l\leq 2n$, that is one where $Z^+_1, Z^-_1,\dots,Z^+_n, Z^-_n$ act non-trivially, but $Z^+_{n+1}, Z^-_{n+1},\dots,Z^+_{2n}, Z^-_{2n}$ act trivially. Then, for all $l=1, \dots, 2n$: 
    \[
         (Z^\pm_l\otimes L Z_l^\mp)(V^\mb{r}_0 \otimes V^\mb{s}_0)=0
    \]
    so 
    \[
        \exp\left(\sum^{2n}_{l=1} \alpha_l (Z_l^+ \otimes L_l Z_l^{-} - Z_l^- \otimes L_l Z_l^{+})\right)  (V^\mb{r}_0 \otimes V^\mb{s}_0)= (\1\otimes \1)(V^\mb{r}_0 \otimes V^\mb{s}_0)
    \]
    and so
    \[
        c^2(V^\mb{r}_0 \otimes V^\mb{s}_0) =(V^\mb{r}_0 \otimes V^\mb{s}_0).
    \]
    We also note that if $\pmb{\alpha} = \mb{0}$ then $\operatorname{SF}_{4n}\operatorname{-mod}$ is a symmetric category.
\end{example}
\begin{example}
\label{N2ExNSSMod}
    Consider now $\operatorname{N}_2$. Then, all the projectives $P_\mb{p}$ braid non-trivially with all $V^\mb{s}_\mb{q}$. Similarly, all simples $S_\mb{p}$ braid symmetrically with all $V^\mb{s}_\mb{q}$. But take now $V^\mb{r}_\mb{p}$, for $\mb{r} = (0, 0, 1, 1)$ then any $Z^+$, $Z^-$ act trivially on $V^\mb{r}_\mb{p}$. Then for any $V \in \operatorname{N}_2\operatorname{-mod}$
    \[
        c^2(V^\mb{r}_\mb{p} \otimes V) = (V^\mb{r}_\mb{p} \otimes V). 
    \]
    These $V^\mb{r}_\mb{p}$ are therefore transparent, so the M\"uger centre $\mathcal{Z}_2(\operatorname{N}_2\operatorname{-mod})$ is non-trivial (in fact non-semisimple), which gives another proof that $\operatorname{N}_2$ is non-factorizable, by Defintion \ref{factorizability} 3.
\end{example}

\noindent From these examples we can deduce that following about $\mathcal{Z}_2(\operatorname{N}_2\operatorname{-mod})$.
\begin{remark}
\label{catEmbeddingProp}
    There is a faithfully full braided monoidal functor
    $
        F: \operatorname{SF}_{2t_1} \operatorname{-mod} \rightarrow \mathcal{Z}_{(2)} (\operatorname{N}_{2} \operatorname{-mod}),  
    $
    where $\operatorname{SF}_{2 t_2}$ carries the braiding induced by the R-matrix $R_{\matfont{z}}$ of Proposition \ref{SFprop}, that is in the case $\pmb{\alpha} = \mb{0}$. 
\end{remark}

\subsection{Further indecomposable modules}
Building on the above remark, the projective modules and their submodules are by no means the only indecomposable modules of $H(\mb{m}, t, \matfont{u}, \matfont{d})$. For instance, the indecomposable modules of $\operatorname{SF}_{2}$, which can be seen as the bosonization of the exterior algebra on two generators $\bigwedge V_2$, for a 2-dimensional $V_2 \in \text{SVect}_{\C}$ (a Hopf algebra in the category of super-vector spaces $\text{SVect}_{\C}$) are fully classified. This is a consequence of the result of \cite{zaitsev_nikolenko_1971}, described in more modern terms in \cite{runkel_2014}. While these authors work with the unbosonized algebra exterior algebra $\bigwedge V_2$ directly, the adaptation to the bosonized case is straightforward, by making even elements of the modules to be acted upon by $L$ with the eigenvalue $1$, and odd by $-1$. More precisely, besides the projective modules described in Example \ref{SFProjectives}, we find the following two families cf. \cite{runkel_2014}[Theorem 2.6].
\begin{definition}
\begin{enumerate}
    \item  Let $\Xi^+_{k, \varepsilon, \delta}$ be the representation of $\operatorname{SF}_{2n}$ where $k \in \Z_{\geq 0}$ and $\varepsilon, \delta \in \{0, 1\}$, and the underlying vector space is $\C^{2k +\varepsilon +\delta +1}$. Fix a basis a basis of the "even" vectors $e^{[\varepsilon]}, e^1, \dots, e^k, e^{[\delta]}$ where $e^{\varepsilon}$ appears only if $\varepsilon = 1$ and $e^{\delta}$ appears only if $\delta = 1$. Fix also a basis of "odd" vectors $o^1, \dots, o^{k+1}$. Then we have the following module structure
    \begin{align*}
         L e^\alpha = e^\alpha, \;\; \text{for } \alpha = 1, \dots, k, \varepsilon, \delta, &&  L o^\beta = -o^\beta, \;\; \text{for } \beta = 1, \dots, k+1.
    \end{align*}
    and for $\gamma=1, \dots, k$
    \begin{align*}
        Z^+ e^{\gamma} &= o^{\gamma+1}, && Z^- e^{\gamma} = o^{\gamma} && Z^\pm o^\beta = 0,  \\
        Z^+ e^{[\varepsilon]} &= o^1, && Z^- e^{[\varepsilon]} = 0, && Z^+ e^{[\delta]} = o^1, \;\;Z^- e^{[\delta]} = o^{k+1}. 
    \end{align*}
    %\[
     %   Z^+ e^{\gamma} = o^{\gamma+1}, \;\; Z^- e^{\gamma} = o^{\gamma} \;\; Z^\pm o^\beta = 0,  
    %\]
    %\[
     %   Z^+ e^{[\varepsilon]} = o^1, \;\; Z^- e^{[\varepsilon]} = 0,\;\; Z^+ e^{[\delta]} = o^1, \;\; Z^- e^{[\delta]} = o^{k+1}.  
    %\]
    The representations $\Xi^-_{k, \varepsilon, \delta}$ are similar as $\Xi^+_{k, \varepsilon, \delta}$, but with the action of $L$
    \begin{align*}
         L e^\alpha = -e^\alpha, \;\; \text{for } \alpha = 1, \dots, k, \varepsilon, \delta, &&  L o^\beta = o^\beta, \;\; \text{for } \beta = 1, \dots, k+1.
    \end{align*}
    \item Let $\Eta^+_{\mu, n}$ be the representation of $\operatorname{SF}_{2n}$, where $\mu \in \C^\times$ and $n \in \Z_{> 0}$, and the underlying vector space is $\C^{2n}$. Fix the basis of "even" vectors $e^1, \dots, e^n$ and "odd" vectors $o^1, \dots, o^n$. Then as in Item 1. for $\alpha=1, \dots, n$
    \begin{align*}
        L e^\alpha = e^\alpha,&& L o^\alpha = -o^\alpha
    \end{align*}
    but now $Z^+$ acts by "Jordan block" and $Z^-$ "diagonally"
    \begin{align*}
         Z^+ e^\alpha = \mu o^\alpha + o^{\alpha+1}, && Z^- e^n = \mu o^n, && Z^- e^\alpha = o^\alpha,&& Z^\pm o^\alpha = 0.
    \end{align*}
    %\[
     %   Z^+ e^\alpha = \mu o^\alpha + o^{\alpha+1}, \;\; Z^- e^n = \mu o^n, \;\; Z^- e^\alpha = o^\alpha,\;\; Z^\pm o^\alpha = 0.
    %\]
    \noindent The representations  $\Eta^-_{\mu, n}$ are the same as  $\Eta^+_{\mu, n}$, but with the action of $L$
    \begin{align*}
        L e^\alpha = -e^\alpha,&& L o^\alpha = o^\alpha.
    \end{align*}
     %\[
     %   L e^\alpha = -e^\alpha,\;\; L o^\alpha = o^\alpha.
    %\]
 \end{enumerate}
 The shapes of the two families of modules are given by
% https://q.uiver.app/#q=WzAsMjQsWzAsMSwiZV57W1xcZXBzaWxvbl19Il0sWzAsMiwiZV4xIl0sWzAsNCwiXFxkb3RzIl0sWzAsNSwiZV5rIl0sWzAsNiwiZV57W1xcZGVsdGFdfSJdLFsyLDIsIm9eMSJdLFsyLDUsIm9eayJdLFsyLDYsIm9ee2srMX0iXSxbMiwzLCJvXjIiXSxbNCwyLCJlXjEiXSxbNiwyLCJvXjEiXSxbNCw0LCIuLi4iXSxbNiw0LCIuLi4iXSxbNCw2LCJlXmsiXSxbMCwzLCJlXjIiXSxbMiw0LCIuLi4iXSxbNiw2LCJvXmsiXSxbNiwzLCJvXjIiXSxbNCwzLCJlXjIiXSxbNCw1LCJlXntrLTF9Il0sWzYsNSwib157ay0xfSJdLFsxLDddLFsxLDAsIlxcWGlfe1xcZXBzaWxvbiwgXFxkZWx0YSwga30iXSxbNSwwLCJcXEV0YV97XFxtdSwga30iXSxbMSw1LCJaXi0iLDFdLFsxLDgsIlpeKyIsMV0sWzMsNiwiWl4tIiwxXSxbMyw3LCJaXisiLDFdLFs0LDcsIlpeLSIsMV0sWzAsNSwiWl4rIiwxXSxbMCwwLCJMIiwxXSxbNSw1LCJMIiwxXSxbMSwxLCJMIiwxXSxbMTQsMTUsIlpeKyIsMV0sWzE0LDgsIlpeLSIsMV0sWzIsNiwiWl4rIiwxXSxbOSwxNywiWl4rIiwxXSxbOSwxMCwiWl4tIiwxXSxbOSw5LCJMIiwxXSxbMTAsMTAsIkwiLDFdLFs0LDQsIkwiLDFdLFsxOSwyMCwiWl4tIiwxXSxbMTMsMTYsIlpeLSIsMV0sWzE5LDE2LCJaXisiLDFdLFszNiwxMCwiIiwxLHsic2hvcnRlbiI6eyJzb3VyY2UiOjIwfSwibGV2ZWwiOjF9XSxbNDMsMjAsIiIsMSx7InNob3J0ZW4iOnsic291cmNlIjoyMH0sImxldmVsIjoxfV1d
\[\begin{tikzcd}[ampersand replacement=\&]
	\& {\Xi_{k, \epsilon, \delta}} \&\&\&\& {\Eta_{\mu, k}} \\
	{e^{[\epsilon]}} \\
	{e^1} \&\& {o^1} \&\& {e^1} \&\& {o^1} \\
	{e^2} \&\& {o^2} \&\& {e^2} \&\& {o^2} \\
	\dots \&\& {\dots} \&\& {\dots} \&\& {\dots} \\
	{e^k} \&\& {o^k} \&\& {e^{k-1}} \&\& {o^{k-1}} \\
	{e^{[\delta]}} \&\& {o^{k+1}} \&\& {e^k} \&\& {o^k} \\
	\& {}
	\arrow["L"{description}, from=2-1, to=2-1, loop, in=55, out=125, distance=10mm]
	\arrow["{Z^+}"{description}, from=2-1, to=3-3]
	\arrow["L"{description}, from=3-1, to=3-1, loop, in=55, out=125, distance=10mm]
	\arrow["{Z^-}"{description}, from=3-1, to=3-3]
	\arrow["{Z^+}"{description}, from=3-1, to=4-3]
	\arrow["L"{description}, from=3-3, to=3-3, loop, in=55, out=125, distance=10mm]
	\arrow["L"{description}, from=3-5, to=3-5, loop, in=55, out=125, distance=10mm]
	\arrow["{Z^-}"{description}, from=3-5, to=3-7]
	\arrow[""{name=0, anchor=center, inner sep=0}, "{Z^+}"{description}, from=3-5, to=4-7]
	\arrow["L"{description}, from=3-7, to=3-7, loop, in=55, out=125, distance=10mm]
	\arrow["{Z^-}"{description}, from=4-1, to=4-3]
	\arrow["{Z^+}"{description}, from=4-1, to=5-3]
    \arrow["{Z^-}"{description}, from=4-5, to=4-7]
	\arrow["{Z^+}"{description}, from=5-1, to=6-3]
	\arrow["{Z^-}"{description}, from=6-1, to=6-3]
	\arrow["{Z^+}"{description}, from=6-1, to=7-3]
	\arrow["{Z^-}"{description}, from=6-5, to=6-7]
	\arrow[""{name=1, anchor=center, inner sep=0}, "{Z^+}"{description}, from=6-5, to=7-7]
	\arrow["L"{description}, from=7-1, to=7-1, loop, in=55, out=125, distance=10mm]
	\arrow["{Z^-}"{description}, from=7-1, to=7-3]
	\arrow["{Z^-}"{description}, from=7-5, to=7-7]
	\arrow[shorten <=5pt, from=0, to=3-7]
	\arrow[shorten <=4pt, from=1, to=6-7]
\end{tikzcd}\]
\end{definition}
\begin{remark}
    It was also shown in \cite{zaitsev_nikolenko_1971} that for more than two generators the exterior algebra, and so $\operatorname{SF}_{2n}$ for $n>1$, becomes \textit{wild}, which means the indecomposable modules cannot be classified. With this in mind, we do not attempt this for general $H(\mb{m}, t, \matfont{d}, \matfont{u})$. Note, however, that the modules corresponding $\Xi^\pm_{k, \varepsilon, \delta}$ and $\Eta^\pm_{\mu, n}$ can be easily constructed for many examples with more than two nilpotent generators, simply by setting the action to be trivial for all $X_k$, but a pair $X^\pm_l$ for some $l$, provided such a pair is present.
\end{remark} 
With this remark in mind, we now construct modules analogous to $\Xi^\pm_{k, \varepsilon, \delta}$ for the algebra $\operatorname{N}_2$ of Example \ref{ex2}. Moreover, we will show these examples to be transparent in the category $\operatorname{N}_2\operatorname{-mod}$, thereby showing the M\"uger center $\mathcal{Z}_{(2)}(\operatorname{N}_2\operatorname{-mod})$ of $\operatorname{N}_2\operatorname{-mod}$ contains infinitely many non-isomorphic, non-semisimple, indecomposable modules. 
\begin{example}
\label{XiandEtamodules}
Let $\Xi^\mb{v}_{k, \varepsilon, \delta}$ be representations of $\operatorname{N}_2$ at $t_1 = t_2 =1$ for $\mb{v} \in \Z^3_4$ and $k \in \Z_{\geq 0}$ and $\varepsilon, \delta \in \{0, 1\}$, with the underlying vector space $\C^{2k +\varepsilon +\delta +1}$. Fix a basis a basis of the vectors $e^{\varepsilon}, e^1, \dots, e^k, e^{[\delta]}$ where $e^{\varepsilon}$ appears only if $\varepsilon = 1$ and $e^{\delta}$ appears only if $\delta = 1$. Fix also a basis of vectors $o^1, \dots, o^{k+1}$. Then we have the following module structure
\begin{align*}
     K_a e^\alpha = i^{v_a+2\alpha} e^\alpha, \;\; \text{for } \alpha = 1, \dots, k, \varepsilon, \delta, && K_a o^\beta = i^{v_a+2\beta-1}o^\beta, \;\; \text{for } \beta = 1, \dots, k+1,
\end{align*}
    where $i^\epsilon := i^0 = 1$ and $i^\delta := i^{k+1}$, and for $\gamma=1, \dots, k$
    \begin{align*}
         X^+ e^{\gamma} &= o^{\gamma+1}, && X^- e^{\gamma} = o^{\gamma} && X^\pm o^\beta = 0, \\
          X^+ e^{[\varepsilon]} &= o^1, && X^- e^{[\varepsilon]} = 0,&& X^+ e^{[\delta]} = o^1, \;\;X^- e^{[\delta]} = o^{k+1}.
    \end{align*}
    \noindent and for all $\alpha, \beta$
    \[
        Z^\pm e^\alpha = 0, \;\; Z^\pm o^\beta = 0.  
    \]
% Let $\Eta^\mb{v}_{\mu, n}$, for $\mb{v} \in \Z^3_4$ and $\mu \in \C^\times$ and $n \in \Z_{> 0}$, with the underlying vector space $\C^{2n}$. Fix the basis of vectors $e^1, \dots, e^n$ and vectors $o^1, \dots, o^n$. Then as in item 1) for $\alpha=1, \dots, n$
%     \[
%         K_a e^\alpha = i^{v_a} e^\alpha,\;\; K_a o^\alpha = i^{v_a+1}o^\alpha
%     \]
%     $X^+$ acts by "Jordan block" and $X^-$ "diagonally"
%     \[
%         X^+ e^\alpha = \mu o^\alpha + o^{\alpha+1}, \;\; X^- e^n = \mu o^n, \;\; X^- e^\alpha = o^\alpha,\;\; X^\pm o^\alpha = 0.
%     \]
%     and for all $\alpha$
%     \[
%         Z^\pm e^\alpha = 0, \;\; Z^\pm o^\beta = 0.  
%     \]
It is easy to check $\Xi^\mb{v}_{k, \varepsilon, \delta}$ are $\operatorname{N}_2$-modules. In particular, we have 
\[
    K_a o^\beta = K_a X^- e^\beta = K_a X^+ e^{\beta-1} = i^{v_a+2\beta-1} o^\beta. 
\]
\end{example}
\begin{remark}
    \noindent We note that an analogous family of $\Eta^\mb{v}_{\mu, n}$ does not exist, if we choose the action of $X^+$ to be non-zero and of $Z^-$ to be trivial. The reason is that since we would have $X^+ e^\alpha= \mu o^\alpha +o^{\alpha+1}$ and $X^- e^\alpha = \mu o^\alpha$. If we choose $K_a e^\mb{\alpha} = i^{v_a+2\alpha}$, we would need to have $K_a o^{\alpha} = i^{v_a+2\alpha+1} \mu o^\alpha  = i^{v_a+2\alpha-1} \mu o^\alpha$, which is false. In particular, this would require $K_a X^\pm = - X^\pm K_a$ for any $K_a$.
\end{remark}
\begin{remark}
    Another way to realize Zaitsev-Nikolenko-type modules in $\operatorname{N}_2\operatorname{-mod}$ is to consider the subalgebra $\langle L, Z^\pm_l \rangle$ that is in fact isomorphic to $\operatorname{SF}_{2t_2}$ and thus it is easy to construct representations of the shape both $\Xi^+_{k, \varepsilon, \delta}$ and $\Eta^+_{\mu, n}$.
\end{remark}
\noindent We close this section with another characterisation of $\mathcal{Z}_{(2)}(\operatorname{N}_2\operatorname{-mod})$
\begin{proposition}
\label{NenicuNonSSMuger}
    Let $(H, R_{\matfont{z}}R_{\pmb{\alpha}})$ be a Nenciu type Hopf algebra generated by grouplike $K_a$, and nilpotent $X_k$, $Z^\pm_l$ generators for $a = 1, ..., s$, $k = 1, ...t_1$ and $l = 1, ..t_2$ as in the beginning of Section 5.4. Then the M\"uger centre $\mathcal{Z}_{(2)}(H\operatorname{-mod})$ is a non-semisimple subcategory of $H\operatorname{-mod}$. 
\end{proposition}
\begin{proof}
    Consider the Hopf subalgebra $A = \langle K_a, X_k| a=1, ..., s; k = 1, ..., t_1 \rangle \subset H$, carrying the R-matrix $R_{\matfont{z}}$. It is a non-semismiple triangular Hopf algebra, and $A\operatorname{-mod}$ is a non-semisimple symmetric monoidal category. Since, all relations are diagonal and, by assumption, neither of $X_k$ is present in $R_{\pmb{\alpha}}$, there exists a quotient map of Hopf algebras
    \begin{align*}
        q: H  \twoheadrightarrow A, &&
           K_a \mapsto K_a, &&
           X_k \mapsto X_k, &&
           Z^\pm_l \mapsto 0,
    \end{align*}
    such that $q(R) = R_{\matfont{z}}$. The map $q$ induces a faithfully full braided monoidal functor $I: A\operatorname{-mod} \hookrightarrow H\operatorname{-mod}$. As it is faithfully full braided, we have $I(A\operatorname{-mod}) \subset \mathcal{Z}_{(2)}(H\operatorname{-mod})$, and $I$ is in fact faithfully full symmetric monoidal onto its image. Since $A\operatorname{-mod}$ is a non-semisimple symmetric monoidal categiory and it symmetrically embeds into $\mathcal{Z}_{(2)}(H\operatorname{-mod})$, the later has to be non-semisimple as well.
\end{proof}
\begin{remark}
    The statement remains true if $\pmb{\alpha} = \mb{0}$, as then $\operatorname{N}_2 \operatorname{-mod}$ is a symmetric category.
\end{remark}
\begin{remark}
    Proposition \ref{NenicuNonSSMuger} can be also proved by finding a suitable non-semisimple  indecomposable module. The proof suggests a good example would be any indecomposable projective module over $A$. For the Hopf algebra $\operatorname{N_2}$ this is illustrated in Example \ref{N2ExNSSMod} and we also constructed the modules $\Xi^\mb{v}_{k, \varepsilon, \delta}$ for $\operatorname{N}_2$ that serve the purpose. This in fact shows that $\mathcal{Z}_{(2)}(\operatorname{N_2-mod})$ is non-semisimple, with infinitely many indecomposable objects.
\end{remark}

\section{Representation theory of the non-factorizable biproducts of $u_q \mathfrak{sl}_2$}
In this section we will consider some properties of the representation theory of the second family of Hopf algebras presented in \cite{faes2025non}. We will see many of them will follow from that of their components $U$ and $H$.

\subsection{Projective modules of $u_q \mathfrak{sl}_2 \ltimes H(\mb{m}, t, \mb{d}, \mb{u})$}
We will describe the projective modules of $U \ltimes H$ in terms of those of $u_q \mathfrak{sl}_2$. We recall first the well-known basis of $Z(u_q \mathfrak{sl}_2)$ described in \cite{gainutdinov_tipunin_2009}\footnote{More precisely \cite{gainutdinov_tipunin_2009} describe the basis of $Z(\Bar{U}_q \mathfrak{sl}_2)$, but the one of $Z(u_q \mathfrak{sl}_2)$ is obtained straightforwardly by multiplying each element with the central idempotent $\frac{\1+K^{r'}}{2}$ and taking the quotient by the ideal $\langle K^{r'}-\1 \rangle$.} 

\begin{proposition}
\label{uqsl2center}
The centre $Z(u_q \mathfrak{sl}_2)$ has a basis consisting,  for all $\sigma=1, ..., r'-1$, of the nilpotent elements
\begin{align*}
    \mb{w}^+_\sigma &= \zeta_\sigma \sum^{\sigma-1}_{n=1}\sum^n_{i=0}\sum^{r'-1}_{j=0}([i]!)^2 q^{j(\sigma-1-2n)} \begin{bmatrix} \sigma-n+i-1\\ i \end{bmatrix} \begin{bmatrix} n \\ i \end{bmatrix} F^{r'-1-i}E^{r'-1-i}K^j\\
    \mb{w}^-_\sigma &= \zeta_\sigma \sum^{p-\sigma-1}_{n=1}\sum^n_{i=0}\sum^{r'-1}_{j=0}(-1)^{i+j}([i]!)^2 q^{j(r'-s-1-2n)} \begin{bmatrix} r'-\sigma-n+i-1\\ i \end{bmatrix} \begin{bmatrix} n \\ i \end{bmatrix} F^{r'-1-i}E^{r'-1-i}K^j
\end{align*}
and the central idempotents
\begin{align*}
    \mb{e}_0 &= \zeta_0 \sum^{r'-1}_{n=1}\sum^n_{i=0}\sum^{r'-1}_{j=0}(-1)^{i+j}([i]!)^2 q^{j(r'-1-2n)} \begin{bmatrix} r'-\sigma-n+i-1\\ i \end{bmatrix} \begin{bmatrix} n \\ i \end{bmatrix} F^{r'-1-i}E^{r'-1-i}K^j\\
    \mb{e}_\sigma &= \frac{q^\sigma + q^{-\sigma}}{[\sigma]^2}(\mb{w}^+_\sigma + \mb{w}^-_\sigma)\\
    &+\zeta_\sigma\sum^{r'-2}_{m=0}\sum^{r'-1}_{j=0}\left(\sum^{\sigma-1}_{n=0} q^{j(\sigma-1-2n)} \text{B}^+_{n, r'-1-m}(\sigma) + \sum^{r'-\sigma-1}_{n=0} q^{j(r'-\sigma-1-2k)} \text{B}^-_{k, r'-1-m}(r'-\sigma) \right) F^m E^m K^j,\\
    \mb{e}_{r'} &= \zeta_{r'} \sum^{r'-1}_{n=0}\sum^n_{i=0} \sum^{r'-1}_{j=0}([i]!)^2 q^{j(r'-1-2n)} \begin{bmatrix} r'--n+i-1\\ i \end{bmatrix} \begin{bmatrix} n \\ i \end{bmatrix} F^{r'-1-i}E^{r'-1-i}K^j
\end{align*}
where $\text{B}^\pm_{n, m}$ are non-zero numbers and 
\begin{align*}
    \zeta_0 &= \frac{(-1)^{r'-1}}{r}\frac{1}{([r'-1]!)^2}, \\
    \zeta_\sigma &= \frac{(-1)^{r'-\sigma-1}}{r}\frac{[\sigma]^2}{([r'-1]!)^2}, \;\; \; 1\leq \sigma\leq r'-1,\\
    \zeta_{r'} &= \frac{1}{r'} \frac{1}{([r'-1]!)^2}.
\end{align*}
\end{proposition}

Recall from \cite[Eq. (8.6)]{beliakova_blanchet_gainutdinov_2021}, that the full system of idempotents for $u_q \mathfrak{sl}_2$ has the form
\begin{equation}
    \label{sl2fullsystem}
    I_{p, \sigma} = \varphi^K_p \mb{e}_\sigma, \;\; p \in \Z_{r'}, \;\; 1 \leq \sigma \leq r'-1, \;\; p-\sigma \equiv 1\mod 2,
\end{equation}
where $\mb{e}_s$ are central idempotents introduced in Proposition \ref{uqsl2center}, and 
$
    \varphi^K_p = \frac{1}{r'}\sum^{r'-1}_{j=0} q^{-2pj} K^j,
$
is the projector to the $p$th  eigenspace of $K$, analogous to Definition \ref{subIdems}. It was shown in \cite{gainutdinov_semikhatov_tipunin_feigin_2006} that the corresponding indecomposable projective modules, $\mathcal{P}^+_\sigma := u_q \mathfrak{sl}_2 I_{\sigma-1, \sigma}$ and $\mathcal{P}^-_{r'-\sigma} := u_q \mathfrak{sl}_2 I_{-\sigma-1, \sigma}$, have the following shape
\[
    % https://q.uiver.app/#q=WzAsNCxbMSwwLCJcXG1hdGhjYWx7U31eXFxwbV9zIl0sWzAsMSwiXFxtYXRoY2Fse1N9XlxcbXBfe3InLXN9Il0sWzIsMSwiXFxtYXRoY2Fse1N9XlxcbXBfe3InLXN9Il0sWzEsMiwiXFxtYXRoY2Fse1N9XlxccG1fcyJdLFswLDEsIkYiLDJdLFswLDIsIkUiXSxbMSwzLCJFIiwyXSxbMiwzLCJGIl1d
\begin{tikzcd}
	& {\mathcal{S}^\pm_\sigma} \\
	{\mathcal{S}^\mp_{r'-\sigma}} && {\mathcal{S}^\mp_{r'-\sigma}}, \\
	& {\mathcal{S}^\pm_\sigma}
	\arrow["F"', from=1-2, to=2-1]
	\arrow["E", from=1-2, to=2-3]
	\arrow["E"', from=2-1, to=3-2]
	\arrow["F", from=2-3, to=3-2]
\end{tikzcd}
\] 
where ${\mathcal{S}^\pm_\sigma}$ are the simple modules corresponding to $\mathcal{P}^\pm_{\sigma}$. 
We now consider the case $U \ltimes H$. Let $H$ have grouplike generators $K_a$ for $a=1, ..., s$. The corresponding full system of idempotents was found in Proposition \ref{propIdems}. 
\begin{proposition}
\label{productFullSystem}
    Let $U \ltimes H$ with $U =u_q \mathfrak{sl}_2$ and let $H = H(\mb{m}, t, \mb{d}, \mb{u})$ have the full system of idempotents $\{\omega_\mb{p}\}_{\mb{p}\in \Z_\mb{m}}$. Then the full system of idempotents of $U \ltimes H$ is 
    \[
       \Omega_{p, \sigma, \mb{p}}:= I_{p, \sigma} \omega_\mb{p}, \;\;\;\text{for} \;\; p=0, ..., r'-1, \;\; \mb{p} \in \Z_\mb{m}, \;\; p-s \equiv 1 \mod 2.
    \]
\end{proposition}
\begin{proof}
    We first note that $I_{p, \sigma}$ and $\omega_\mb{p}$ commute for all $p, \sigma$ and $\mb{p}$. Since both systems are separately orthogonal, we find
    \begin{align*}
        \Omega_{p_1, \sigma_1, \mb{p}_1}\Omega_{p_2, \sigma_2, \mb{p}_2} &= I_{p_1, \sigma_1}\omega_{\mb{p}_1} I_{p_2, \sigma_2}\omega_{\mb{p}_2} =  
         I_{p_1,\sigma_1} I_{p_2, \sigma_2}\omega_{\mb{p}_2} \omega_{\mb{p}_1} \\
         &= \delta_{p_1, p_2} \delta_{\sigma_1, \sigma_2} \prod^s_{a=1} \delta_{(\mb{p}_{1})_a, (\mb{p}_{2})_a}.
    \end{align*}
    Similarly, since both systems decompose the unit we have
    \[
        \sum_{p\in \Z_{r'}, \; \sigma=1, ..., r'-1, \; \mb{p} \in \Z_\mb{m}} \Omega_{p, \sigma, \mb{p}} = \sum_{p\in \Z_{r'}, \; \sigma=1, ..., r'-1} I_{p, \sigma} (\sum_{\mb{p} \in \Z_\mb{m}} \omega_\mb{p}) = \1_U (\1_H) = \1. 
    \]
    Finally, for primitivity, we first note that $\mb{e}_\sigma$ are primitive idempotentents. Now, similarly to the reasoning in \cite{beliakova_blanchet_gainutdinov_2021} the products $\varphi^k_p \omega_\mb{p}$ are projectors to the eigenspaces of the abelian group algebra $\C[K, \mb{K}]$ that we showed in Proposition \ref{propIdems} to be primitive. Since the two pieces commute and are primitive, the result will be primitive as well. 
\end{proof}
\noindent We can now find the indecomposable projective modules for $U \ltimes H$. 
\begin{definition}
\label{productProjective}
    For the setup above, define 
    $
        \mb{P}^+_{\sigma, \mb{p}}:= (U\ltimes H) \Omega_{\sigma-1, \sigma, \mb{p}}
    $
    and
    $
        \mb{P}^-_{r'-\sigma, \mb{p}}:= (U\ltimes H) \Omega_{-\sigma-1, \sigma, \mb{p}}.
    $
    with the basis 
$$
\left\{{}^{\pm}\mathrm{a}^\mb{r}_n, {}^{\pm}\mathrm{~b}^\mb{r}_n\right\}_{0 \leq n \leq \sigma-1} \cup\left\{{}^{\pm}\mathrm{x}^\mb{r}_m, {}^{\pm}\mathrm{y}^\mb{r}_m\right\}_{0 \leq m \leq r'-\sigma-1},
$$
where $\left\{{}^{\pm}\mathrm{b}^\mb{r}_n\right\}_{0 \leq n \leq \sigma-1}$ is the basis corresponding to the top module ${}^{\pm}\mathcal{S}_\sigma^{a, \mb{r}}$ in (1.2), $\left\{{}^{\pm}\mathrm{a}^\mb{r}_n\right\}_{0 \leq n \leq \sigma-1}$ to the bottom ${}^{\pm}\mathcal{S}_\sigma^{a, \mb{r}},\left\{{}^{\pm}\mathrm{x}^\mb{r}_m\right\}_{0 \leq k \leq r'-\sigma-1}$ to the right ${}^{\pm}\mathcal{S}_{r'-\sigma}^{-a, \mb{r}}$, and $\left\{{}^{\pm}\mathrm{y}^\mb{r}_m\right\}_{0 \leq k \leq r'-\sigma-1}$ to the left module ${}^{\pm}\mathcal{S}_{r'-\sigma}^{-a, \mb{r}}$. The superscript $\mb{r}\in \Z^t_2$ keeps track of the $X_k$ that acted on the generator.\\
The top-left superscript $\pm$ distinguishes between $\mb{P}^+_{\sigma, \mb{p}}$ and $\mb{P}^-_{r'-\sigma, \mb{p}}$, and for notational brevity we will concentrate on the positive variant, as they are analogous.\\

\noindent The $U \ltimes H$-action on $\mb{P}^+_{\sigma, \mb{p}}$ is given by
\begin{align*}
&K \mathrm{x}^\mb{r}_m=-(-1)^\mb{|r|}a q^{r'-\sigma-1-2 k} \mathrm{x}^\mb{r}_m,& K \mathrm{y}^\mb{r}_m =-(-1)^\mb{|r|}a q^{r'-\sigma-1-2 k} \mathrm{y}^\mb{r}_m, &\;\; 0 \leq m \leq r'-\sigma-1,\\
&K \mathrm{a}^\mb{r}_n=(-1)^\mb{|r|}a q^{s-1-2 n} \mathrm{a}^\mb{r}_n,& K \mathrm{b}^\mb{r}_n =(-1)^\mb{|r|}a q^{\sigma-1-2 n} \mathrm{~b}^\mb{r}_n, &\;\; 0 \leq n \leq \sigma-1,
\end{align*}
\begin{align*}
    E\mathrm{x}^\mb{r}_m &= -a[k][r'-\sigma-k] \mathrm{x}^\mb{r}_{k-1}, \;\; 0 \leq m \leq r'-\sigma-1, \;\; \left(\text{with } \mathrm{x}^\mb{r}_{-1} \equiv 0\right)\\
    E \mathrm{y}^\mb{r}_m &= \begin{cases}-a[k][r'-\sigma-k] \mathrm{y}^\mb{r}_{k-1} +  \mathrm{x}^\mb{r}_{m-1}, & 1 \leq m \leq r'-\sigma-1 \\
\mathrm{a}^\mb{r}_{\sigma-1}, & m=0\end{cases} \\
    E \mathrm{a}^\mb{r}_n &=a[n][\sigma-n] \mathrm{a}^\mb{r}_{n-1}, \;\; 0 \leq m \leq \sigma-1, \;\; \left(\text {with } \mathrm{x}^\mb{r}_{-1} \equiv 0\right)\\
    E \mathrm{b}^\mb{r}_n &= \begin{cases}a[n][\sigma-n] \mathrm{b}^\mb{r}_{n-1}+\mathrm{a}_{n-1}, & 1 \leq n \leq \sigma-1, \\
\mathrm{x}^\mb{r}_{r'-\sigma-1}, & n=0,\end{cases}\\
&\\
    F \mathrm{x}^\mb{r}_m &= \begin{cases}\mathrm{x}^\mb{r}_{k+1}, & 0 \leq k \leq r'-\sigma-2, \\
\mathrm{a}^\mb{r}_0, &k=r'-\sigma-1,\end{cases} \\
    F \mathrm{y}^\mb{r}_m &=\mathrm{y}^\mb{r}_{k+1},\;\; 0 \leq m \leq r'-\sigma-1 \quad\left(\text {with } \quad \mathrm{x}^\mb{r}_{r'-\sigma} \equiv 0\right),\\
    F \mathrm{a}^\mb{r}_n &=\mathrm{a}^\mb{r}_{n+1}, \;\; 0 \leq n \leq \sigma-1 \quad\left(\text {with } \mathrm{a}^\mb{r}_\sigma \equiv 0\right)\\
    F \mathrm{b}^\mb{r}_n &= \begin{cases}  \mathrm{b}^\mb{r}_{n+1}, & 1 \leq n \leq \sigma-2, \\ \mathrm{y}_0, & n = \sigma-1\end{cases}
\end{align*}
Moreover, for any $\mathrm{z}^\mb{r}_l\in \{\mathrm{a}^\mb{r}_n, \mathrm{b}^\mb{r}_n, \mathrm{x}^\mb{r}_m, \mathrm{y}^\mb{r}_m\}$, we have
\begin{align*}
    X_k \mathrm{x}^\mb{r}_m &= \begin{cases}
        (-1)^{m-r'+s}\pmb{\xi}^{\mb{u}_k\cdot \sum_{1\leq j<k}r_j \mb{d}_j}\mathrm{x}^{\Bar{\mb{r}}}_m, &\text{ if } r_k = 0 \\
        0, \text{ if } r_k = 1,\\
    \end{cases}\\
    X_k \mathrm{y}^\mb{r}_m &= \begin{cases}
        (-1)^{m+s-1}\pmb{\xi}^{\mb{u}_k\cdot \sum_{1\leq j<k}r_j \mb{d}_j}\mathrm{y}^{\Bar{\mb{r}}}_m, &\text{ if } r_k = 0 \\
        0, \text{ if } r_k = 1,\\
    \end{cases}\\
    X_k \mathrm{a}^\mb{r}_n &= \begin{cases}
        (-1)^{n+1}\pmb{\xi}^{\mb{u}_k\cdot \sum_{1\leq j<k}r_j \mb{d}_j}\mathrm{a}^{\Bar{\mb{r}}}_n, &\text{ if } r_k = 0 \\
        0, \text{ if } r_k = 1,\\
    \end{cases}\\
         X_k \mathrm{b}^\mb{r}_n &= \begin{cases}
        (-1)^n \pmb{\xi}^{\mb{u}_k\cdot \sum_{1\leq j<k}r_j \mb{d}_j} \mathrm{b}^{\Bar{\mb{r}}}_n, &\text{ if } r_k = 0 \\
        0, \text{ if } r_k = 1,\\
        \end{cases}
\end{align*}
\noindent where $\Bar{\mb{r}} \in \Z^t_2$ is such that $\Bar{r}_i = r_i$ for $i\neq k$ and $\Bar{r}_k=1$. Finally,
\[
    K_a \mathrm{z}^\mb{r}_l = \xi^{\mb{r}\cdot \mb{d}_a}\xi^{-p_a} \mathrm{z}^\mb{r}_l, 
\]
where $\mathrm{z}^\mb{r}_l \in \{\mathrm{b}^\mb{r}_n, \mathrm{x}^\mb{r}_m, \mathrm{y}^\mb{r}_m, \mathrm{a}^\mb{r}_n\}$.
The case $\mb{P}^-_{r'-\sigma, \mb{p}}$ has analogous $H$-component action and the $U$ action as in \cite[Appendix C.2.2]{gainutdinov_semikhatov_tipunin_feigin_2006} (for all occurrences of $q$ replace, for instance, $q^{r'-\sigma-1-2 k}$ with $-q^{r'-\sigma-1-2 k}$ etc.).
\end{definition}
\noindent For the case $t = 2$, and $X_1 = X^+$, $X_2 = X^-$, such a $\mb{P}^\pm_{\sigma, \mb{p}}$ has the following shape.

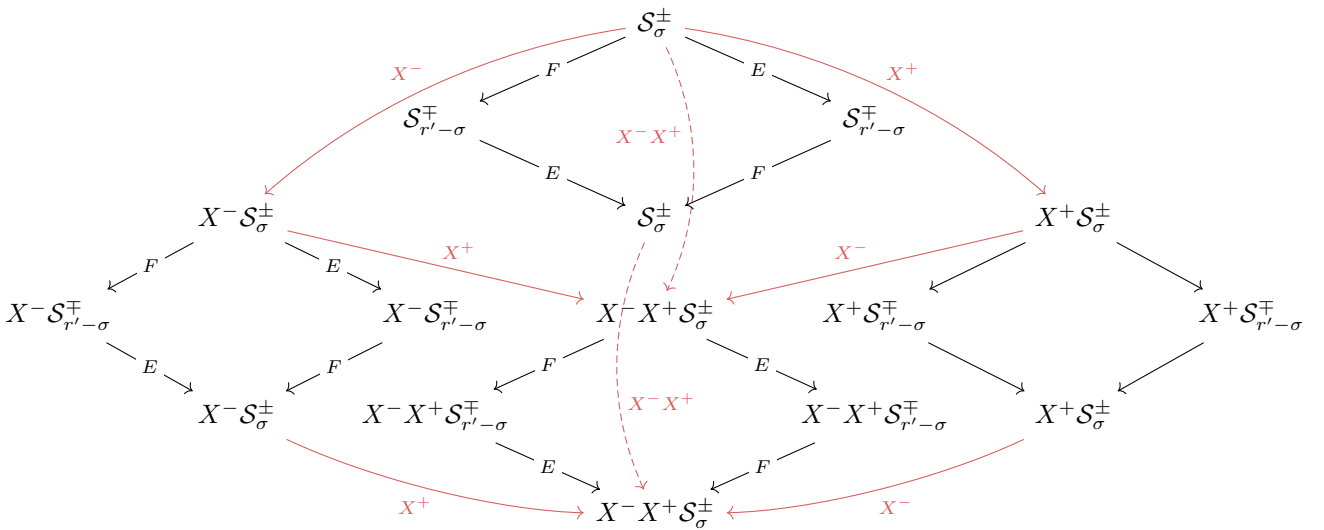
\[\noindent\makebox[\textwidth]{\begin{tikzcd}[ampersand replacement=\&]
	\&\&\& {\mathcal{S}^\pm_\sigma} \\
	\&\& {\mathcal{S}^\mp_{r'-\sigma}} \&\& {\mathcal{S}^\mp_{r'-\sigma}} \\
	\& {X^-\mathcal{S}^\pm_\sigma} \&\& {\mathcal{S}^\pm_\sigma} \&\& {X^+\mathcal{S}^\pm_\sigma} \\
	{X^-\mathcal{S}^\mp_{r'-\sigma}} \& {} \& {X^-\mathcal{S}^\mp_{r'-\sigma}} \& {X^-X^+\mathcal{S}^\pm_\sigma} \& {X^+\mathcal{S}^\mp_{r'-\sigma}} \& {} \& {X^+\mathcal{S}^\mp_{r'-\sigma}} \\
	\& {X^-\mathcal{S}^\pm_\sigma} \& {X^-X^+\mathcal{S}^\mp_{r'-\sigma}} \&\& {X^-X^+\mathcal{S}^\mp_{r'-\sigma}} \& {X^+\mathcal{S}^\pm_\sigma} \\
	\&\&\& {X^-X^+\mathcal{S}^\pm_\sigma}
	\arrow["F"{description}, from=1-4, to=2-3]
	\arrow["E"{description}, from=1-4, to=2-5]
	\arrow["{X^-}"', color={rgb,255:red,214;green,92;blue,92}, curve={height=18pt}, from=1-4, to=3-2]
	\arrow["{X^+}", color={rgb,255:red,214;green,92;blue,92}, curve={height=-18pt}, from=1-4, to=3-6]
	\arrow["{X^-X^+}"'{pos=0.3}, color={rgb,255:red,214;green,92;blue,92}, curve={height=-18pt}, dashed, from=1-4, to=4-4]
	\arrow["E"{description}, from=2-3, to=3-4]
	\arrow["F"{description}, from=2-5, to=3-4]
	\arrow["F"{description}, from=3-2, to=4-1]
	\arrow["E"{description}, from=3-2, to=4-3]
	\arrow["{X^+}", color={rgb,255:red,214;green,92;blue,92}, from=3-2, to=4-4]
	\arrow["{X^-X^+}"{pos=0.7}, color={rgb,255:red,214;green,92;blue,92}, curve={height=18pt}, dashed, from=3-4, to=6-4]
	\arrow["{X^-}"', color={rgb,255:red,214;green,92;blue,92}, from=3-6, to=4-4]
	\arrow[from=3-6, to=4-5]
	\arrow[from=3-6, to=4-7]
	\arrow["E"{description}, from=4-1, to=5-2]
	\arrow["F"{description}, from=4-3, to=5-2]
	\arrow["F"{description}, from=4-4, to=5-3]
	\arrow["E"{description}, from=4-4, to=5-5]
	\arrow[from=4-5, to=5-6]
	\arrow[from=4-7, to=5-6]
	\arrow["{X^+}"', color={rgb,255:red,214;green,92;blue,92}, curve={height=12pt}, from=5-2, to=6-4]
	\arrow["E"{description}, from=5-3, to=6-4]
	\arrow["F"{description}, from=5-5, to=6-4]
	\arrow["{X^-}", color={rgb,255:red,214;green,92;blue,92}, curve={height=-12pt}, from=5-6, to=6-4]
\end{tikzcd}}\]
\noindent There are analogous red arrows linking left and right simples of the copies of $\mathcal{P}^\pm_\sigma$, as well as loops for the action of $K$ and $K_a$ for each weight space, suppressed for legibility. Analogously to the explanation in Remark \ref{reorderingRem}, we do not include the reordering constants from subsequent actions of $X_k$ in the diagram. 
\begin{remark}
    It is important to note that the projective modules of the usual trivial product Hopf algebra $U \otimes H$ have similar shape and basis. The key difference is the $(-1)^\mb{r}$ factor involved in the action of $K$, resulting from the fact $K$ and $X_k$ anticommute, as well as reordering signs between the actions of $E$ and $X_k$, where we count the monomial $v \otimes \1_H$ for the starting point, where $v$ is the highest weight vector of the underlying $U$-representation.
\end{remark} 

Similar strategy can be used to construct an indecomposable module $V\ltimes W \in (U\ltimes H)\operatorname{-mod}$ from any indecomposable modules $V\in U\operatorname{-mod}$ and $W \in U\operatorname{-mod}$, by taking the tensor product of the underlying spaces and adjusting the signs of $K$ weight spaces for coherence with the action of the $X_k$. We consider an example.
\begin{example}
\label{prodExample1}
    Consider $U\ltimes H = u_q \mathfrak{sl}_2 \ltimes \operatorname{SF}_2$. Let $\mathcal{M}^+_1(2)$ be an M-shaped module of $u_q \mathfrak{sl}_2$, see \cite[Appendix A.1.2]{gainutdinov_semikhatov_tipunin_feigin_2006}, with the basis $\mathrm{a}, \mathrm{x}_m, \mathrm{y}_m$, for $m=1, ..., r'-2$ with the action
    \begin{align*}
         K\mathrm{a} &= q^{0} \mathrm{a}, \\ 
         K\mathrm{x}_m &= -q^{r'-2-2m} \mathrm{x}, && K\mathrm{y}_m = -q^{r'-2-2m} \mathrm{y}\\
         E\mathrm{a} &= \mathrm{x}_{r'-2}, && \hspace{0.3cm} F\mathrm{a} = \mathrm{y}_{0}\\
         E \mathrm{x}_m &= -[m][r'-1-m] \mathrm{x}_{m-1},  && F\mathrm{x}_m = \mathrm{x}_{m+1},\\
          E\mathrm{y}_m &= -[k][r'-1-m] \mathrm{x}_{m-1},  && F \mathrm{y}_{m} = \mathrm{y}_{m-1}   
    \end{align*}
     where $\mathrm{x}_{-1} = \mathrm{x}_{r'-1} = \mathrm{x}_{-1} = \mathrm{x}_{r'-1} =0$. \\
     \noindent We take the product with $P_0$ of $SF_2$. The resulting basis is $\mathrm{a}^\mb{r}, \mathrm{x}^\mb{r}_m, \mathrm{y}^\mb{r}_m$ for $\mb{r}\in \Z^2_2$, with the action
     \begin{align*}
          K\mathrm{a}^\mb{r} &= \mathrm{a}^\mb{r}, \\
          K\mathrm{x}^\mb{r}_m &= -q^{r'-2-2m} (-1)^{|\mb{r}|} \mathrm{x}^\mb{r}_m, && K\mathrm{y}_m = -q^{r'-2-2m} (-1)^{|\mb{r}|} \mathrm{y}^\mb{r}_m \\
          E\mathrm{a}^\mb{r} &= \mathrm{x}^\mb{r}_{r'-2}, && F\mathrm{a}^\mb{r} \hspace{0.3cm}= \mathrm{y}^\mb{r}_{0}\\
          E \mathrm{x}^\mb{r}_m &= -[m][r'-1-m] \mathrm{x}^\mb{r}_{m-1}, && F\mathrm{x}^\mb{r}_m = \mathrm{x}^\mb{r}_{m+1},\\
          E\mathrm{y}^\mb{r}_m &= -[k][r'-1-m] \mathrm{y}^\mb{r}_{m-1}, && F \mathrm{y}^\mb{r}_{m} = \mathrm{y}^\mb{r}_{m-1} \\
     \end{align*}
Moreover, for any $\mathrm{z}^\mb{r}_l\in \{\mathrm{a}^\mb{r}, \mathrm{x}^\mb{r}_m, \mathrm{y}^\mb{r}_m\}$, and $Z_k = Z^\pm$ we have
\begin{align*}
    Z_k \mathrm{x}^\mb{r}_m &= \begin{cases}
        (-1)^{m-r'+1}(-1)^{\sum_{1\leq j<k}r_j}\mathrm{x}^{\Bar{\mb{r}}}_m, \text{ if } r_k = 0 \\
        0, \text{ if } r_k = 1,
    \end{cases}\\
    Z_k \mathrm{y}^\mb{r}_m &= \begin{cases}
        (-1)^{m}(-1)^{\sum_{1\leq j<k}r_j}\mathrm{y}^{\Bar{\mb{r}}}_m, \text{ if } r_k = 0 \\
        0, \text{ if } r_k = 1,
    \end{cases}\\
    Z_k \mathrm{a}^\mb{r} &= \begin{cases}
        (-1)^{\sum_{1\leq j<k}r_j}\mathrm{a}^{\Bar{\mb{r}}}, \text{ if } r_k = 0 \\
        0, \text{ if } r_k = 1, \\
    \end{cases}
\end{align*}
\noindent where $\Bar{\mb{r}} \in \Z^t_2$ is such that $\Bar{r}_i = r_i$ for $i\neq k$ and $\Bar{r}_k=1$. Finally, 
\begin{align*}
    L \mathrm{a}^\mb{r} &= (-1)^{|\mb{r}|} \mathrm{a}^{\mb{r}}\\
    L \mathrm{z}^\mb{r}_m &= (-1)^{|\mb{r}|} \mathrm{z}^\mb{r}_m,
\end{align*}
where $\mathrm{z}^\mb{r}_m \in \{\mathrm{x}^\mb{r}_m, \mathrm{y}^\mb{r}_m\}$. 
\end{example}

\subsection{Some braiding properties}
As in the previous cases, $U \ltimes H$ can inherit from $H = H(\mb{m}, t, \mb{d}, \mb{u})$ the infinitely generated, non-semisimple M\"uger center. It is easy to see the following.
\begin{proposition}
\label{extInfiniteMuger}
    Let $U \ltimes H$ be as before. If $\mathcal{Z}_{(2)}(H\operatorname{-mod})$ is non-semisimple then so is $\mathcal{Z}_{(2)}(U\ltimes H\operatorname{-mod})$.
\end{proposition}
\begin{proof}
    Let $\{V_i\}_{i\in I}\subset H\operatorname{-mod}$ be a family of transparent, non-semisimple modules over $H$, for a (possibly infinite) indexing set $I$, which exists by hypothesis. Assume also that all of them are finite dimensional with $\dim V_i= n_i$. For any $V_i$, let $\{v^\alpha_i\}_{1\leq \alpha \leq n_i}$ be its generating set. Since all grouplike generators $K_a$ belong to an abelian group algebra $\C[\mb{K}]$, we can assume that $K_a v^\alpha_i = \xi_a^{\gamma^\alpha_{i, a}}$ for some $\gamma^\alpha_{i, a} \in \Z_{m_a}$. We now construct a family $\{\Bar{V}_i\}_{i\in I} \subset (U\ltimes H)\operatorname{-mod}$, by assigning first the action $K \Bar{v}^\alpha_i = \Bar{v}_i$, for $\{\Bar{v}^\alpha_i\}_{1\leq \alpha \leq n_i}$, the generating set of $\Bar{V}_i$ induced from $\{v^\alpha_i\}_{1\leq \alpha \leq n_i}$ of $V_i$, as well as $Ev^\alpha_i = F v^\alpha_i = 0$. Now, clearly $\Bar{V}_i$ is generated by the action of all $X_k$. What remains to be defined is the action of $K$ on the basis generated from $\{\Bar{v}^\alpha_i\}_{1\leq \alpha \leq n_i}$. We do it inductively, and the base case has already been defined. Now, let $\Bar{w}\in \Bar{V}$, be an element with a $K$ eigenvalue $(-1)^{s_0}$ for $0\leq s_0 \leq s$, of the form
    \[
        \Bar{w} = \sum_{i, \;\mb{r} \in \Z^t_2, \; |\mb{r}|\leq s_0} \kappa_{i, \mb{r}} \mb{X}^\mb{r} \Bar{v}_i, 
    \]
    where $\kappa_{i, \mb{r}} \in \C$ are constants. This means that at the inductive step $s_0$ we have basis elements obtained from action of $s_0$-many $X_k$ generators on the generating set. This will indeed have the $K$-eigenvalue $(-1)^{s_0}$. Now, there is a $k_0\in \{1, ..., t\}$ such that $X_{k_0}\Bar{w}\neq 0$, otherwise we are done. Then $K(X_{k_0}\Bar{w}) = (-1)^{s_0} X_{k_0}\Bar{w}$, so we assign this $K$ eigenvalue to all basis elements $X_{k_0}\mb{X}^r\Bar{v}_i$. This way we build the action of $K$ coherent with the multiplication rules of $U \ltimes H$.\\
    Now, we already know the $R_\mb{z} R_\mb{a}$ part of the R-matrix acts symmetrically on any pair $V_i \otimes W \in H\operatorname{-mod}$. We chose the trivial action of $E, F$ so the action of $D R_\mb{z} \Theta R_\mb{a}$ on any product $\Bar{V}_i \otimes \Bar{W} \in (U\ltimes H)$ will also be symmetric. Since $V_i$ were all indecomposable, $\Bar{V}_i$ are as well. 
\end{proof}
\noindent One non-semisimple module would suffice for the proof, but we kept it more general to immediately retrieve the following.
\begin{corollary}
    If $\mathcal{Z}_{(2)}(H\operatorname{-mod})$ has infinitely many finite dimensional non-semisimple indecomposable objects then so does $\mathcal{Z}_{(2)}(U\ltimes H\operatorname{-mod})$.
\end{corollary}

\noindent We now give an example of this fact.
\begin{example}
\label{prodExample2}
    Consider $U\ltimes \operatorname{N}_2$ of Example \ref{extexample2}. In the category $H\operatorname{-mod}$ there is a family of modules $\Xi^\mb{v}_{\epsilon, \delta, k}$, with $\mb{v} \in \Z^3_4$, constructed in Example \ref{XiandEtamodules}. Since the generators $Z^\pm$ do not act, all $\Xi^\mb{v}_{\epsilon, \delta, k}$ are transparent in $H\operatorname{-mod}$. We can apply the strategy of Proposition \ref{extInfiniteMuger} in two ways - by assigning $K$ eigenvalue $1$ to the $e^\alpha$ and $(-1)$ to $o^\beta$ or vice-versa. We denote these modules ${}^+\Xi^\mb{v}_{\epsilon, \delta, k}$ and ${}^-\Xi^\mb{v}_{\epsilon, \delta, k}$ respectively.
\end{example}

\begin{remark}
\label{catEmbeddingPropProd}
    There is a faithfully full braided monoidal functor
    $
        I: \operatorname{SF}_{2t_1} \operatorname{-mod} \rightarrow \mathcal{Z}_{(2)} (U \ltimes \operatorname{N}_{2} \operatorname{-mod}),  
    $
    where $\operatorname{SF}_{2 t_1}$ carries the braiding induced by the R-matrix $R_{\matfont{z}}$ of Proposition \ref{SFprop}, that is in the case $\pmb{\alpha} = \mb{0}$. 
\end{remark}

\section{Examples in the context of 4d TQFTs}

In this final section we consider our examples in the context of 4d TQFTs introduced in CGHP and extend our biproduct construction to Hopf algebras over fields of finite characteristic. We, however, also prove that none of this allows to produce a powerful invariant of 4-manifolds.  

\subsection{Chromatic non-degenerate Hopf algebras}
Let $\mathcal{C} = H\operatorname{-mod}$ be the representation category of a unimodular non-factorizable ribbon Hopf algebra $H$ over a field $\K$ of characteristic $p>2$. The category $\mathcal{C}$ is, by definition, a unimodular non-factorizable ribbon category, as well as $\K$-linear and finite. We recall the following data from CGHP.
\begin{definition}
    Let $P_{\mathbbm{1}}$ be the projective cover of the unit $\mathbbm{1} \in \mathcal{C}$ (the trivial representation of $H$). Let also $\lambda\in H^*$ be the left integral,  $\Lambda \in H$ be the cointegral and $M\in H\otimes H$ the monodromy matrix. The we define the following maps $\Lambda_{P_{\mathbbm{1}}}, \Delta^{P_{\mathbbm{1}}}_0: P_{\mathbbm{1}} \rightarrow P_{\mathbbm{1}}$
    \begin{align*}
        \Lambda_{P_{\mathbbm{1}}}: P_{\mathbbm{1}} \rightarrow \Lambda P_{\mathbbm{1}}, &&
        \Delta^{P_{\mathbbm{1}}}_0: P_{\mathbbm{1}} \rightarrow \lambda(S(M'))M'' P_{\mathbbm{1}}.
    \end{align*}
    That is, the maps induced by the (regular) action of $\Lambda$ and $\lambda(S(M'))M''$, respectively.
\end{definition}
\begin{remark}
    We note that CGHP, the definition of $\Delta^{P_\mathbbm{1}}_0$ stands
    $
         \Delta^{P_{\mathbbm{1}}}_0: P_{\mathbbm{1}} \rightarrow \lambda_R(M')M'' P_{\mathbbm{1}},
    $
    where $\lambda_R \in H^*$ is the right integral, while we use the left integral. But, as noted in Definition \ref{defunimod}, we have that $\lambda_R = \lambda_L \circ S$, so the two definitions are equivalent.
\end{remark}
Since by unimodularity $\Lambda_{P_{\mathbbm{1}}} \neq 0$ (we have $P_{\mathbbm{1}} = P^*_{\mathbbm{1}}$, and $\Lambda_{P_{\mathbbm{1}}}$ is the non-zero endomorphism of $P_{\mathbbm{1}}$ sending the head to the socle), it is easy to see that chromatic compactness implies chromatic non-degeneracy.
With this notation, we introduce the following notions
\begin{definition}[CGHP, Definition 1.7 (2), (3)]
\label{ChromNonDegenCompact}
    Let $\mathcal{C}$ be the category. Then $\mathcal{C}$ is 
    \begin{enumerate}
        \item \textit{chromatic non-degenerate} if $\Delta^{P_{\mathbbm{1}}}_0 \neq 0$,
        \item \textit{chromatic compact} if there exists a scalar $\zeta \in \K^\times$ such that $\Delta^{P_{\mathbbm{1}}}_0 = \zeta \Lambda_{P_{\mathbbm{1}}}$.
    \end{enumerate}
\end{definition}
\noindent By the usual convention we will call a Hopf algebra chromatic non-degenerate (resp. compact) if its representation category is chromatic non-degenerate (resp. compact). We finally point out that if $H$ is \textit{anomalous} in the sense of \cite[Definition 2.20]{faes2025non}, then $\mathcal{C}$ is \textit{twist-degenerate} in the sense of \cite[Definition 1.7 (1)]{costantino2023skein}.\\

\subsection{Examples}
In this subsection we relate the results from the previous sections to the construction in CGHP. We start with the following observation.
\begin{proposition}
\label{snfChromNonCompact}
    Let $H$ be a unimodular, strongly non-factorizable ribbon Hopf algebra. Then it is neither chromatic non-degenerate, nor chromatic compact.
\end{proposition}
\begin{proof}
    By definition of non-factorizability, we have $\lambda(S(M'))M'' = 0$, which implies $\Delta^{P_{\mathbbm{1}}}_0 = 0$, and so $H$ is not chromatic non-degenerate. But we established that $\Lambda_{P_{\mathbbm{1}}} \neq 0$. So for $\Delta^{P_{\mathbbm{1}}}_0 = \zeta \Lambda_{P_{\mathbbm{1}}}$ to hold, we need $\zeta = 0$, which negates chromatic compactness as well.
\end{proof}
\noindent In particular, this holds for all non-factorizable examples of Nenciu type and their biproducts with $u_q \mathfrak{sl}_2$.\\
However, these arguments do not hold if we drop the strong non-factorizability condition, as exhibited by another class of examples.
\begin{proposition}
\label{UHGStabProp}
    Let $H=\operatorname{SF}_{2n}$ and $U=u_q \mathfrak{sl}_2$ carrying the structures from the Proposition \ref{SFextProp}, such that all entries in $\pmb{\alpha}$ are non-zero. Then $U\ltimes H$ is chromatic compact. 
\end{proposition}
\begin{proof}
    From \cite[Lemma B.3]{beliakova_derenzi_2023} for the order of $q$, $r \equiv 0 \mod 4$, which is our usual case of interest, we have
    \[
        M_U = \frac{1}{r''}\sum^{r'-1}_{a, b=0} \sum^{r''-1}_{c, d = 0} \frac{\{1\}^{a+b}}{[a]![b]!} q^{\frac{a(a-1)+b(b-1)}{2}+2b^2+4dc} K^{-b+2c} E^a F^b  \otimes K^{b+2d} E^b F^a.
    \]
    Thus, we can deduce for $M$, the monodromy matrix of $U\ltimes H$ %not sure about the constants
    \[
        \lambda(S(M'))M'' = \frac{1}{\sqrt{r''}} \frac{\{1\}^{r'-1}}{[r'-1]!} 2^{2n}\prod^{n}_{l} \left(\alpha^2_l\right) \sum^{r''-1}_{d = 0} q^{r''-4d} K^{2d-1} E^{r'-1} F^{r'-1} \prod^{n}_{l} \left( Z^+_l Z^-_l \right).
    \]
    Now this element contains the top nilpotent element in the monomial basis, so the action on $P_{\mathbbm{1}}$ sends the generating idempotent $\Omega_{0,1,\mb{0}}$ to the socle of the module, with a multiplicative non-zero constant in front. The constant corresponds to the part of the expression corresponding to the piece of the R-matrix generated by the grouplike elements. The constant is non-zero as long as all $\alpha_l$ are, which we assumed (for the non-vanishing of the $U$-piece see \cite[Lemma 9.2]{beliakova_derenzi_2023}). \\ %Since $r$ and $p$ are coprime, $p$ is also odd, the characteristic of the field does not affect the result.
    If we compare it with the action cointegral from the Proposition \ref{SFextProp}, we clearly get an analogous situation - the top nilpotent element send the generator to the socle and the trailing grouplikes provide a non-zero constant in front, by definition of the integral. Hence $U\ltimes H$ is chromatic compact, and so chromatic non-degenerate.
\end{proof}

\begin{proof}
    This follows by applying $\lambda$ to $\lambda(S(M'))M''$, which is non-zero. 
\end{proof}
Proposition \ref{UHGStabProp} and Corollary \ref{UHGStabCor} show that the two types of strongly non-factorizable Hopf algebras described in previous sections will not produce 4d TQFTs capable of distinguishing exotic structures, regardless of the underlying field. On the other hand the weakly non-factorizable example remains susceptible to stabilisation. Moreover, the discussion in the introduction of CGHP suggests that invariants corresponding to categories with semisimple M\"uger centre are not expected to be powerful. However, we can remedy some these problems (bar stabilisation) by introducing another factor to the biproduct and changing the underlying field to one of finite characteristic.

\subsection{Unimodular ribbon Hopf algebras over fields of finite characteristic}

Fix a field $\K$ of prime characteristic $\operatorname{char} \K = p>2$. We start with showing that the constructions we used so far are not affected by this change. For the examples in this paper, we verified these statements directly using Mathematica.

\begin{proposition}
\label{nenciu_finitefield}
    Let $H = H(\mb{m}, t, \mb{d}, \mb{u})$ be a defined as in Section 2, but with the assumption that for all $a=1, ..., s$ the order of the root of unity $\xi_a$ is coprime to $p$. Let $\K$ be a field of finite characteristic containing all $\xi_a$. Then all results regarding $H$, including the unimodular and ribbon structures hold verbatim as in \cite{faes2025non}.
\end{proposition}
\begin{proof}
    This follows as all coefficients appearing in the computations involving $H$ are given in terms of $\pmb{\xi}$, that is a family of roots of unity of even orders. Some care is required when computing exponential terms appearing in the definitions of quasitriangular and ribbon structures, as exponential expressions and normalisation factors appear. For this recall the expontentials are formal and indicate the enumeration of all the products of the summands in the exponent. In case of the monodromy matrix, the integer $2$ can appear but we assumed $p>2$. Finally, since we assumed $p$ and $r'$ coprime, the normalisation factors in the R-matrix also remain unaltered. Thus, no change is incurred by this change of the field.
\end{proof}
\noindent For Examples \ref{auxExample} and \ref{ex2} this was verified directly using Mathematica. \\
We also have the analogous result for the small quantum group.
\begin{proposition}
\label{biproduct_finitefield}
     Let $U = u_q \mathfrak{sl}_2$ be the small quantum group at an even root of unity coprime to $p$. Let $\K$ be a field of finite characteristic containing $q$. Then all results regarding $U$, including the unimodular and ribbon structures hold verbatim.
\end{proposition}
\begin{proof}
    % It is suggested in \cite[Chapter 36]{lusztig2010introduction}, that the definition of $u_q \mathfrak{sl}_2$ is unaffected by passing to the field of finite characteristic. 
    To verify this, notice that all the operations in the Hopf algebra are given by linear combinations of the elements of the PBW basis and its tensor powers with coefficients containing only quantum integers defined in terms of $q$, whose order $r'$ is coprime to $p$. Thus, the operations and relations between them are unaffected by passing to the field of characteristic $p$. We similarly verify the unimodularity and ribbon axioms.
\end{proof}
\noindent For the case $r=8$, the statement was verified directly. It should be noted that we cannot drop the coprimeness assumption - otherwise, for instance, the R-matrices would no longer be invertible.

Let now $S_{p}$ be the symmetric group on $p$ letters, and let $\K[S_{p}]$ be its group algebra. Then it follows from Maschke's theorem (see for instance \cite{benson_2008}, \cite{drozd_kirichenko_2012}) that $\K[S_p]\operatorname{-mod}$ is a non-semisimple category. In complete analogy to \cite[Definition 6.6 and Proposition 6.7]{costantino2023skein} we have the following
\begin{proposition}
\label{CGHPgroupAlg}
    The symmetric category $\K[S_p]\operatorname{-mod}$ is chromatic non-degenerate but not chromatic compact.
\end{proposition}
\noindent We can improve on this result using the following strategy. 
\begin{definition}
    Let $S_{p}$ be as above and define the following action on $\operatorname{SF}_{2p}$
    \begin{align*}
        \diamond: S_{p} \times \operatorname{SF}_{2p} \rightarrow \operatorname{SF}_{2p}, &&
                    \sigma \times L \mapsto L, && 
                    \sigma \times Z^\pm_l \mapsto Z^\pm_{\sigma(l)} ,
    \end{align*}
    where $l=1, ..., p$.
\end{definition}
We now extend this action to the Hopf-algebraic one of $\K[S_p]\operatorname{-mod}$ on $u_q \mathfrak{sl}_2 \ltimes \operatorname{SF}_{2p}$. For convenience set $H = \operatorname{SF}_{2p}$, $U = u_q \mathfrak{sl}_2$ and $G = \K[S_{p}]$, where the order of $q$ is coprime to $p$.
\begin{definition}
    Let $\diamond: G \otimes U\ltimes H \rightarrow U\ltimes H$, be the left action given on generators by
    \begin{align*}
        \sigma \otimes Z^\pm_l &\mapsto Z^\pm_{\sigma(l)} && \sigma \otimes L \mapsto L && l = 1, ..., p\\
        \sigma \otimes K &\mapsto K && \sigma \otimes E \mapsto E  && \sigma \otimes F \mapsto F.
    \end{align*}
\end{definition}
Hence, we can define the Hopf algebra $U \ltimes H \rtimes G$, where "$\rtimes$" indicates the usual semi-direct product of Hopf algebras, given by the action $\diamond$, so that we introduce only the relations $\sigma Z_l  = Z_{\sigma(l)}\sigma$.

\begin{definition}
\label{UHGdef}
    Let $U \ltimes H \rtimes G$ be the Hopf algebra over $\K$ generated by $\1, E, F, K, L, Z_l, \sigma_t$, for $l=1, ..., p$ and $t=2, ..., p!$. The relations between $E, F, K, L, Z_l$ are inherited from $U\ltimes H$ and between $\sigma_t$ from $\K[S_{p}]$, with $\sigma_1\equiv \1$ omitted. With the usual association of $\1_U \otimes \1_H \otimes \sigma_t \equiv \sigma_t$, the new relations are
    \begin{align*}
         \sigma_t Z^\pm_l = Z^\pm_{\sigma_t(l)} \sigma_t && \sigma_t h = h \sigma_t,  
    \end{align*}
    for any other $h \in \{E, F, K, L\}$. The Hopf structure for $\sigma_t$ is
    \begin{align*}
        \epsilon(\sigma_t) = 1, && \Delta(\sigma_t) = \sigma_t \otimes \sigma_t, && S(\sigma_t) = \sigma^{-1}_t, 
    \end{align*}
    while for other generators it is inherited from $U\ltimes H$. 
\end{definition}
\begin{remark}
    Clearly our choice of $\sigma_t$ is not a minimal generating set for $\K[S_{p}]$, but gives a simpler description for our purposes.
\end{remark} 
\noindent It is also easy to prove the following.
\begin{proposition}
    \label{UHGBasis}
The Hopf algebra $U \ltimes H \rtimes G$ has a monomial basis
    \[
       \{E^e F^f K^k \mb{K}^\mb{w}\mb{X}^\mb{r} \sigma_t|e, f, k =1, \dots, r', \mb{w} \in \Z_\mb{m}, \mb{r} \in \Z^t_2, t=1, ..., p! \}.
    \]
\end{proposition}
\begin{proof}
    This set is clearly generating. Freeness follows from Proposition \ref{extMonomialBasis} and the choice of action $\diamond$. Namely, the action is free, and while the resulting relations with the $H$ piece are not diagonal, the result of commuting a $\sigma_t$ through a monomial remains a monomial. Thus, after fixing the position of $\sigma_t$ to be the rightmost one, we get the free property for the basis. 
\end{proof}

\begin{proposition}
    \label{UHGProp}
    The algebra $U \ltimes H \rtimes G$ is 
    \begin{enumerate}
        \item unimodular with a two-sided cointegral
        \[
            \Lambda:=\frac{\{1\}^{r'-1}}{\sqrt{r''}[r'-1]!}\sum^{r'-1}_{a=0}\sum^1_{b=0}  \sum^{p!}_{t=1} E^{r'-1} F^{r'-1}K^a L^b \prod^{p}_{l=1} Z_l^+ Z_l^- \sigma_t,
        \]
        and a left integral expressed on the monomial basis by 
        \[
            \lambda(  E^e F^f K^a L^b (Z_l^+)^g (Z_l^-)^h \sigma_t) :=
            \\ \frac{\sqrt{r''}[r'-1]!}{\{1\}^{r'-1}} \delta_{a, r'-1}\delta_{b, 0}\delta_{e, r'-1}\delta_{f, r'-1}\delta_{g, 1}\delta_{h, 1} \delta_{t, 1},
        \]
        \item quasitriangular, with the R-matrix $R:=R_{\matfont{z}}D\Theta \Bar{R}_{\pmb{\alpha}}$, where $D, \Theta$ were defined in Proposition \ref{uqsl2prop},  $R_{\matfont{z}}$ in Proposition \ref{SFprop}, and
        \[
            \Bar{R}_{\pmb{\alpha}} := \exp\left( \sum^{p}_{l=1}\alpha_l( Z_l^+ \otimes \Bar{L} Z_l^- - Z_l^- \otimes \Bar{L} Z_l^+)\right)
        \]
        for $\Bar{L} = K^{r''}L$,
        \item ribbon, with the ribbon element
        \[
            v := \frac{1-i}{\sqrt{r'}} \sum^{r'-1}_{a=0}\sum^{r''-1}_{b=0} \frac{\{-1\}^a}{[a]!} q^{-\frac{(a+3)a}{2} + 2b^2}  E^aF^a K^{-a-2b-1} L\exp \left( -2 \sum^{p}_{l=1}\alpha_l Z_l^+ Z_l^-  \right)
        \]
        corresponding to the pivotal element $g=K$,
        \item twist degenerate, that is 
        $$\lambda(v)=0.$$
    \end{enumerate}
\end{proposition}
\begin{proof}
    Most of the result follows from the Theorem \ref{extMainThm} for $U\ltimes H$ . Some additional checks involve the following:
    \begin{enumerate}
        \item It is easy to see that $\Lambda \sigma_t = \epsilon(\sigma_t)\Lambda = \Lambda$, as the rightmost sum is just the cointegral of $\K[S_p]$. For two-sidedness, we have that first
        \[
            S(\Lambda) = \frac{\{1\}^{r'-1}}{\sqrt{r''}[r'-1]!}\sum^{r'-1}_{a=0}\sum^1_{b=0}  \sum^{p!}_{t=1} \sigma_t E^{r'-1} F^{r'-1}K^a L^b \prod^{p}_{l=1} Z_l^+ Z_l^- 
        \]
        by two-sidedness of $\Lambda_{U\ltimes H}$ and the fact that all elements of $\K[S_p]$ are contained in the sum and become permuted when the antipode is taken. Then 
         \[
            S(\Lambda) = \frac{\{1\}^{r'-1}}{\sqrt{r''}[r'-1]!}\sum^{r'-1}_{a=0}\sum^1_{b=0}  \sum^{p!}_{t=1} E^{r'-1} F^{r'-1}K^a L^b \prod^{p}_{l=1} Z_l^+ Z_l^- \sigma_t, 
        \]
        follows from the fact that $\prod^{p}_{l=1} Z_l^+ Z_l^-$ commutes with all $\sigma_t$, because it is the pairs $Z_l^+ Z_l^-$ that are permuted rather than individual generators $Z^\pm_l$, which results in no additional signs. Obviously, $\sigma_t$ commute with all the remaining terms by construction.
        \item Since the form of the R-matrix is unchanged, it is left to show that (QT5) holds for $\sigma_t$. But we have 
        \begin{align*}
             R_{\matfont{z}}D\Theta \exp&\left( \sum^{p}_{l=1}\alpha_l( Z_l^+ \otimes \Bar{L} Z_l^- - Z_l^- \otimes \Bar{L} Z_l^+)\right) \Delta(\sigma_t)= \\
             &=\exp\left( \sum^{p}_{l=1}\alpha_l( Z_l^+ \otimes \Bar{L} Z_l^- - Z_l^- \otimes \Bar{L} Z_l^+)\right) (\sigma_t \otimes \sigma_t)  R_{\matfont{z}}D\Theta\\
             &=(\sigma_t \otimes \sigma_t) \exp\left( \sum^{p}_{l=1}\alpha_{\sigma_t(l)}( Z_{\sigma_t(l)}^+ \otimes \Bar{L} Z_{\sigma_t(l)}^- - Z_{\sigma_t(l)}^- \otimes \Bar{L} Z_{\sigma_t(l)}^+)\right) R_{\matfont{z}}D\Theta\\
             &=\Delta^{cop}(\sigma_t) \exp\left( \sum^{p}_{l=1}\alpha_{l}( Z_{l}^+ \otimes \Bar{L} Z_{l}^- - Z_{l}^- \otimes \Bar{L} Z_{l}^+)\right) R_{\matfont{z}}D\Theta,
        \end{align*}
        where only the sum of (commuting) elements in the exponent is permuted with no overall effect.
        \item The fact that $v$ remains central, that is commutes with all $\sigma_t$, follows by the same argument as in Item 2.
        \item This follows from twist degeneracy of $u_q \mathfrak{sl}_2$ (see for instance \cite{beliakova_derenzi_2023}).
    \end{enumerate}
\end{proof}
In what follows, we always assume that the tuple $\pmb{\alpha}$ has only non-zero entries. Having defined it, we discuss some properties of $U \ltimes H \rtimes G$. First of all, we note the following fact.

\begin{proposition}
\label{chromNon-Degen}
    The category $U \ltimes H \rtimes G\operatorname{-mod}$ is chromatic non-degenerate but not chromatic compact.
\end{proposition}
\begin{proof}
    Since the R-matrix is the same as for $U\ltimes H$, chromatic non-degeneracy follows from Proposition \ref{UHGStabProp}. For the second part, consider the action of $\Lambda$ on the projective cover of the unit $P_{\mathbbm{1}} \in U \ltimes H \rtimes G\operatorname{-mod}$. Let $e \in U \ltimes H \rtimes G$ be the corresponding idempotent (which always exists for a finite dimensional $\K$-algebra). Clearly $\sigma_t e = \epsilon(\sigma_t)e = e$. But 
    \[
        \left(\sum^{p!}_{t=1} \sigma_t \right) e = p!\,e = 0.
    \]
    Thus, $\Lambda_{P_{\mathbbm{1}}} = 0$ and  $U \ltimes H \rtimes G\operatorname{-mod}$ is not chromatic compact.
\end{proof}
\begin{proposition}
\label{UHGmuger}
    There is an embedding of symmetric categories $I: \K[S_p]\operatorname{-mod} \rightarrow \mathcal{Z}_{(2)}(U \ltimes H \rtimes G)\operatorname{-mod}$. In particular, $\mathcal{Z}_{(2)}(U \ltimes H \rtimes G)\operatorname{-mod}$ is a non-semisimple symmetric category. 
\end{proposition}
\begin{proof}
    The fully faithful braided monoidal functor $I$ is induced by the projection of the semi-direct product factor $U \ltimes H \rtimes G \twoheadrightarrow G$. Clearly any representation of $G$ can be made one of $U \ltimes H \rtimes G$ by trivially extending the action. Since there is no contribution of any $\sigma_t$ in the R-matrix, this subcategory category remains symmetric.\\
    Non-semisimplicity is a consequence of the Machke's theorem.
\end{proof}
\begin{remark}
    Proposition \ref{UHGmuger} is an improvement on both Remark \ref{SF-notSNF} and Proposition \ref{CGHPgroupAlg}. The R-matrix of $U\ltimes H$ is not symmetric but the algebra is only slightly factorizable. On the other hand the R-matrix is only symmetric (trivial) for $G$ alone. Both components are twist degenerate, and so is the result.
\end{remark}

\begin{remark}
\label{StabRem}
    It is explained in CGHP that for the TQFT defined there to detect exotic pairs of 4-manifolds, $\mathcal{C}$ has to be chromatic non-degenerate but not chromatic compact. Moreover, it is required that $\lambda(M')\lambda(S(M''))$ be non-invertible in $\K$. The latter condition is a consequence of the fact that $\lambda(M')\lambda(S(M''))$ is the invariant of $S^2\times S^2$ and by the famous theorem of Gompf (see for example \cite{gompf_stipsicz_1999}) for any pair of homeomorphic compact oriented 4-manifolds $X$ and $Y$, there exists an integer $l \in \Z_{>0}$ such that $X \#^l S^2\times S^2$ and $Y \#^l S^2\times S^2$ are diffeomorphic, and so cannot be distinguished if this invariant is invertible. This phenomenon is usually known as $ S^2\times S^2$-\textit{stabilization} or just \textit{stabilization}. For the analogue for 4-dimensional 2-handlebodies up to 2-deformations see \cite[Appendix C]{beliakova_derenzi_2023}
\end{remark}
\begin{proposition}
\label{UHGStabCor}
    Let $H=\operatorname{SF}_{2n}$ and $U=u_q \mathfrak{sl}_2$ carrying the structures from the Proposition \ref{SFextProp}, such that all entries in $\pmb{\alpha}$ are non-zero. Then $\lambda(S(M'))\lambda(M'')\neq 0$, so the resulting invariant is susceptible to stabilisation. 
\end{proposition}

\begin{proposition}
\label{StabProblem}
    Let $U\ltimes H \rtimes G$ be as above and let $M$ be its monodromy matrix. Then $\lambda(S(M'))\lambda(M'')\neq 0$. Thus, the resulting invariant is susceptible to stablilisation.
\end{proposition}
\begin{proof}
    The component $G$ does not participate non-trivially in the R-matrix (hence the monodromy) and does not affect the element of the monomial basis picked up by $\lambda$. Thus, the proof follows by the same argument as in Proposition \ref{UHGStabProp} and Corollary \ref{UHGStabCor} 
\end{proof}
\begin{remark}
    One could adapt a similar strategy for $u_q \mathfrak{g} \rtimes S_p$, where $u_q(\mathfrak{g})$ is the small quantum group for a Lie algebra $\mathfrak{g}$ larger than $\mathfrak{sl}_2$, taken over $\K$ with suitable assumptions on $p$ and $q$. Then the action of $S_p$ could be defined to permute $E_i$ generators among themselves, and analogously $F_i$. From the Hopf algebraic perspective we adapt in this paper, provided the unimodular and ribbon properties remain fulfilled, it is reasonable to expect the resulting product will also be chromatic non-degenerate but not chromatic compact, but will remain susceptible to stabilisation.
\end{remark}

\bibliography{SNFRepTh}

\Addresses

\end{document}